\documentclass{amsart}

\usepackage{amsthm}
\usepackage{amsfonts}
\usepackage{amsmath}
\usepackage{amssymb}
\usepackage{mathrsfs}
\usepackage{fullpage}
\usepackage{stmaryrd}
\usepackage{hyperref}
\usepackage{verbatim}
\usepackage{graphicx}
\usepackage{tikz}

\theoremstyle{plain}

\theoremstyle{definition}

\theoremstyle{remark}

\newcommand{\Begin}[2]{\begin{#1}\label{#2}}

\newcommand{\bPi}{\mathbf{\Pi}}
\newcommand{\bSigma}{\mathbf{\Sigma}}
\newcommand{\bDelta}{\mathbf{\Delta}}

\newcommand{\bbP}{\mathbb{P}}
\newcommand{\bbQ}{\mathbb{Q}}
\newcommand{\bbR}{\mathbb{R}}

\newcommand{\bbO}{\mathbb{O}}

\newcommand{\CM}{\mathcal{M}}

\newcommand{\SCRL}{\mathscr{L}}

\newcommand{\forces}{\Vdash}
\newcommand{\analytic}{{\bSigma_1^1}}

\newcommand{\coanalytic}{{\bPi_1^1}}

\newcommand{\borel}{{\bDelta_1^1}}

\newcommand{\cantorspace}{{{}^\omega 2}}
\newcommand{\bairespace}{{{}^\omega\omega}}
\newcommand{\finBinarySequence}{{{}^{<\omega}2}}
\newcommand{\finNaturalSequence}{{{}^{<\omega}\omega}}

\newcommand{\ZFC}{\mathsf{ZFC}}
\newcommand{\ZF}{\mathsf{ZF}}

\newcommand{\reals}{\bbR}
\newcommand{\Los}{\L{}o\'s}
\newcommand{\OD}{\mathrm{OD}}
\newcommand{\HOD}{\mathrm{HOD}}
\newcommand{\AD}{\mathsf{AD}}
\newcommand{\AC}{\mathsf{AC}}
\newcommand{\DC}{\mathsf{DC}}
\newcommand{\ON}{\mathrm{ON}}
\newcommand{\dom}{\mathrm{dom}}
\newcommand{\rang}{\mathrm{rang}}
\newcommand{\powerset}[1]{{\mathscr{P}(#1)}}

\newcommand{\WO}{{\mathrm{WO}}}

\newcommand{\decode}{\mathsf{decode}}
\newcommand{\GC}{\mathsf{GC}}
\newcommand{\fail}{\mathsf{fail}}
\newcommand{\block}{\mathsf{block}}
\newcommand{\joint}{\mathsf{joint}}
\newcommand{\ot}{\mathrm{ot}}
\newcommand{\cut}{\mathsf{cut}}
\newcommand{\enum}{\mathsf{enum}}
\newcommand{\bbA}{\mathbb{A}}
\newcommand{\degrees}{\mathcal{D}}
\newcommand{\add}{\mathsf{add}}
\newcommand{\funct}{\mathsf{funct}}
\newcommand{\symreals}{{\dot\reals_\mathrm{sym}}}

\begin{document}

\title{An Introduction to Combinatorics of Determinacy}

\author{William Chan}
\address{Department of Mathematics, University of North Texas, Denton, TX 76203}
\email{William.Chan@unt.edu}

\thanks{May 11, 2020. The author was supported by NSF grant DMS-1703708.}

\maketitle


\section{Introduction}\label{introduction}

This article is an introduction to combinatorics under the axiom of determinacy. The main topics are partition properties and $\infty$-Borel codes. To illustrate the important ideas, the article will focus on the simplest settings. This will mean that one will work at simple pointclass such as $\analytic$, small cardinals such as $\omega_1$ or $\omega_2$, or very natural models of determinacy such as $L(\reals)$.

Despite the simplicity, there are still many interesting combinatorial questions in these settings. One purpose of this article is to serve as a reference for some of the notations and backgrounds that will be helpful for some results of Jackson, Trang, and the author on combinatorics around $\omega_1$, $\omega_2$, and $\reals$ under $\AD$ and $\AD^+$. (Section 2 to 6 were written for \cite{Definable-Combinatorics-at-First-Uncountable-Cardinal} and \cite{More-Definable-Combinatorics-Around-First-Second-Uncountable-Cardinal}. Section 7 and 8 will be helpful for reading \cite{Applications-Infinity-Borel-Codes-Definability-Definable-Cardinals} and \cite{L(R)-With-Determinacy-Satisfies-Suslin-Hypothesis}.)

The article should be accessible with very basic knowledge of descriptive set theory and general set theory. Some familiarity with determinacy and games, the pointclass $\analytic$, and the bounding principle for $\analytic$ should suffice for most of the article. Knowledge of club sets, ultrafilters, measures, and basic constructibility will also be helpful. There will be occasional mentions of the theory of prewellorderings and scales such as the Moschovakis coding lemma or the Kunen-Martin theorem. In a few places, one will refer to the result of Kechris that $\AD$ implies $L(\reals) \models \DC_\reals$. The last two sections will require more familiarity with topics in general set theory such as $\HOD$ and iterated forcing.

The article has essentially two parts. The first part deals with various partition properties which are very powerful tools of $\AD$. The second part introduces the $\infty$-Borel code and the Vop\'enka forcing which are useful $\AD^+$ tool. Except possibly for a few minor observations, no results in this article are due to the author. 

Section 2 introduces the partition property in its various forms. It also defines the associated measures induced by the partition properties. 

Section 3 develops the notation and theory around good codings of functions. This section proves Martin's criterion for establishing partition properties. 

Section 4 gives some examples of good coding of functions on $\omega_1$. This is used to give a simple proof of the weak partition property on $\omega_1$ and two proofs of the strong partition property on $\omega_1$. One will give Martin's original proof of the strong partition property using sharps. Then one gives a proof due to Kechris of the strong partition property on $\omega_1$ which uses category arguments and the simple generic coding function for $\omega_1$.

Section 5 develops the theory around the Kunen function which is named after the eponymous Kunen tree. Here it will be shown that every function from $\omega_1$ to $\omega_1$ has a Kunen function by using the original Kunen tree. The Kunen functions are especially useful for choosing representatives of certain ultrapowers. This will be used to establish the identity of $\omega_2$ as the ultrapower ${}^{\omega_1}\omega_1 \slash \mu$, where $\mu$ is the club ultrafilter on $\omega_1$. Then one will prove a useful sliding lemma and use this to establish the weak partition property on $\omega_2$. 

Section 6 is dedicated to proving a result of Martin and Paris which shows that $\omega_2$ is not a strong partition cardinal. However, this section presents a proof of Jackson which produces an explicit partition without a homogeneous subset. 

Section 7 proves some results of Woodin about the structure of $L(\reals)$ and the nature of $\infty$-Borel codes in $L(\reals)$. In particular, it will be shown that in $L(\reals)$, there is an ``ultimate'' $\infty$-Borel code that can be used to generate all other $\infty$-Borel codes for sets of reals. This is quite useful for cardinality results in $L(\reals)$. 

Section 8 will present some further descriptive set theoretic applications due to Woodin of the Vop\'enka forcing. The countable section uniformization and the Woodin's perfect set dichotomy theorem will be proved here.

The author would like to thank Stephen Jackson for the numerous conversations concerning the topics that appear here. The author would also like to thank Thilo Wienert for carefully reading this article and suggesting many comments. 

\section{Partition Properties}\label{partition properties}

The following is the usual partition property:

\Begin{definition}{usual partition property}
Let $\kappa$ be an ordinal, $\lambda \leq \kappa$, and $\gamma < \kappa$, then $\kappa \rightarrow (\kappa)^\lambda_\gamma$ indicates that for every $P : [\kappa]^\lambda \rightarrow \gamma$, there is some $\beta < \gamma$ and some $A \subseteq \kappa$ with $|A| = \kappa$ so that for all $f \in [A]^\lambda$, $P(f) = \beta$. 

The most frequent situation is when $\gamma = 2$.

$\kappa$ is a strong partition cardinal if and only if $\kappa \rightarrow (\kappa)^\kappa_2$. 

$\kappa$ is a weak partition cardinal if and only if $\kappa \rightarrow (\kappa)^\alpha_2$ for all $\alpha < \kappa$. 
\end{definition}

\Begin{definition}{more general uniform cofinality}
Let $\kappa$ be an ordinal. Let $S : \kappa \rightarrow \kappa$. Let $Y^S = \{(\alpha,\beta) : \alpha < \kappa \wedge \beta < S(\alpha)\}$. 

Let $f : \kappa \rightarrow \kappa$. $f$ is said to have uniform cofinality $S$ if and only if there is a function $g : Y^S \rightarrow \kappa$ with the following properties 

\noindent (1) For all $\alpha < \kappa$ and all $\beta < \gamma < S(\alpha)$, $g(\alpha,\beta) < g(\alpha,\gamma)$ 

\noindent (1) For all $\alpha < \kappa$, $f(\alpha) = \sup\{g(\alpha,\beta) : \beta < S(\alpha)\}$. 

Let $\mu$ be a measure on $\kappa$. $f : \kappa \rightarrow \kappa$ is said to have uniform cofinality $S$ $\mu$-almost everywhere if and only if there is a $g : Y^S \rightarrow \kappa$ as above so that for $\mu$-almost all $\alpha < \kappa$, $f(\alpha) = \sup\{g(\alpha,\beta) : \beta < S(\alpha)\}$.

A function $f : \kappa \rightarrow \kappa$ has uniform cofinality $\omega$ if and only if $f$ has uniform cofinality $S$ where $S(\alpha) = \omega$ for all $\alpha < \kappa$. 
\end{definition}

\Begin{definition}{correct type}
Let $f : \kappa \rightarrow \kappa$ be a function. $f$ is discontinuous everywhere if and only if for all $\alpha < \kappa$, $f(\alpha) > \sup\{f(\beta) : \beta < \alpha\}$. 

(Jackson) A function $f : \kappa \rightarrow \kappa$ has the correct type if and only if $f$ has uniform cofinality $\omega$ and is discontinuous everywhere.

Let $[\kappa]^\kappa_*$ denote the subset of $[\kappa]^{\kappa}$ consisting of functions of the correct type. 
\end{definition}

There are partition properties formulated for functions of the correct type in which the homogeneous set can be chosen to be club. Club homogeneous sets are conceptually useful in various construction. Functions of the correct type seem to naturally appear in many proofs establishing the partition properties under $\AD$. 

It should be noted that throughout the survey an asterisk, $*$, will be used to denote the corresponding concept that involves functions of the correct type when there is an ordinary version of this concept.

\Begin{definition}{partition properties}
Let $\kappa$ be an ordinal. Let $\lambda \leq \kappa$ and $\gamma< \kappa$. Let $\kappa \rightarrow_* (\kappa)^\lambda_\gamma$ assert that for all $P : [\kappa]^{\lambda}_* \rightarrow \gamma$, there exists a club $C \subseteq \kappa$ and a $\beta < \lambda$ so that for all $f \in [C]^\lambda_*$, $P(f) = \beta$. Such a club $C$ is said to be homogeneous for $P$ taking value $\beta$ (for functions of the correct type). 
\end{definition}

\Begin{definition}{enumeration of sets}
Let $\kappa$ be an ordinal. If $A \subseteq \kappa$ is such that $|A| = \kappa$, then let $\enum_A : \kappa \rightarrow A$ denote the increasing enumeration of $A$. (In context, when one writes $\enum_A$, it should be clear what $\kappa$ is.)
\end{definition}

\Begin{fact}{equivalence of partition properties}
Let $\kappa$ be an ordinal and $\lambda \leq \kappa$. Then $\kappa \rightarrow_* (\kappa)^\lambda_2$ implies $\kappa \rightarrow (\kappa)^\lambda_2$. $\kappa \rightarrow (\kappa)^{\omega\cdot\lambda}_2$ implies $\kappa \rightarrow_* (\kappa)^\lambda_2$.
\end{fact}

\begin{proof}
Assume $\kappa \rightarrow_* (\kappa)^\lambda_2$. 

Let $P : [\kappa]^\lambda \rightarrow 2$. Then there is some $C \subseteq \omega_1$ club so that $C$ is homogeneous for $P$ for function of the correct type. Let $D = \{\alpha < \kappa : (\exists \beta \in \mathrm{Lim})(\alpha = \enum_C(\beta + \omega))\}$, where $\mathrm{Lim}$ refers to the class of limit ordinals. Note $|D| = \kappa$. 

The next claim is that every $f \in [D]^\lambda$ is a function of the correct type. Let $g : \kappa \times \omega \rightarrow \kappa$ be defined as follows: Suppose $\gamma < \kappa$. Let $\beta_\gamma$ be the unique limit ordinal $\beta$ so that $f(\gamma) = \enum_C(\beta_\gamma + \omega)$. Then for each $n \in \omega$, define $g(\gamma,n) = \enum_C(\beta_\gamma + n)$. Then it is clear that $f(\gamma) = \sup\{g(\gamma, n) : n \in \omega\}$. This shows that $f$ has uniform cofinality $\omega$. Suppose $\gamma < \kappa$. There is a unique limit ordinal $\beta_\gamma$ so that $f(\gamma) = C(\beta_\gamma + \omega)$. Then for all $\epsilon < \gamma$, $f(\epsilon) \leq \enum_C(\beta_\gamma) < f(\gamma)$. That is, $\sup\{f(\epsilon) : \epsilon < \gamma\} \leq \enum_C(\beta_\gamma) < f(\gamma)$. $f$ is discontinuous everywhere.

Thus every $f \in [D]^{\lambda}$ belongs to $[C]^\lambda_*$. Since $C$ is homogeneous for $P$ for functions of the correct type, $D$ is homogeneous for $P$ in the ordinary sense. This establishes $\kappa \rightarrow (\kappa)^\lambda_2$. 

Now suppose $\kappa \rightarrow (\kappa)^{\omega\cdot\lambda}_2$.

Suppose $P : [\kappa]^\lambda \rightarrow 2$. Let $\block : [\kappa]^{\omega\cdot\lambda} \rightarrow [\kappa]^\lambda$ be defined by $\block(f)(\gamma) = \sup\{\omega\cdot\gamma + n : n \in \omega\}$. Define $P' : [\kappa]^{\omega\cdot\lambda} \rightarrow 2$ by $P'(f) = P(\block(f))$. Let $A \subseteq \omega_1$ be such that $|A| = \kappa$ and $A$ is homogeneous for $P'$ in the ordinary sense. Without loss of generality, suppose $P'(h) = 0$ for all $h \in [A]^{\omega\cdot\lambda}$. Let $C$ be the collection of limit points of $A$. That is, $C = \{\alpha \in \kappa : \alpha = \sup (C \cap \alpha)\}$. $C$ is a club subset of $\omega_1$. 

Suppose $f \in [C]^{\lambda}_*$. Since $f$ is of the correct type, let $g : \lambda \times \omega \rightarrow \omega_1$ witness that it has uniform cofinality $\omega$. Let $\gamma < \kappa$. Since $f$ is discontinuous everywhere, $\sup\{f(\alpha) : \alpha < \gamma\} < f(\gamma)$. Therefore there is some $N_\gamma \in \omega$ so that for all $n \geq N_\gamma$, $g(\gamma, n) > \sup\{f(\alpha) : \alpha < \gamma\}$. Since $f(\gamma)$ is a limit point of $A$, define by recursion, $h(\omega\cdot \gamma + n)$ to be the least element of $A$ greater than 
$$\max\{g(\gamma, N_\gamma + n), h(\omega\cdot\gamma +  k) : k < n\}.$$
Then $h \in [A]^{\omega\cdot\gamma}$ and $\block(h) = f$. Thus $P(f) = P'(h) = 0$. 

It has been shown that for all $f \in [C]^{\lambda}_*$, $P(f) = 0$. $\kappa \rightarrow_* (\kappa)^\lambda_2$ has been established.
\end{proof}

Partition properties on a cardinal $\kappa$ yield other very interesting properties on $\kappa$. Much of the material of the remainder of this section can be found in \cite{Infinitary-Combinatorics-and-the-Axiom-of-Determinateness}.

\Begin{fact}{partition implies regularity regularity}
$\kappa \rightarrow (\kappa)^2_2$ implies that $\kappa$ is regular.
\end{fact}

\begin{proof}
Suppose $\eta < \kappa$ and $h : \eta \rightarrow \kappa$ is a cofinal function. Define $P : [\kappa]^2 \rightarrow 2$ by 
$$P(\alpha,\beta) = \begin{cases}
0 & \quad (\exists \gamma < \eta)(\alpha \leq h(\gamma) \leq \beta) \\
1 & \quad \text{otherwise}
\end{cases}.$$
Let $A \subseteq \omega_1$ with $|A| = \kappa$ be homogeneous for $P$. If $A$ is homogeneous for $P$ taking value $0$, then one can show that $\ot(\rang(h)) = \kappa$, which is a contradiction. If $A$ is homogeneous for $1$, then one can show that $\rang(h) \subseteq \min A$. This violates $\kappa \rightarrow (\kappa)^2_2$. 
\end{proof}

\Begin{definition}{eta club filter}
Let $\kappa$ be a regular cardinal and $\eta< \kappa$ be a limit ordinal. A set $C \subseteq \kappa$ is $\eta$-closed if and only if for all $f : \eta \rightarrow C$ increasing, $\sup(f) \in C$. $C \subseteq \kappa$ is a $\eta$-club if and only if $C$ is $\eta$-closed and unbounded.

Let $W^\kappa_\eta$ denote the filter of sets containing an $\eta$-club as a subset.
\end{definition}

\Begin{fact}{intersection less than kappa eta club}
Let $\kappa$ be a regular cardinal and $\eta < \kappa$ be a limit ordinal. Let $\delta < \kappa$ and $\langle C_\alpha : \alpha < \delta\rangle$ be a sequence of $\eta$-clubs. Then $\bigcap_{\alpha \in \delta} C_\alpha$ is an $\eta$-club.
\end{fact}

\begin{proof}
Clearly $\bigcap_{\gamma< \delta} C_\gamma$ is $\eta$-closed. One needs to show that $\bigcap_{\gamma<\delta} C_\gamma$ is unbounded.

Fix $\epsilon < \kappa$. Let $\beta_0^0$ be the least element of $C_0$ greater than $\epsilon$. 

Suppose for $\alpha < \eta$ and $\gamma < \delta$, $\beta_\alpha^\xi$ has been defined for all $\xi < \gamma$. Note $\sup\{\beta^\xi_\alpha : \xi < \gamma\} < \kappa$ by the regularity of $\kappa$. Let $\beta_\alpha^\gamma$ be the least element of $C_\gamma$ which is larger than $\sup\{\beta^\xi_\alpha : \xi < \gamma\}$. 

Suppose for $\alpha < \eta$, $\beta^{\gamma}_\nu$ has been defined for all $\nu < \alpha$, and $\gamma < \delta$.  By the regularity of $\kappa$, $\sup\{\beta^\gamma_\nu : \nu < \alpha \wedge \gamma < \delta\} < \kappa$. Let $\beta_\alpha^0$ be the least element of $C_0$ greater than $\sup\{\beta^\gamma_\nu : \nu < \alpha \wedge \gamma < \delta\}$. 

Note that for all $\gamma_0,\gamma_1 \in \delta$, $\sup\{\beta_\alpha^{\gamma_0} : \alpha < \eta\} = \sup\{\beta_\alpha^{\gamma_1} : \alpha < \eta\}$. Let $\lambda$ denote this common value. $\lambda \in C_\gamma$ for all $\gamma < \delta$ since $C_\gamma$ is an $\eta$-club. Thus $\lambda \in \bigcap_{\gamma < \delta} C_\gamma$ and $\lambda > \epsilon$. This shows that $\bigcap_{\gamma < \delta} C_\gamma$ is unbounded.
\end{proof}

In $\ZF$, Fact \ref{intersection less than kappa eta club} does not imply that $W^\kappa_\eta$ is $\kappa$-complete. Suppose $\langle A_\alpha : \alpha < \delta\rangle$ where $\delta < \kappa$ is a sequence in $W^\kappa_\eta$. For each $A_\alpha$, there is an $\eta$-club $C \subseteq A_\alpha$. To apply Fact \ref{intersection less than kappa eta club}, one would need to produce a sequence of $\eta$-clubs, $\langle C_\alpha : \alpha < \delta\rangle$, so that $C_\alpha \subseteq A_\alpha$ for each $\alpha < \delta$. This appears to require some choice principle.

It will be shown next that appropriate partition properties imply that $W^\kappa_\gamma$ is an ultrafilter, is $\kappa$-complete, and is normal.

\Begin{fact}{eta limit set is eta club}
Let $\kappa$ be a regular cardinal and $\eta < \kappa$ be a limit ordinal. Let $A \subseteq \kappa$ be an unbounded set. Let $\mathrm{Lim}^\eta(A) = \{\sup(f) : f \in [A]^\eta\}$. That is, $\mathrm{Lim}^\eta(A)$ is the collection of all $\alpha \in \kappa$ which are the supremum of an $\eta$-increasing sequence through $A$. 

Then $\mathrm{Lim}^\eta(A)$ is an $\eta$-club.
\end{fact}

\begin{proof}
$\mathrm{Lim}^\eta(A)$ is clearly unbounded since $A$ is unbounded.

Suppose $f : \eta \rightarrow \mathrm{Lim}^\eta(A)$ is an increasing function. For each $\alpha < \eta$, let $g(\alpha)$ be the least element $\gamma \in A$ so that $f(\alpha) < \gamma <  f(\alpha + 1)$ which exists since $f(\alpha + 1) = \sup(h)$ for some $h \in [A]^\eta$. Thus $\sup(f) = \sup(g)$ and $g \in [A]^\eta$. Thus $\sup(f) \in \mathrm{Lim}^\eta(A)$. $\mathrm{Lim}^\eta(A)$ is an $\eta$-club.
\end{proof}

\Begin{fact}{partition ultrafilter}
Let $\kappa$ be a regular cardinal and $\eta < \kappa$ be a limit ordinal. Suppose $\kappa \rightarrow (\kappa)^\eta_2$. Then $W^\kappa_\eta$ is an ultrafilter.
\end{fact}

\begin{proof}
Let $A \subseteq \kappa$. Define the partition $P : [\kappa]^\eta \rightarrow 2$ by $P(f) = 1 \Leftrightarrow \sup(f) \in A$. 

By $\kappa \rightarrow (\kappa)^\eta_2$, let $B \subseteq \kappa$ with $|B| = \kappa$ be homogeneous for this partition. Without loss of generality, suppose $B$ is homogeneous taking value $1$. $\mathrm{Lim}^\eta(B)$ is an $\eta$-club by Fact \ref{eta limit set is eta club}. Let $\alpha \in \mathrm{Lim}^\eta(B)$. Let $f \in [B]^\eta$ be such that $\sup(f) = \alpha$. Then $P(f) = 1$ implies that $\alpha = \sup(f) \in A$. This shows that $\mathrm{Lim}^\eta(B) \subseteq A$ which implies $A \in W^\kappa_\eta$. 

If $B$ was homogeneous for $P$ taking value $0$, then the same argument would have shown $\kappa \setminus A \in W^\kappa_\eta$. 
\end{proof}

\Begin{fact}{partition and completeness}
Assume $\kappa$ is a regular cardinal and $\eta < \kappa$ is a limit ordinal. Let $\lambda < \kappa$ be infinite. If $\kappa \rightarrow (\kappa)^\eta_\lambda$ holds, then $W^\kappa_\eta$ is $\lambda^+$-complete. 
\end{fact}

\begin{proof}
Let $\langle A_\alpha : \alpha < \lambda\rangle$ be a sequence in $W^\kappa_\eta$. Suppose $\bigcap_{\alpha < \lambda} A_\alpha \notin W^\kappa_\eta$. Since $W^\kappa_\eta$ is an ultrafilter by Fact \ref{partition ultrafilter}, one may assume that $\bigcap_{\alpha < \lambda}A_\alpha = \emptyset$ by adding one further set to the sequence.

Let $P : [\kappa]^\eta \rightarrow \lambda$ be defined by $P(f)$ is the least $\alpha < \lambda$ so that $\sup(f) \notin A_\alpha$. By $\kappa \rightarrow (\kappa)^\eta_\lambda$, there is some $\alpha^*$ and a $B \subseteq \kappa$ so that $|B| = \kappa$ and $P(f) = \alpha^*$ for all $f \in [B]^\eta$. 

By Fact \ref{eta limit set is eta club}, $\mathrm{Lim}^\eta(B) \in W^\kappa_\eta$. However, $\mathrm{Lim}^\eta(B) \cap A_{\alpha^*} = \emptyset$. This contradicts $A_{\alpha^*} \in W^\kappa_\eta$.
\end{proof}

\Begin{fact}{sum partition imply lower subscript partition}
Assume $\kappa$ is a regular cardinal and $\eta < \kappa$ is a limit ordinal. $\kappa \rightarrow (\kappa)^{\eta + \eta}_2$ implies $\kappa \rightarrow (\kappa)^\eta_\lambda$ for all $\lambda < \kappa$. Therefore, $\kappa \rightarrow (\kappa)^{\eta + \eta}_2$ implies $W^\kappa_\eta$ is a $\kappa$-complete measure on $\kappa$.

If $\kappa$ has the weak partition property, then $W^\kappa_\eta$ is a $\kappa$-complete measure on $\kappa$ for all limit ordinals $\eta < \kappa$.
\end{fact}

\begin{proof}
Let $P : [\kappa]^\eta \rightarrow \lambda$. Define $Q : [\kappa]^{\eta + \eta} \rightarrow 2$ by $Q(f_0\hat{\ }f_1) = 0$ if and only if $P(f_0) = P(f_1)$. By $\kappa \rightarrow (\kappa)^{\eta + \eta}_2$, let $A \subseteq \kappa$ with $|A| = \kappa$ be homogeneous for $Q$. 

Suppose $A$ is homogeneous for $Q$ taking value $1$. Let $Q' : [A]^{\eta + \eta} \rightarrow 2$ be defined by 
$$Q'(f_0\hat{\ }f_1) = \begin{cases}
0 & \quad P(f_0) < P(f_1) \\ 
1 & \quad P(f_0) > P(f_1)
\end{cases}$$
Again by $\kappa \rightarrow (\kappa)^{\eta + \eta}_2$, there is a $B \subseteq A$ with $|B| = \kappa$ which is homogeneous for $Q'$. One can check that if $B$ is homogeneous for $Q'$ taking value $1$, then one would have an infinite descending sequences of ordinals. If $B$ is homogeneous for $Q'$ taking value $0$, then one can produce an injection of $\kappa$ into $\lambda < \kappa$. $B$ can not be homogeneous. Contradiction. Therefore, $A$ must have been homogeneous for $Q$ taking value $0$. 

Now the claim is that $A$ is homogeneous for $P$: Let $f_0,f_1 \in [A]^\eta$. Find $f_2 \in [A]^\eta$ such that $\min(f_2)$ is larger than both $\sup(f_0)$ and $\sup(f_1)$. Then $f_0\hat{\ }f_2 \in [A]^{\eta + \eta}$ and $f_1 \hat{\ } f_2 \in [A]^{\eta + \eta}$. Since $A$ is homogeneous for $Q$ taking value $0$, $Q(f_0\hat{\ } f_2) = 0$ and $Q(f_1\hat{\ }f_2) = 0$. This implies that $P(f_0) = P(f_2) = P(f_1)$. Thus $A$ is homogeneous for $P$.
\end{proof}

\Begin{fact}{partition property and normality}
Let $\kappa$ be a regular cardinal and $\eta < \kappa$ be a limit ordinal. Suppose $\kappa \rightarrow (\kappa)^{\eta + \eta}_2$. Then $\kappa$ is a normal $\kappa$-complete ultrafilter.

If $\kappa$ has the weak partition property, then $W^\kappa_\eta$ is a normal $\kappa$-complete measure on $\kappa$ for each limit ordinal $\eta < \kappa$.
\end{fact}

\begin{proof}
By Fact \ref{sum partition imply lower subscript partition}, $W^\kappa_\eta$ is a $\kappa$-complete ultrafilter.

Let $F : \kappa \rightarrow \kappa$ be a regressive function. That is, $A' = \{\alpha \in \kappa : F(\alpha) < \alpha\} \in W^\kappa_\eta$. Let $A \subseteq A'$ be a $\eta$-club set. Define $P : [\kappa]^\eta \rightarrow 2$ by
$$P(f) = \begin{cases}
0 & \quad F(\sup(f)) < \min f \\
1 & \quad \text{otherwise}
\end{cases}$$
Since $\kappa \rightarrow (\kappa)^{\eta + \eta}_2$ implies $\kappa \rightarrow (\kappa)^\eta_2$, let $B \subseteq A$ with $|B| = \kappa$ be homogeneous for $P$.

Let $\tilde B = \{\enum_B(\eta \cdot \alpha + \eta) : \alpha < \kappa\}$, where $\enum_B : \kappa \rightarrow \kappa$ is the increasing enumeration of $B$. Suppose $f \in [\tilde B]^\eta$. Since $\tilde B \subseteq B \subseteq A$ and $A$ is $\eta$-closed, $\sup(f) \in A \subseteq A'$. Therefore, $F(\sup(f)) < \sup(f)$. Let $\gamma < \eta$ be least so that $F(\sup(f)) < f(\gamma)$. Since $f(\gamma + 1) \in \tilde B$, let $f_0 : \gamma + 1 \rightarrow B$ be such that for all $\nu < \gamma + 1$, $f(\gamma) < f_0(\nu) < f(\gamma + 1)$. Let $f_1 : (\eta - (\gamma + 1)) \rightarrow \tilde B$ be the tail of $f$ above $f(\gamma)$. Note that $f_2 = f_0\hat{\ }f_1$ is an $\eta$-sequence in $B$ with the property that $\min(f_2) > f(\gamma) > F(\sup(f))$ and $\sup(f_2) = \sup(f)$. Thus $F(\sup(f_2)) = F(\sup(f)) < \min(f_2)$. Since $f_2 \in [B]^\eta$, one has that $B$ must be homogeneous for $P$ taking value $0$.

Let $f \in [B]^\eta$. Let $f' : \eta \rightarrow B$ be the increasing enumeration of $\{\min B\} \cup \rang(f)$. Thus $P(f') = 0$ implies that $F(\sup(f)) = F(\sup(f')) < \min(f') = \min(B)$. It has been shown that for all $f \in [B]^\eta$, $F(\sup(f)) < \min(B)$. Thus for all $\alpha \in \mathrm{Lim}^\eta(B)$, $F(\alpha) < \min(B)$. Since $\mathrm{Lim}^\eta(B) \in W^\kappa_\eta$ and $W^\kappa_\eta$ is $\kappa$-complete, there is some $\xi < \min(B)$ and $C \in W^\kappa_\eta$ so that $F[C] = \xi$. It has been shown that $F$ is constant $W^\kappa_\eta$-almost everywhere. Normality has been established.
\end{proof}

Under suitable circumstances on $\kappa$, one can determine all the normal $\kappa$-complete measures on $\kappa$.

\Begin{fact}{characterize normal measure by partition}
Let $\kappa$ be a regular cardinal and $\eta < \kappa$ be a limit ordinal. Assume $\kappa \rightarrow (\kappa)^{\eta}_2$. Then $W^\kappa_\eta = W^\kappa_{\mathrm{cof}(\eta)}$. 

Suppose that $\eta_0 < \eta_1$ are two infinite regular cardinals less than $\kappa$. Then $W^\kappa_{\eta_0} \neq W^\kappa_{\eta_1}$.

Suppose the collection of infinite regular cardinals below $\kappa$ has cardinality less than $\kappa$. Let $\mu$ be any $\kappa$-complete normal measure on $\kappa$. There is some infinite regular cardinal $\eta < \kappa$ so that $\mu$ is equivalent to $W^\kappa_\eta$. 
\end{fact}

\begin{proof}
The partition property $\kappa \rightarrow (\kappa)^{\eta}_2$ implies that $W^\kappa_{\eta}$ and $W^\kappa_{\mathrm{cof}(\eta)}$ are both ultrafilters. 

Suppose $A \subseteq \kappa$ is a $\mathrm{cof}(\eta)$-club. Let $f \in [A]^\eta$. Let $\rho : \mathrm{cof}(\eta) \rightarrow \eta$ be a cofinal increasing sequence. Note that $f \circ \rho \in [A]^{\mathrm{cof}(\eta)}$. Since $A$ is a $\mathrm{cof}(\eta)$-club, $\sup(f \circ \rho) = \sup(f) \in A$. Thus $A$ is also an $\eta$-club. This shows that $W^\kappa_{\mathrm{cof}(\eta)} \subseteq W^\kappa_\eta$. 

Suppose that $\neg(W^\kappa_\eta \subseteq W^\kappa_{\mathrm{cof}(\eta)})$.  Let $A \in W^\kappa_\eta$ be such that $A \notin W^\kappa_{\mathrm{cof}(\eta)}$. Since $W^\kappa_{\mathrm{cof}(\eta)}$ is an ultrafilter, $\kappa \setminus A \in W^\kappa_{\mathrm{cof}(\eta)}$. It has already been shown that $W^\kappa_{\mathrm{cof}(\eta)} \subseteq W^\kappa_\eta$. Therefore, $\kappa \setminus A \in W^\kappa_\eta$. However, $A,\kappa\setminus A \in W^\kappa_\eta$ contradicts the fact that $W^\kappa_\eta$ is an ultrafilter. 

It has been established that $W^\kappa_\eta = W^\kappa_{\mathrm{cof}(\eta)}$.

Now suppose $\eta_0 < \eta_1$ are two regular cardinals less than $\kappa$. Let $A_i = \{\alpha < \kappa : \mathrm{cof}(\alpha) = \eta_i\}$. For $i \in 2$, $A_i$ is an $\eta_i$-club. Therefore, $A_i \in W^\kappa_{\eta_i}$. If $W^\kappa_{\eta_0} = W^\kappa_{\eta_1}$, then $A_0 \cap A_1 = \emptyset \in W^\kappa_{\eta_1}$. This contradicts $W^\kappa_{\eta_1}$ is a filter.

Now suppose the collection $K$ of infinite regular cardinals below $\kappa$ has cardinality less than $\kappa$. Let $\mu$ be a $\kappa$-complete normal measure on $\kappa$. 

Let $T$ be the set of limit ordinals below $\kappa$. Since $\kappa$ and $T$ are in bijection, $\mu$ is equivalent to a measure concentrating on $T$. Therefore, assume for simplicity that $T \in \mu$. For each $\eta \in K$, let $A_\eta = \{\alpha \in \kappa : \mathrm{cof}(\alpha) = \eta\}$. $T = \bigcup_{\eta \in K} A_\eta$. Since $\mu$ is $\kappa$-complete and $|K| < \kappa$, there is some $\eta_0$ so that $A_{\eta_0} \in \mu$. 

The claim is that $\mu = W^\kappa_{\eta_0}$: Let $B \in W^\kappa_{\eta_0}$. Let $C \subseteq B$ be an $\eta_0$-club set. Suppose $B \notin \mu$. Then $\kappa \setminus B \in \mu$. Let $Q = A_{\eta_0} \cap (\kappa \setminus B)$. Note $Q \in \mu$. 

Let $\alpha \in Q$. Note that $\mathrm{cof}(\alpha) = \eta_0$ and $\alpha \notin B$ so $\alpha \notin C$. $\sup(C \cap \alpha) < \alpha$ since otherwise $\alpha \in C$ since $\mathrm{cof}(\alpha) = \eta_0$ and $C$ is an $\eta_0$-club. Let $G : Q \rightarrow \kappa$ be defined $G(\alpha) = \sup(C \cap \alpha)$. $G$ is a regressive function on $Q$. Since $\mu$ is normal, there is some $\gamma < \kappa$ so that $J = \{\alpha : G(\alpha) = \gamma\} \in \mu$. For any $\alpha \in J$, $C \cap \alpha \subseteq \gamma$. Since $J \in \mu$ implies that $J$ is unbounded, this implies that $C \subseteq \gamma$. This is impossible since $C$ is unbounded.
\end{proof}

\section{Good Coding of Functions}\label{good coding of functions}

\Begin{definition}{good coding of functions}
Let $\kappa$ be a regular cardinal and $\lambda \leq \kappa$ be an ordinal. A good coding system for ${}^\lambda\kappa$ consists of $\Gamma$, $\decode$, and $\GC_{\beta,\gamma}$ for each $\beta < \lambda$ and $\gamma < \kappa$ with the following properties:

\noindent (1) $\Gamma$ is a pointclass closed under continuous substitution and $\exists^\reals$. Let $\check \Gamma$ denote the dual pointclass. Let $\Delta = \Gamma \cap \check\Gamma$.

\noindent (2) $\decode : \reals \rightarrow \powerset{\lambda\times\kappa}$. For all $f \in {}^\lambda\kappa$, there is some $x \in \reals$ so that $\decode(x) = f$. 

\noindent (3) For all $\beta < \lambda$ and $\gamma < \kappa$, $\GC_{\beta,\gamma} \subseteq \reals$, $\GC_{\beta,\gamma} \in \Delta$, and $\GC_{\beta,\gamma}$ is defined by $x \in \GC_{\beta,\gamma}$ if and only if
$$\decode(x)(\beta,\gamma) \wedge (\forall \gamma' < \kappa)(\decode(x)(\beta,\gamma') \Rightarrow \gamma = \gamma').$$

For each $\beta < \lambda$, let $\GC_\beta = \bigcup_{\gamma < \kappa} \GC_{\beta,\gamma}$.

\noindent (4) (Boundedness property) Suppose $A \in \exists^\reals\Delta$ and $A \subseteq \GC_\beta$, then there exists some $\delta < \kappa$ so that $A \subseteq \bigcup_{\gamma < \delta}\GC_{\beta,\gamma}$. 

\noindent (5) $\Delta$ is closed under less than $\kappa$ length wellordered unions. 
\newline

Suppose $x \in \reals$, let $\fail(x)$ be the least $\beta < \lambda$ so that $x \notin \GC_\beta$ if it exists. Otherwise, let $\fail(x) = \infty$.

Let $\GC = \bigcap_{\beta < \lambda}\GC_\beta$. Note that if $x \in \GC$, then $\decode(x)$ is the graph of a function in ${}^\lambda\kappa$. If $x \in \GC$, then one will use function notations such as $\decode(x)(\beta) = \gamma$ to indicate $(\beta,\gamma) \in \decode(x)$.
\end{definition}

Assuming $\AD$, (5) follows from the other four conditions. This comes from a pointclass argument. See the end of the proof of Theorem 2.34 of \cite{Structural-Consequences-AD}. Later, one will apply this to $\omega_1$ with the associated pointclass $\analytic$. It is clear that (5) holds in this setting since $\borel$ is closed under countable unions.

\Begin{remark}{good coding remarks}
The meaning of $x \in \GC_{\beta,\gamma}$ is that $x$ is good at $\beta$ in the sense that $\decode(x)$ has successfully mapped $\beta$ to $\gamma$. One interprets $x \in \GC_{\beta}$ to be mean that $x$ is good at $\beta$ in the sense that $\decode(x)$ has successfully mapped $\beta$ to some value. So a reals $x$ belongs to $\GC$ means that $x$ is a good code in the sense that $\decode(x)$ is truly a function from $\lambda$ to $\kappa$.
\end{remark}

\Begin{definition}{good strategies of a coding system}
Let $\kappa$ be a regular cardinal and let $\lambda \leq \kappa$ be an ordinal. Let $(\Gamma,\decode,\GC_{\beta,\gamma} : \beta < \lambda,\gamma<\kappa)$ be a coding system for ${}^\lambda\kappa$. 

Let $S^1$ consists of reals $r$ coding a Lipschitz continuous function $\Xi_r : \reals \rightarrow \reals$ which has the following property:
$$(\forall y)[\fail(y) \leq \fail(\Xi_r(y)) \wedge (\fail(y) < \infty \Rightarrow \fail(y) < \fail(\Xi_r(y)))].$$

Let $S^2$ consists of Player 2 strategies $s$ so that its associated Lipschitz continuous function $\Xi_s : \reals \rightarrow \reals$ has the property that $(\forall x)(\fail(x) \leq \fail(\Xi_s(x))$. 
\end{definition}

\Begin{fact}{club of opposite player}
Let $\kappa$ be a regular cardinal and $\lambda \leq \kappa$ be an ordinal. Let $(\Gamma,\decode,\GC_{\beta,\gamma} : \beta < \lambda, \gamma < \kappa)$ be a coding system for ${}^\lambda\kappa$.

If $r \in S^1$, then there is a club $C \subseteq \kappa$ (obtained uniformly from $r$) so that for all $\delta \in C$, for all $\beta < \min\{\lambda,\delta\}$ and $\gamma < \delta$,
$$\Xi_r\left[\bigcap_{\beta' < \beta}\bigcup_{\gamma' < \gamma}\GC_{\beta',\gamma'}\right] \subseteq \bigcap_{\beta' \leq \beta} \bigcup_{\gamma' < \delta} \GC_{\beta',\gamma'}$$

Suppose $s \in S^2$, then there is a club $C \subseteq \kappa$ so that for all $\delta \in C$, for all $\beta < \min\{\lambda,\delta\}$, for all $\gamma < \delta$, 
$$\Xi_s\left[\bigcap_{\beta' \leq \beta}\bigcup_{\gamma' < \gamma}\GC_{\beta',\gamma'}\right] \subseteq \bigcap_{\beta' \leq \beta}\bigcup_{\gamma' < \delta} \GC_{\beta',\gamma'}.$$
\end{fact}

\begin{proof}
Suppose $r \in S^1$. For each $\beta < \lambda$ and $\gamma < \kappa$, let
$$R_{\beta,\gamma} = \bigcap_{\beta'<\beta}\bigcup_{\gamma'<\gamma}\GC_{\beta',\gamma'}.$$
So if $x \in R_{\beta,\gamma}$, then $\decode(x)$ is good up to $\beta$ by successfully mapping each ordinal less than $\beta$ to some value less than $\gamma$. By property (5) of the good coding system, $\Delta$ is closed under less than $\kappa$ wellordered union and also intersections. Property (3) states that $\GC_{\beta',\gamma'} \in \Delta$ for all $\beta' < \lambda$ and $\gamma' < \kappa$. Thus $R_{\beta,\gamma} \in \Delta$. Then $\Xi_r[R_{\beta,\gamma}]$ is $\exists^\reals\Delta$. 

Claim 1: For all $\beta' \leq \beta$, $\Xi_r[R_{\beta,\gamma}] \subseteq \GC_{\beta'}$. 

To see Claim 1: Note that $x \in R_{\beta,\gamma}$ implies that $\fail(x) \geq \beta$. Since $r \in S^1$, $\fail(\Xi_r(x)) > \beta$. Thus for all $\beta ' \leq \beta$, $\Xi_r(x) \in \GC_{\beta'}$. This shows Claim 1.

Therefore $\Xi_r[R_{\beta,\gamma}]$ is a $\exists^\reals\Delta$ subset of $\GC_{\beta'}$ for each $\beta' \leq \beta$. By the boundedness property (4) of a good coding system, there is some $\epsilon_{\beta'} < \kappa$ so that $\Xi_r[R_{\beta,\gamma}] \subseteq \bigcup_{\gamma' < \epsilon_{\beta'}} \GC_{\beta',\gamma'}$. 

Since $\kappa$ is regular, let $\Upsilon(\beta,\gamma)$ be the least $\epsilon < \kappa$ so that $\beta < \epsilon$, $\gamma < \epsilon$, and for all $\beta' \leq \beta$, $\Xi_r[R_{\beta,\gamma}] \subseteq \bigcup_{\gamma' < \epsilon}\GC_{\beta',\gamma'}$. Thus $\Upsilon : \lambda \times \kappa \rightarrow \kappa$ is a well defined function. 

Let $C = \{\delta : (\forall \beta < \delta)(\forall \gamma < \delta)(\Upsilon(\min\{\lambda,\beta\},\gamma) < \delta)\}$. 

Claim 2: $C$ is a club. 

To see Claim 2: Let $\alpha < \kappa$. Let $\alpha_0 = \alpha$. If $\alpha_n$ has been defined, let $\alpha_{n + 1} = \Upsilon(\min\{\lambda,\alpha_n\},\alpha_n)$. Let $\alpha_\infty = \sup\{\alpha_n : n \in \omega\}$. By definition of $\Upsilon$, $\alpha_n < \alpha_{n + 1}$ for all $n$. Let $\beta < \min\{\lambda, \alpha_\infty\}$ and $\gamma < \alpha_\infty$. Then there is some $n$ so that $\beta < \alpha_n$ and $\gamma < \alpha_n$. Then $\Upsilon(\min\{\lambda,\beta\},\gamma) \leq \Upsilon(\min\{\lambda,\alpha_n\},\alpha_n) = \alpha_{n + 1} < \alpha_\infty$. This shows that $\alpha_\infty \in C$. As $\alpha < \alpha_\infty$, $C$ is unbounded. It is straightforward to show that $C$ is closed. Claim 2 has been shown.

Fix a $\delta \in C$. Pick a $\beta < \min\{\lambda,\delta\}$ and a $\gamma < \delta$. Suppose $x \in R_{\beta,\gamma}$. Let $\beta' \leq \beta$. Then 
$$\Xi_r(x) \in \bigcup_{\gamma' < \Upsilon(\beta,\gamma)} \GC_{\beta',\gamma'} \subseteq \bigcup_{\gamma' < \delta}\GC_{\beta',\gamma'}.$$ 
Since $\beta' \leq \beta$ was arbitrary,
$$\Xi_r(x) \in \bigcap_{\beta'\leq\beta} \bigcup_{\gamma'<\delta}\GC_{\beta',\gamma'}.$$
This completes the proof of the result for $S^1$. 

The argument for $S^2$ is similar.
\end{proof}

\Begin{remark}{club of opposite remark}
Let $r \in S^1$. Let $C$ be the club produced by Fact \ref{club of opposite player}. This club $C$ is called the club on which the opposite player (Player 2) has taken control of the output. 

This means the following: In most applications, $r$ codes a Lipschitz function coming from a Player 1 winning strategy. (Hence the notation $S^1$.) Let $\delta \in C$. Fact \ref{club of opposite player} states that for any code $y$ so that $\decode(y)$ successfully defines a function up to some $\beta < \min\{\lambda,\delta\}$ and takes value some $\gamma$ strictly below $\delta$, then the Player 1 strategy coded by $r$ when played against $y$ produces some code $x$ which successfully defines a function at and below $\beta$ and still takes value below $\delta$. 

This observation is used in a game to show that although Player 1 may have a winning strategy which determines that a final sequence will land in a certain payoff set, there is a club $C$ so that if the opposite player (Player 2) plays suitable codes for functions through this club, Player 2 actually determines the value of the final sequence. A similar statement holds if Player 2 has the winning strategy in this game.

This is the main idea of the following game of Martin to establish the partition properties and of an earlier game of Solovay to show the club filter is an ultrafilter.
\end{remark}

\Begin{definition}{sup of block function}
Suppose $\kappa$ is a regular cardinal and $\lambda$ is such that $\omega\cdot\lambda < \kappa$. Suppose $f \in {}^{\omega\cdot\lambda}\kappa$. Let $\block : {}^{\omega\cdot\lambda}\kappa \rightarrow {}^\lambda\kappa$ be defined by $\block(f)(\alpha) = \sup\{f(\omega\cdot\alpha + k) : k \in \omega\}$.

Suppose $f,g \in {}^{\omega\cdot\lambda}\kappa$. Let $\joint : {}^{\omega\cdot\lambda}\kappa \times {}^{\omega\cdot\lambda}\kappa \rightarrow {}^\lambda\kappa$ be defined by 
$$\joint(f,g)(\alpha) = \sup\{f(\omega\cdot\alpha + k),g(\omega\cdot\alpha + k) : k \in \omega\}.$$
\end{definition}

\Begin{theorem}{martin partition criterion}
(Martin) Assume $\ZF + \AD$. Suppose $\lambda, \kappa$ are ordinals such that $\omega\cdot\lambda \leq \kappa$. Suppose there is a good coding system $(\Gamma,\decode,\GC_{\beta,\gamma} : \beta \in \omega\cdot\lambda,\gamma < \kappa)$ for ${}^{\omega\cdot\lambda}\kappa$. Then $\kappa \rightarrow_* (\kappa)^\lambda_2$ holds.
\end{theorem}

\begin{proof}
Let $P : [\kappa]^\lambda_* \rightarrow 2$. Consider the following game
$$\begin{tikzpicture}
\node at (0,0) {$G$};

\node at (1,.5) {I};
\node at (2,.5) {$x_0$};
\node at (4,.5) {$x_1$};
\node at (6,.5) {$x_2$};
\node at (8,.5) {$x_3$};
\node at (11,.5) {$x$};

\node at (1,-.5) {II};
\node at (3,-.5) {$y_0$};
\node at (5,-.5) {$y_1$};
\node at (7,-.5) {$y_2$};
\node at (9,-.5) {$y_3$};
\node at (11,-.5) {$y$};
\end{tikzpicture}$$
where Player 1 and 2 take turns playing integers. Player 1 produces the real $x$, and Player 2 produces the real $y$. Player 1 wins the game if and only the conjunction of the following two conditions hold:

\noindent (1) $\fail(x) > \fail(y) \vee \fail(x) = \fail(y) = \infty$.

\noindent (2) $(\fail(x) = \fail(y) = \infty) \Rightarrow P(\joint(\decode(x),\decode(y))) = 0$.

(Case I) Suppose Player 1 has a winning strategy $\sigma$ in this game.

The strategy $\sigma$ induces a Lipschitz function $\Xi_\sigma$. Explicitly, $\Xi_\sigma(y)$ is just the response of Player 1 using $\sigma$ when Player 2 plays using $y$. Note that $\sigma \in S^1$. 

Let $C$ be the club from Fact \ref{club of opposite player}. On $C$, Player 2 takes control of the output in the sense below. Let $D$ be the limit points of $C$. 

Let $f \in [D]^{\lambda}_*$. Since $f$ is of the correct type, let $F : \lambda \times \omega \rightarrow \kappa$ be a witness to $f$ having uniform cofinality $\omega$. For each $\alpha < \lambda$, let $\nu_\alpha = \sup\{f(\xi) : \xi < \alpha\}$. Let $g(\omega\cdot\alpha)$ be the least element of $C$ greater than $\max\{\nu_\alpha,F(\alpha,0)\}$. Suppose $g(\omega\cdot\alpha + k)$ has been defined for some $k \in \omega$. Let $g(\omega\cdot\alpha + k + 1)$ be the least element of $C$ greater than $\max\{g(\omega\cdot\alpha+k),F(\alpha, k + 1)\}$. Now $g \in [C]^{\omega\cdot\lambda}$, $g$ is discontinuous everywhere, and $\block(g) = f$. 

By property (2) of the good coding system, there is some $y \in \GC$ so that $\decode(y) = g$. Let Player 2 play this $y$ against Player 1 using $\sigma$. Since $y \in \GC$, $\fail(y) = \infty$. Thus $\fail(\Xi_\sigma(r)) = \infty$.

For $\beta < \omega\cdot\lambda$, let $\epsilon_\beta = \sup\{g(\alpha) : \alpha < \beta\}$. Since $g$ is discontinuous everywhere, $\epsilon_\beta < g(\beta)$. Then 
$$\Xi_\sigma(y) \in \Xi_\sigma\left[\bigcap_{\beta' < \beta}\bigcup_{\gamma' < \epsilon_\beta}\GC_{\beta',\gamma'}\right] \subseteq \bigcap_{\beta'\leq\beta}\bigcup_{\gamma' < g(\beta)}\GC_{\beta',\gamma'}$$
by the definition of $C$. Hence $\decode(\Xi_\sigma(y))(\beta) < g(\beta)$. Thus for all $\alpha < \lambda$,
$$\sup\{\decode(\Xi_\sigma(y))(\omega\cdot\alpha + k), \decode(y)(\omega\cdot\alpha + k) : k \in \omega\} = \sup\{g(\omega\cdot\alpha + k) : k \in \omega\} = f(\alpha).$$
This shows that $\joint(\decode(\Xi_\sigma(y)),\decode(y)) = f$. Since $\sigma$ is a Player 1 winning strategy, one has that $P(f) = 0$. Since $f \in [D]^{\lambda}_*$ was arbitrary, $D$ is homogeneous for $P$ taking value $0$.

(Case II) Suppose Player 2 has a winning strategy $\tau$. 

Note that Player 2 wins if and only if the disjunction of the following holds:

\noindent (1) $\fail(x) < \fail(y) \vee (\fail(x) = \fail(y) < \infty)$.

\noindent (2) $\fail(x) = \fail(y) = \infty \wedge P(\joint(\decode(x),\decode(y))) = 1$. 

Therefore, $\tau \in S^2$. Let $C$ be the club coming from Fact \ref{club of opposite player}. $C$ is the club for which Player 1 takes control of the output. One may assume that $C$ consists of only limit ordinals. Let $D$ be the limit points of $C$.

Let $f \in [D]^\lambda_*$. As before, there is a $g \in [C]^{\omega\cdot\lambda}$ such that $\block(g) = f$. Let $x \in \GC$ be such that $\decode(x) = g$. Since $\fail(x) = \infty$, $\fail(\Xi_\tau(x)) = \infty$ and $\Xi_\tau(x) \in \GC$ as well.

Let $\beta < \omega\cdot\lambda$. Note that $g(\beta) + 1 < g(\beta + 1)$
$$\Xi_\tau(x) \in \Xi_\tau\left[\bigcap_{\beta' \leq \beta}\bigcup_{\gamma' < g(\beta) + 1} \GC_{\beta',\gamma'}\right] \subseteq \bigcap_{\beta'\leq\beta}\bigcup_{\gamma' < g(\beta + 1)} \GC_{\beta',\gamma'}$$
by definition of $C$. Hence $\decode(\Xi_\tau(x))(\beta) < g(\beta + 1)$. Thus for all $\alpha < \lambda$, 
$$\sup\{\decode(x)(\omega\cdot\alpha + k), \decode(\Xi_\tau(x))(\omega\cdot\alpha + k) : k \in \omega\} = \sup\{g(\omega\cdot\alpha + k) : k \in \omega\} = f(\alpha).$$
This shows that $\joint(\decode(x),\decode(\Xi_\tau(x))) = f$. Since $\tau$ is a Player 2 winning strategy, one has that $P(f) = 1$. $D$ is homogeneous for $P$ taking value $1$.

The proof is complete.
\end{proof}

\Begin{theorem}{club code uniformization}
(Almost everywhere uniformization on codes) Assume $\ZF + \AD$. Let $\kappa$ be a regular cardinal and $\lambda < \kappa$. Let $(\Gamma,\decode,\GC_{\beta,\gamma} : \beta < \omega\cdot\lambda, \gamma < \kappa)$ be a good coding system for ${}^{\omega\cdot\lambda}\kappa$. Let $R \subseteq [\kappa]^\lambda_* \times \reals$ be a relation. 

There is a club $C \subseteq \kappa$ and a Lipschitz continuous function $F : \reals \rightarrow \reals$ so that for all $x \in \GC$ with $\decode(x) \in [C]^{\omega\cdot\lambda}$ and $\block(\decode(x)) \in [C]^\lambda_* \cap \dom(R)$, $R(\block(\decode(x)),F(x))$.
\end{theorem}

\begin{proof}
Define a game $G$ as follows:
$$\begin{tikzpicture}
\node at (0,0) {$G$};

\node at (1,.5) {I};
\node at (2,.5) {$x_0$};
\node at (4,.5) {$x_1$};
\node at (6,.5) {$x_2$};
\node at (8,.5) {$x_3$};
\node at (11,.5) {$x$};

\node at (1,-.5) {II};
\node at (3,-.5) {$y_0,z_0$};
\node at (5,-.5) {$y_1,z_1$};
\node at (7,-.5) {$y_2,z_2$};
\node at (9,-.5) {$y_3,z_3$};
\node at (11,-.5) {$y,z$};
\end{tikzpicture}$$
Player 2 wins if and only if the disjunction of the following hold:

\noindent (1) $\fail(x) < \fail(y) \vee \fail(x) = \fail(y) < \infty$.

\noindent (2) $(\fail(x) = \fail(y) = \infty) \wedge (\joint(\decode(x),\decode(y)) \in \dom(R) \Rightarrow R(\joint(\decode(x),\decode(y)),z))$.

Claim 1: Player 2 has a winning strategy. 

To prove Claim 1: Suppose not. By $\AD$, Player 1 has a winning strategy $\sigma$. Note that Player 1 winning condition is the conjunction of the follow:

\noindent (1) $\fail(y) < \fail(x) \vee \fail(x) = \fail(y) = \infty$.

\noindent (2) $\fail(x) = \fail(y) = \infty \Rightarrow [\joint(\decode(x),\decode(y)) \in \dom(R) \wedge \neg R(\joint(\decode(x),\decode(y)),z)]$.

Let $\Xi_\sigma$ be the associated continuous function. Note that
$$(\forall y)(\forall z)[\fail(y) \leq \fail(\Xi_\sigma(y,z)) \wedge (\fail(y) < \infty \Rightarrow \fail(y) < \fail(\Xi_\sigma(y,z)))].$$
That is, if one fixes any $z$, the associated function $F_z(y) = \Xi_\sigma(y,z)$ belongs to $S^1$. With a small modification to the argument of Fact \ref{club of opposite player}, there is a club $C$ with the property that for all $\delta \in C$, for all $z \in \reals$, for all $\beta < \min\{\lambda,\delta\}$, and $\gamma < \delta$
$$F_z\left[\bigcap_{\beta'<\beta}\bigcup_{\gamma' < \gamma}\GC_{\beta',\gamma'}\right] \subseteq \bigcap_{\beta'\leq\beta}\bigcup_{\gamma' < \gamma}\GC_{\beta',\gamma'}.$$
That is, $C$ is a club for which Player 2 controls the output against $F_z$ for every $z \in \reals$.

Let $D$ be the limit points of $C$. Let $f \in [D]^{\lambda}_*$. As in the proof of Theorem \ref{martin partition criterion}, let $g \in [D]^{\omega\cdot\lambda}$ be such that $\block(g) = f$. Let $y \in \GC$ be such that $\decode(y) =g$. By the same argument as in Theorem \ref{martin partition criterion}, one has that for all $z$, $\joint(\decode(\Xi_\sigma(y,z)),\decode(y)) = f$. Since $\sigma$ is a Player 1 winning strategy, one must have that $f \in \dom(R)$. Since $f \in \dom(R)$, there is some $z^*$ so that $R(f,z^*)$. If Player 2 plays $(y,z^*)$, then Player 1 loses using $\sigma$. This contradicts $\sigma$ being a Player 1 winning strategy. Claim 1 has been shown.

Let $\tau$ be a Player 2 winning strategy in this game. Let $G(x) = \pi_1[\Xi_\tau(x)]$, where $\pi_1$ refers to the projection onto the first coordinate. Then $G$ is a Lipschitz continuous function satisfying the condition from the definition of $S^2$. By Fact \ref{club of opposite player}, let $C$ be the club on which Player 2 takes controls of the output. Again one may assume $C$ consists of only limit ordinals. Let $D$ be the limit points of $C$.

Let $f \in [D]^{\lambda}_* \cap \dom(R)$. Let $g \in [C]^{\omega\cdot\lambda}$ be such that $\block(g) = f$. Let $x \in \GC$ be such that $\decode(x) = g$. By the argument in Theorem \ref{martin partition criterion}, $\joint(\decode(x),\decode(G(x))) = f$. Since $\tau$ is a Player 2 winning strategy, one must have that $R(f,\pi_2(\Xi_r(x)))$. 

Define $F(x) = \pi_2(\Xi_r(x))$. Then $F$ is a Lipschitz function with the desired uniformization property.
\end{proof}

\Begin{theorem}{absorbing function on function spaces}
Let $\kappa$ be a regular cardinal and $\lambda \leq \kappa$. Suppose $(\Gamma, \decode, \GC_{\beta,\gamma} : \beta < \lambda, \gamma < \kappa)$ is a good coding system for ${}^\lambda\kappa$. Let $M \models \AD$ be an inner model containing all the reals and within $M$, $(\Gamma,\decode,\GC_{\beta,\gamma} : \beta < \lambda,\gamma<\kappa)$ is a good coding system.

Then for any $\Phi: [\kappa]^\lambda \rightarrow \kappa$, there is a club $D$, necessarily in $M$ by the Moschovakis coding lemma, so that $\Phi \upharpoonright [D]^\lambda_* \in M$.
\end{theorem}

\begin{proof}
The hypothesis implies that $\kappa < \Theta^M$. Let $\pi : \reals \rightarrow \kappa$ be a surjection in $M$. Define a relation $R \subseteq [\kappa]^\lambda_* \times \reals$ by
$$R(f,x) \Leftrightarrow \Phi(f) = \pi(x)$$
Let the Lipschitz function $F : \reals \rightarrow \reals$ and club $C \subseteq \kappa$ be the objects given by Theorem \ref{club code uniformization}. Let $D$ be the set of limit points of $C$. 

Let $f \in [D]^\lambda_*$. Let $x \in \GC$ be any function such that $\decode(x) \in [C]^{\omega\cdot\lambda}_*$ and $\block(\decode(x)) = f$. Then $R(\block(\decode(x)), F(x))$. This means that $\Phi(f) = \pi(F(x))$. 

Note that for any $x$ and $x'$ so that $\decode(x),\decode(x') \in [C]^{\omega\cdot\lambda}_*$ and $\block(\decode(x)) = \block(\decode(x'))$, one has that $\pi(F(x)) = \pi(F(x')) = \Phi(\block(\decode(x)))$. 

By the Moschovakis coding lemma and the fact that $\pi \in M$, $C$ and hence $D$ belongs to $M$. So within $M$, one can define $\Phi \upharpoonright [D]^{\lambda}_*$ as follows, $\Phi(f) = \gamma$ if and only if there exists some $x \in \GC$ so that $\decode(x) \in [C]^{\omega\cdot\lambda}_*$, $\block(\decode(x)) = f$, and $\pi(F(x)) = \gamma$. By the above, this is well defined and works.
\end{proof}

One will say that a cardinal $\kappa$ is ``reasonable'' if one is in the situation where there exists a good coding system:

\Begin{definition}{reasonable cardinals}
(Jackson) Let $\kappa$ be a regular cardinal and $\lambda \leq \kappa$. $\kappa$ is $\lambda$-reasonable if and only if there is a good coding system for ${}^\lambda\kappa$. 
\end{definition}

\section{Reasonableness at $\omega_1$}\label{reasonableness at omega1}

\Begin{definition}{WO definitions}
Fix some recursive pairing function $\pi : \omega \times \omega \rightarrow \omega$. A real $x \in \cantorspace$ codes a relation $R_x$ defined as follows: $R_x(m,n) \Leftrightarrow x(\pi(m,n)) = 1$. 

Let the domain of $x$ be $\dom(x) = \{n : (\exists m)(R_x(m,n) \vee R_x(m,n))\}$. 
Let $\mathrm{LO}$ be the set of reals $x$ so that $R_x$ is a linear ordering on its domain. Let $\WO$ be the set of reals $x$ so that $R_x$ is a well ordering. $\mathrm{LO}$ is an arithmetic set of reals. $\mathrm{WO}$ is $\Pi_1^1$. 

If $x \in \WO$, then $\ot(x)$ is the order type of $x$. If $\beta < \ot(x)$, let $n_\beta^x$ denote the unique element of $\dom(x)$ whose rank according to $R_x$ is $\beta$. 

If $\alpha < \omega_1$, then let $\WO_\alpha = \{x \in \WO : \ot(x) = \alpha\}$. $\WO_{<\alpha} = \{x \in \WO : \ot(x) < \alpha\}$. $\WO_{\leq \alpha}$, $\WO_{>\alpha}$, and $\WO_{\geq \alpha}$ are defined similarly. $\WO_\alpha$, $\WO_{<\alpha}$, and $\WO_{\leq \alpha}$ are $\borel$ (in any element of $\WO_\alpha$).
\end{definition}

\Begin{definition}{pi11 norm on WO}
The ordertype function $\ot : \WO \rightarrow \omega_1$ is a $\Pi_1^1$ norm. Let $\preceq$ denote the induced prewellordering: $x \preceq y$ if and only if $\ot(x) \leq \ot(y)$. 

Being a $\Pi_1^1$-norm implies that there is a $\Sigma_1^1$ relation $\leq_{\Sigma_1^1}$ and a $\Pi_1^1$ relation $\leq_{\Pi_1^1}$ so that 
$$y \in \WO \Rightarrow (\forall x)[(x \in \WO \wedge x \preceq y) \Leftrightarrow (x \leq_{\Sigma_1^1} y) \Leftrightarrow (x \leq_{\Pi_1^1} y)].$$
\end{definition}

It is useful to have a concrete coding of a nice collection of club subsets of $\omega_1$.

\Begin{fact}{strategy outputting large WO gives club}
Let $\tau$ be a Player 2 strategy with the property that for all $w \in \mathrm{WO}$, $\tau(w) \in \mathrm{WO} \wedge \mathrm{ot}(w) < \mathrm{ot}(\tau(w))$. Let $C_\tau = \{\eta : (\forall w \in \mathrm{WO}_{<\eta})(\mathrm{ot}(\tau(w)) < \eta)\}$. $C_\tau$ is a club.
\end{fact}

\begin{proof}
Let $R_\alpha = \{\tau(w) : w \in \mathrm{WO}_{<\alpha}\}$. $R_\alpha \subseteq \mathrm{WO}$ is $\analytic$. By boundedness, $\sup \{\mathrm{ot}(v) : v \in R_\alpha\} < \omega_1$. 

Let $\Phi(\alpha) = \sup\{\mathrm{ot}(\tau(w)) : w \in \mathrm{WO}_{<\alpha}\} + 1$. By the above, $\Phi(\alpha)$ is defined. 

It is clear that $C_\tau$ is closed.

Let $\alpha < \omega_1$. Let $\alpha_0 = \alpha$. Let $\alpha_{n + 1} = \Phi(\alpha_n)$. Let $\eta = \sup\{\alpha_n : n \in \omega\}$. Let $\beta < \eta$. There is some $n$ so that $\beta < \alpha_n$. Let $w \in \mathrm{WO}_\beta \subseteq \mathrm{WO}_{<\alpha_n}$. Thus $\mathrm{ot}(\tau(w)) < \Phi(\alpha_n) = \alpha_{n + 1} < \eta$. So $\eta \in C_\tau$. This shows that $C_\tau$ is unbounded.
\end{proof}

\Begin{definition}{reals coding clubs}
Let $\mathsf{clubcode}_{\omega_1}$ be the set of Player 2 strategies so that for all $w \in \mathrm{WO}$, $\tau(w) \in \mathrm{WO} \wedge \mathrm{ot}(w) < \mathrm{ot}(\tau(w))$.

Note that $\mathsf{clubcode}_{\omega_1}$ is a $\Pi_2^1$ set of reals.
\end{definition}

\Begin{fact}{membership in a club}
Assume that $\tau \in \mathsf{clubcode}_{\omega_1}$. The relation $S(w)$ defined by $w \in \WO \wedge \ot(w) \in C_\tau$ is $\mathbf{\Pi}_1^1$. 

Assume $\AC_\omega^\reals$. If $\alpha < \omega_1$, then the relation $T_\alpha(w)$ defined by $w \in \WO_{<\alpha} \wedge \ot(w) \in C_\tau$ is $\borel$.
\end{fact}

\begin{proof}
Note that $S(w)$ holds if and only if $w \in \WO \wedge (\forall v)(v <_{\Sigma_1^1} w \Rightarrow \tau(v) <_{\Pi_1^1} w)$, where $<_{\Sigma_1^1}$ and $<_{\Pi_1^1}$ are $\Sigma_1^1$ and $\Pi_1^1$ relations, respectively, coming from the ordertype function, $\ot$, being a $\Pi_1^1$-norm.

Note that $\alpha \cap C_\tau$ is a countable set. Using $\AC_\omega^\reals$, let $r \in \omega$ be such that $\{\ot(r_n) : n \in \omega\} = \alpha \cap C_\tau$. (Here $r_n(m) = r(\langle n,m\rangle)$, where $\langle \cdot, \cdot \rangle : \omega \times \omega \rightarrow \omega$ is a recursive bijective pairing function.) One can see that $T_\alpha$ is $\borel$ using as parameters $r$ and some element of $\WO_\alpha$. 
\end{proof}

\Begin{fact}{coding a club by a increasing strategy}
Assume $\ZF + \AD$. Let $C \subseteq \omega_1$ be a club. There is a $\tau \in \mathsf{clubcode}_{\omega_1}$ so that $C_\tau \subseteq C$. 
\end{fact}

\begin{proof}
Consider the game $G_C$:
$$\begin{tikzpicture}
\node at (0,0) {$G_C$};

\node at (1,.5) {I};
\node at (2,.5) {$w_0$};
\node at (4,.5) {$w_1$};
\node at (6,.5) {$w_2$};
\node at (8,.5) {$w_3$};
\node at (11,.5) {$w$};

\node at (1,-.5) {II};
\node at (3,-.5) {$v_0$};
\node at (5,-.5) {$v_1$};
\node at (7,-.5) {$v_2$};
\node at (9,-.5) {$v_3$};
\node at (11,-.5) {$v$};
\end{tikzpicture}$$
Player 2 wins this game if and only if 
$$w \in \mathrm{WO} \Rightarrow (v \in \mathrm{WO} \wedge \mathrm{ot}(v) > \mathrm{ot}(w) \wedge \mathrm{ot}(v) \in C).$$

Claim 1: Player 2 has the winning strategy. 

Suppose $\sigma$ is a Player 1 winning strategy. Note that $\sigma[\reals] \subseteq \mathrm{WO}$ and is $\analytic$. By boundedness, $\gamma = \sup \{\mathrm{ot}(w) : w \in \sigma[\reals]\} < \omega_1$. Since $C$ is unbounded, let $\delta \in C$ be such that $\delta > \gamma$. Let $v \in \mathrm{WO}_\delta$. If Player 2 plays $v$ against $\sigma$, $v \in \mathrm{WO}$, $\mathrm{ot}(v) = \delta \in C$, and $\mathrm{ot}(\sigma(v)) \leq \gamma < \delta = \mathrm{ot}(v)$. Player 2 has won. This contradicts $\sigma$ being a Player 1 winning strategy. 

Let $\tau$ be a Player 2 winning strategy. It is clear that $\tau \in \mathsf{clubcode}_{\omega_1}$. 

Claim 2: $C_\tau \subseteq C$. 

Suppose $\eta \in C_\tau$. Let $\beta < \eta$. Let $w \in \mathrm{WO}_\beta$. Then $\beta < \mathrm{ot}(\tau(w)) < \eta$ and $\mathrm{ot}(\tau(w)) \in C$. Since $\beta < \eta$ was arbitrary and $C$ is closed, $\eta \in C$.
\end{proof}

\Begin{fact}{analytic set of club codes explicit club}
Suppose $A \subseteq \mathsf{clubcode}_{\omega_1}$ is $\analytic$. Then uniformly from $A$, there is a club $C$ so that for all $\tau \in A$, $C \subseteq C_\tau$. 
\end{fact}

\begin{proof}
For each $\beta < \omega_1$, let $R_\beta = \{\tau(w) : \tau \in A \wedge w \in \mathrm{WO}_{<\beta}\}$. $R_\beta \subseteq \mathrm{WO}$ and is $\analytic$. By the boundedness principle, $\sup\{\mathrm{ot}(v) : v \in R_\beta\} < \omega_1$. 

Let $\Phi(\beta) = \sup\{\mathrm{ot}(v) : v \in R_\beta\} + 1$. 

Let $C = \{\eta : (\forall \beta < \eta)(\Phi(\beta) < \eta)\}$. As before, one can check that $C$ is a club subset of $\omega_1$. 

Fix $\tau \in A$. Let $\eta \in C$. Suppose $\beta < \eta$. Let $w \in \mathrm{WO}_\beta$. Then $\beta < \mathrm{ot}(\tau(w)) < \Phi(\beta) <  \eta$. Thus $\eta \in C_\tau$. 
\end{proof}

\Begin{fact}{omega1 club selection}
Assume $\ZF + \AD$. Let $\langle \mathcal{A}_\alpha : \alpha < \omega_1\rangle$ be such that each $\mathcal{A}_\alpha$ is a nonempty $\subseteq$-downward closed collection of club subsets of $\omega_1$. Then there is a sequence $\langle C_\alpha : \alpha < \omega_1 \rangle$ with each $C_\alpha \subseteq \omega_1$ a club subset of $\omega_1$ and $C_\alpha \in \mathcal{A}_\alpha$. 

In particular: Let $\mu$ be the club measure on $\omega_1$. Let $\langle A_\alpha : \alpha < \omega_1\rangle$ be a sequence of sets in $\mu$, that is each $A_\alpha$ contains a club subset of $\omega_1$. Then there is a sequence $\langle C_\alpha : \alpha < \omega_1\rangle$ with each $C_\alpha \subseteq \omega_1$ a club subset of $\omega_1$ so that $C_\alpha \subseteq A_\alpha$. 
\end{fact}

\begin{proof}
Consider the game 
$$\begin{tikzpicture}
\node at (0,0) {$G$};

\node at (1,.5) {I};
\node at (2,.5) {$w(0)$};
\node at (4,.5) {$w(1)$};
\node at (6,.5) {$w(2)$};
\node at (8,.5) {$w(3)$};
\node at (11,.5) {$w$};

\node at (1,-.5) {II};
\node at (3,-.5) {$z(0)$};
\node at (5,-.5) {$z(1)$};
\node at (7,-.5) {$z(2)$};
\node at (9,-.5) {$z(3)$};
\node at (11,-.5) {$z$};
\end{tikzpicture}$$
where Player 2 wins if and only if $w \in \WO \Rightarrow (z \in \mathrm{clubcode}_{\omega_1} \wedge C_z \in \mathcal{A}_{\ot(w)})$.

Claim 1: Player 2 has a winning strategy in this game. 

To prove this, suppose otherwise that Player 1 has a winning strategy $\sigma$. Note that $\sigma[\reals] \subseteq \WO$. Since $\sigma[\reals]$ is a $\analytic$ subset of $\WO$, by the boundedness principle, there is a $\zeta \in \omega_1$ so that for all $w \in \sigma[\reals]$, $\ot(w) < \zeta$. By $\AC_\omega^\reals$, pick $\langle C_\alpha : \alpha < \zeta\rangle$ with the property that for all $\alpha < \zeta$, $C_\alpha \in \mathcal{A}_\alpha$. Let $C = \bigcap_{\alpha < \zeta} C_\alpha$ which is also a club. By Fact \ref{coding a club by a increasing strategy}, there is some $z \in \mathrm{clubcode}_{\omega_1}$ so that $C_z \subseteq C$. Note that since for each $\alpha < \zeta$, $C_z \subseteq C_\alpha$,  $C_\alpha \in \mathcal{A}_\alpha$, and $\mathcal{A}_\alpha$ is $\subseteq$-downward closed, one has that $C_z \in \mathcal{A}_\alpha$ for each $\alpha < \zeta$. Thus Player 2 wins against $\sigma$ by playing $z$. Contradiction.

Thus let $\tau$ be a Player 2 winning strategy in this game. For each $\alpha < \omega_1$, let $P_\alpha = \tau[\WO_\alpha]$. Note that $P_\alpha \subseteq \mathrm{clubcode}_{\omega_1}$ with the property that for all $z \in P_\alpha$, $C_z \in \mathcal{A}_\alpha$. Note that $P_\alpha$ is $\analytic$. By Fact \ref{analytic set of club codes explicit club}, there is a uniform procedure to obtain from the set $P_\alpha$, a club $C_\alpha$ with the property that $C_\alpha \subseteq C_z$ for all $z \in P_\alpha$. In particular since $\mathcal{A}_\alpha$ is $\subseteq$-downward closed, $C_\alpha \in \mathcal{A}_\alpha$. This completes the proof.
\end{proof}

\Begin{fact}{omega1 is bounded reasonable}
(Martin) Assume $\ZF + \AD$. Let $\lambda < \omega_1$. Then $\omega_1$ is $\lambda$-reasonable. 
\end{fact}

\begin{proof}
Note that $\omega_1$ is regular by $\AC_\omega^\reals$. Let $\Gamma = \analytic$. Note that $\borel$ is closed under countable unions. Thus (1) and (5) hold.

For each $x \in \reals$, let $x_n \in \reals$ be defined by $x_n(k) = x(\langle n,k\rangle)$, where $\langle \cdot,\cdot \rangle : \omega \times \omega \rightarrow \omega$ denotes a fixed recursive bijective pairing function. 

Fix some $w \in \WO_\lambda$ with $\dom(w) = \omega$. 

Define $\decode(x)(\alpha,\beta)$ if and only if $x_{n^w_\alpha} \in \WO_\beta$. Suppose $f \in {}^\lambda\omega_1$. By $\AC_\omega^\reals$, for each $\beta < \lambda$, let $u_\beta \in \WO_{f(\beta)}$. Let $x \in \reals$ be such that $x_{n^w_\beta} = u_\beta$. Thus $\decode(x) = f$. This shows (2). 

Let $\GC_{\beta,\gamma} = \{x : x_{n^w_{\beta}} \in \WO_\gamma\}$. Note that $\GC_{\beta,\gamma}$ is $\borel$. It satisfies the property in (3).

Suppose $A \in \analytic$ and $A \subseteq \GC_\beta$. Let $R = \{x_{n^w_\beta} : x \in A\}$. Then $R$ is a $\analytic$ subset of $\WO$. By the usual boundedness principle, there is some $\delta < \omega_1$ so that for all $v \in R$, $\ot(v) < \delta$. Then $A \subseteq \bigcup_{\gamma < \delta}\GC_{\beta,\gamma}$. 

It has been shown that $(\analytic, \decode, \GC_{\beta,\gamma} : \beta < \lambda, \gamma < \omega_1)$ is a good coding system for ${}^\lambda\omega_1$. 
\end{proof}

\Begin{corollary}{omega1 is weak partition}
(Martin) Assume $\ZF + \AD$. $\omega_1$ is a weak partition cardinal, i.e. for all $\lambda < \omega_1$, $\omega_1 \rightarrow_* (\omega_1)^{\lambda}_2$. 
\end{corollary}

\Begin{corollary}{normal measure on omega1}
Assume $\ZF + \AD$. The club filter $\mu = W^{\omega_1}_\omega$ is the unique countably complete normal measure on $\omega_1$. 
\end{corollary}

\begin{proof}
By Fact \ref{characterize normal measure by partition}.
\end{proof}

\Begin{fact}{borelness of function passing through club}
Let $\alpha < \omega_1$. Fix a good coding system for ${}^{\omega\cdot\alpha} \omega_1$. Suppose $z \in \mathsf{clubcode}_{\omega_1}$. Let $D$ be the set of limit points of $C_z$. Let $f \in [D]^{\alpha}_*$. The set of $\{x \in \reals : \decode(x) \in [C_z]^{\omega\cdot\alpha} \wedge \block(\decode(x)) = f\}$ is a $\borel$ set
\end{fact}

\begin{proof}
Fix any $x^*$ so that $\block(\decode(x^*)) = f$. Since $\omega_1$ is regular, $\sup(f) = \gamma < \omega_1$. Using $x^*$ as a parameter, one can now check that the desired set is $\borel$ by unraveling the definition of the coding system in Theorem \ref{omega1 is bounded reasonable} and using Fact \ref{membership in a club}.
\end{proof}

The original method of establishing that $\omega_1$ is $\omega_1$-reasonable involves sharps for reals which will be concisely reviewed below. See \cite{A-Nonconstructible-Delta13-Real} for more details. The reader who would prefer to avoid sharps can skip ahead to Theorem \ref{kechris reasonableness proof} to see the argument of Kechris that uses the Banach-Mazur game and category.

Let $\SCRL = \{\dot \in, \dot E\}$ be a language where $\dot \in$ is a binary relation symbol and $\dot E$ is a unary relation symbol.

For each $\SCRL$-formula $\varphi(w,v_1,...,v_{k})$, whose free variable are exactly those listed, then let $h_\varphi(v_1,...,v_k)$ be a formal function which will be called the Skolem function associated to $\varphi$.

Define the language $\SCRL^I = \{\dot \in, \dot E, \dot c_n : n \in \omega\}$ where for each $n \in \omega$, $\dot c_n$ is a distinct constant symbol.

Let $\mathsf{skolem}$ denote the smallest class of function $h(v_1,...,v_k)$ containing $h_\varphi$ for all $\varphi$ and closed under composition.

$\SCRL^S$ consist of $\dot \in$, $\dot E$, constants $c_n$ for all $n$, and new distinct constant symbols $h(c_{i_1},...,c_{i_k})$ for all $i_1 < ... < i_k$ in $\omega$ and $h \in \mathsf{skolem}$

One will assume one has fixed a recursive coding of formulas and terms of $\SCRL^S$. One will identify terms or formulas with their associated integer code.

Let $\dot \prec$ denote the $\SCRL$-formula that defines the canonical $L[\dot E]$ wellordering. Every formula $\phi$ of $\SCRL^S$ can be converted into a $\SCRL^I$ formula $\tilde \phi$ that roughly amounts to recursively replacing each Skolem term $h_\varphi(v_1,...,v_k)$ with its intended meaning that it should represent the $\dot \prec$-least solution to $\varphi$ if it exists and $\emptyset$ otherwise.

\Begin{definition}{arithmetic syntactical statement in sharp}
Let $A(T)$ assert the following:

\noindent (1) $T$ is the set of integer code for a complete and consistent theory extending $\ZF + V = L[\dot E]$. The statement $(\forall x)(x \in \dot E \Rightarrow x \in \omega)$ belongs to $T$.

\noindent (2) (Indiscernibility) For each $\SCRL$-formula $\varphi(v_1,...,v_k)$ and increasing sequences of integers $(a_1,...,a_k)$ and $(b_1,...,b_k)$, the following sentence belongs to $T$:
$$``\varphi(c_{a_1},...,c_{a_k}) \Leftrightarrow \varphi(c_{b_1},...,c_{b_k})".$$

\noindent (3) (Unbounded) Let $h(v_1,...,v_k) \in \mathsf{skolem}$. The following sentence belongs to $T$: 
$$``h(c_1,...,c_k) \in \mathrm{ON} \Rightarrow h(c_1,...,c_k) < c_{k + 1}".$$ 
(This is actually a $\SCRL^S$-sentence $\phi$. Precisely, one means that $\tilde \phi \in T$.)

\noindent (4) (Remarkability) Suppose $h(v_0,...,v_k,w_1,...,v_j) \in \mathsf{skolem}$. Then the following sentence belongs to $T$:
$$``[h(c_0,...,c_k,c_{a_1},...,c_{a_j}) \in \mathrm{ON} \wedge h(c_0,...,c_k,c_{a_1},...,c_{a_j}) < c_{a_1}]$$ 
$$\Rightarrow h(c_1,...,c_k,c_{a_1},...,c_{a_j}) = h(c_1,...,c_k,c_{b_1},...,c_{b_k})"$$
for all increasing sequences of integers $(a_1,...,a_j)$ and $(b_1,...,b_j)$ such that $k < a_1$ and $k < b_1$.
\newline

Note that $A(T)$ is an arithmetical statement asserting syntactical conditions on the $\SCRL^I$-theory $T$.
\end{definition}

Now assume $T$ is a set of integers so that $A(T)$.

Let $K$ be a linear ordering. For each $a \in K$, let $c_a$ be distinct new constant symbols. Let $\SCRL^{S,K}$ consists of $\dot \in$, $\dot E$, and new constant symbols $h(c_{a_1},...,c_{a_k})$ for each $h \in \mathsf{skolem}$ and increasing tuple $(a_1,...,a_k)$ in $K$.

An equivalence relation is defined on the new constants of $\SCRL^{S,K}$ by have $h_1(c_{a_1},...,c_{a_n}) \sim h_2(c_{a_1},...,c_{a_n})$ if and only $\tilde \phi \in T$ where $\phi$ is the $\SCRL^S$ formula $h_1(c_1,...,c_n) = h_2(c_1,...,c_n)$. 

The membership relation can be defined on Skolem constant by referring to $T$ in a similar manner.

By taking quotients, one obtains a structure denoted $\Gamma(T,K)$.

$K$ embeds into $\Gamma(T,K)$ by the map $j(a) = [c_a]_\sim$.

\Begin{fact}{wellfoundedness of indiscernible models}
Let $T \subseteq \omega$ be such that $A(T)$. Then $\Gamma(T,\alpha)$ is wellfounded for all ordinals $\alpha$ if and only if $\Gamma(T,\alpha)$ is wellfounded for all $\alpha < \omega_1$.
\end{fact}

Suppose $w \in \WO$ and $(\omega,<_w)$ is its associated wellordering relation on $\omega$. One can check that in this case one can find a structure on $\omega$ which is isomorphic to $\Gamma(T,\mathrm{ot}(w))$. This structure is recursive in $T$ and $w$ and produced uniformly from $T$ and $w$. In such case that $w \in \WO$, one will denote $\Gamma(T,w)$ to be this structure on $\omega$ recursive in $T$ and $w$. 

Suppose $M$ is an $\dot \in$ structure on $\omega$. For each $k \in \omega$, there is a recursive function $\mathsf{ST}$ so that $\mathsf{ST}(M,k)$ is a structure on $\omega$ isomorphic to $(k,\dot \in^M)$. 

Since $\Gamma(T,w)$ is considered as a structure on $\omega$, for each $k \in \Gamma(T,w)$ such that $\Gamma(T,w) \models k$ is an ordinal, $\mathsf{ST}(\Gamma(T,w),k)$ gives (uniformly) a structure on $\omega$ isomorphic to ordinal $k$ in $\Gamma(T,w)$.

Note that for each $\beta < \omega_1$, the set of $(T,w)$ such that $\beta$ is an initial segment of the ordinal of $\Gamma(T,w)$ is a $\borel$ set (using any element of $\WO_\beta$ as a parameter).

Consider formulas $\varphi(v,w)$ stating $v \in \mathrm{ON} \wedge \psi(v,w)$ where $\psi$ is some other $\SCRL$-formula. After fixing an enumeration of such formulas, one has an enumeration $\langle t_n : n \in \omega\rangle$ of Skolem constants whose intention is to name ordinals.

For each $\beta < \omega_1$ and $\gamma < \omega_1$, the set of $(T,w,n)$ such that $\beta + \omega$ is an initial segment of $\Gamma(T,w)$ and $t_n^{\Gamma(T,w)}(\beta) < \gamma$ is a $\borel$ set (using any elements of $\WO_\beta$ and $\WO_\gamma$). Similarly for $t_n^{\Gamma(T,w)}(\beta) = \gamma$. 

Let $\mathrm{B}(T)$ assert that for all $\alpha < \omega_1$, $\Gamma(T,\alpha)$ is wellfounded. This is equivalent to 
$$(\forall w)(w \in \WO \Rightarrow \Gamma(T,w) \text{ is wellfounded}).$$ 
This is $\Pi_2^1$. 

\Begin{definition}{sharp definition}
A set $T$ is a sharp of a real if and only if $A(T) \wedge B(T)$. Thus the statement that $T$ is a sharp of a real is $\Pi_2^1$. 

$T$ is the sharp of a real $x$ if and only if $T$ is a sharp of a real and for all $n \in \omega$, $``n \in \dot E" \in T$ if and only if $n \in x$. 

It can be shown that if there is a $T$ such that $T$ is a sharp of $x$, then this $T$ is unique. Therefore $x^\sharp$ will denote this unique $T$ for the real $x$. 
\end{definition}

Solovay showed that every subset of $\omega_1$ is constructible from a real under $\AD$:

\Begin{fact}{solovay subset of omega1 constructible}
(Solovay) Assume $\ZF + \AD$. There is an single formula $\theta(x,\alpha)$ so that for all $A \subseteq \omega_1$, there is an $x \in \reals$ so that $\alpha \in A \Leftrightarrow L[x] \models \theta(x,\alpha)$. 
\end{fact}

\begin{proof}
Let $A \subseteq \omega_1$. Consider the game
$$\begin{tikzpicture}
\node at (0,0) {$S_A$};

\node at (1,.5) {I};
\node at (2,.5) {$x(0)$};
\node at (4,.5) {$x(1)$};
\node at (6,.5) {$x(2)$};
\node at (8,.5) {$x(3)$};
\node at (11,.5) {$x$};

\node at (1,-.5) {II};
\node at (3,-.5) {$y(1)$};
\node at (5,-.5) {$y(2)$};
\node at (7,-.5) {$y(3)$};
\node at (9,-.5) {$y(4)$};
\node at (11,-.5) {$y$};
\end{tikzpicture}$$

Player 2 wins if the disjunction of the following holds

\noindent (1) $x \notin \WO$.

\noindent (2) $x \in \WO \Rightarrow (\forall n)((y_n \in \WO) \wedge (\exists \gamma \geq \ot(x))(\{\ot(y_n) : n \in \omega\} = A \cap \gamma))$. (Here $y_n$ is defined by $y_n(k) = y(\langle n,k\rangle)$.) 

Suppose Player 1 has a winning strategy $\sigma$. Then for all $y \in \reals$, $\sigma(y) \in \WO$. Then $E = \{v : (\exists y)(v = \sigma(y))\}$ is an $\analytic$ subset of $\WO$. Thus there is some $\gamma < \omega_1$ so that $\ot(v) < \gamma$ for all $v \in E$. Let $y \in \reals$ be such that for all $n$, $y_n \in \WO$ and $\{\ot(y_n) : n \in \omega\} = A \cap \gamma$. Then Player 1 using $\sigma$ loses if Player 2 plays this $y$. 

By $\AD$, Player 2 must have a winning strategy $\tau$. Note that $\alpha \in A$ if and only if $L[\tau]$ satisfies that $1_{\mathrm{Coll}(\omega,\alpha + 1)}$ forces that for all $w \in \WO_{\alpha + 1}$, there is some $n$ so that $\ot(\tau(w)_n) = \alpha$.  The formula $\theta$ is defined to assert this statement. The real associated to $A$ is of course $\tau$.
\end{proof}

\Begin{fact}{sharp code for functions}
Let $f : \omega_1 \rightarrow \omega_1$. There is some $x \in \reals$ and an $n \in \omega$ so that for all $\alpha < \omega_1$, $f(\alpha) = t_n^{\Gamma(x^\sharp,\beta)}(\alpha)$ whenever $\beta > \alpha$. 
\end{fact}

\begin{proof}
Consider $A = \{\langle \alpha,\beta\rangle : f(\alpha) = \beta\}$ where $\langle \cdot, \cdot \rangle$ denotes some constructible pairing function. By Fact \ref{solovay subset of omega1 constructible}, there is some real $x$ so that $\langle \alpha,\beta\rangle \in A$ if and only if $L[x] \models \theta(x,\langle\alpha,\beta\rangle)$.

Consider the formula $\psi(\beta,\alpha)$ if and only if $\theta(\dot E,\langle\alpha,\beta\rangle)$. There is some $n$ so that $t_n$ is the Skolem function associated to this formula. Fix two ordinals $\alpha < \beta$. Since $\Gamma(x^\sharp,\beta)$ is of the form $L_{\beta'}[x]$ for some $\beta' \geq \beta$ and $\Gamma(x^\sharp,\beta) \prec_{\omega} L[x]$, $t_n^{\Gamma(x^\sharp,\beta)}(\alpha) = f(\alpha)$.
\end{proof}

It should be noted that coding functions by simply coding their graph as a subset of $\omega_1$ is generally not good enough. For each $\alpha$, if one had to search for the corresponding $\beta$ which is the value of $f(\alpha)$, the coding system will be too complex. In this case, one uses the Skolem term to output the value of $f$ when given $\alpha$. This needs to be handled in the proof of Kechris as well. See Remark \ref{solovay function coding remark}. There it is handled by a modified version of the Solovay game used above.

\Begin{theorem}{omega1 is omega1 reasonable}
(Martin) Assume $\ZF + \AD$. $\omega_1$ is $\omega_1$-reasonable.
\end{theorem}

\begin{proof}
If $x \in \reals = \bairespace$, let $\cut(x) \in \bairespace$ be defined by $\cut(x)(k) = x(k + 1)$.

One will now define a good coding system for ${}^{\omega_1}\omega_1$. 

The pointclass of the coding system is $\analytic$. 

For each $x \in \reals$, let $\decode(x)(\alpha,\beta)$ holds if and only if $A(\cut(x)) \wedge t_{x(0)}^{\Gamma(\cut(x),\alpha + \omega)}(\alpha) = \beta$. Define $\GC_{\beta,\gamma}$ by $x \in \GC_{\beta,\gamma}$ if and only if $\decode(x)(\beta,\gamma)$ holds.

Suppose $f : \omega_1 \rightarrow \omega_1$. Let $z$ be the real and $n$ be the natural number obtained by applying Fact \ref{sharp code for functions} to $f$. Let $x \in \reals$ be such that $x(0) = n$ and $\cut(x) = z^\sharp$. Then $\decode(x) = f$. Such an $x$ will be called a sharp code for $f$.

Let $\beta < \omega_1$ and $\gamma < \omega_1$, note that $\GC_{\beta,\gamma}$ consists of those $x \in \reals$ such that the following holds 

\noindent (i) $A(\cut(x))$ 

\noindent (ii) $\beta$ is in the wellfounded part of $\Gamma(\cut(x),\beta + \omega)$

\noindent (iii) $t_{x(0)}^{\Gamma(\cut(x),\beta + \omega)}(\beta) = \gamma$.

By the discussion earlier, $\GC_{\beta,\gamma}$ is $\borel$. Note that since $t_{x(0)}$ is a Skolem function, $\tau_{x(0)}^{\Gamma(\cut(x),\beta + \omega)}$ is a function on its domain.

Now suppose $A \subseteq \GC_{\beta}$ is $\analytic$. Let $E = \{v \in \reals : (\exists x)(x \in A \wedge v = \mathsf{ST}(\Gamma(\cut(x),\beta + \omega),t_{x(0)}^{\Gamma(\cut(x),\beta + \omega)})\}$. Then $E$ is a $\analytic$ subset of $\WO$. By the boundedness lemma, there is some $\delta < \omega_1$ so that for all $v \in E$, $\ot(v) < \delta$. Thus $A \subseteq \bigcup_{\gamma < \delta} \GC_{\beta,\gamma}$.  

It has been shown that $(\analytic,\decode,\GC_{\beta,\gamma} : \beta < \omega_1,\gamma < \omega_1)$ is a good coding system for ${}^{\omega_1}\omega_1$. 
\end{proof}

Next, one will give the argument of Kechris that $\omega_1$ is a strong partition cardinal. This proof uses category and the Kechris-Woodin generic coding idea. However at $\omega_1$, the generic coding function is essentially trivial and exists without even $\AD$. For many other purposes, the generic coding function is very useful.

Note that in Martin's proof, the indiscernibility models are used to make complexity computations. In the argument of Kechris, category quantifiers will be used to ensure sets have the correct complexity.

First some results on the Banach-Mazur games. See \cite{Classical-Descriptive-Set-Theory} Section 21.C and 21.D for the proofs of these results:

\Begin{fact}{banach-mazur games}
(Banach-Mazur Game) Assume $\ZF$. Let $A \subseteq \bairespace$. Define $G^*(A)$ by 
$$\begin{tikzpicture}
\node at (0,0) {$G^*_A$};

\node at (1,.5) {I};
\node at (2,.5) {$s_0$};
\node at (4,.5) {$s_2$};
\node at (6,.5) {$s_4$};

\node at (1,-.5) {II};
\node at (3,-.5) {$s_1$};
\node at (5,-.5) {$s_3$};
\node at (7,-.5) {$s_5$};
\end{tikzpicture}$$
Player 1 and 2 alternatingly play nonempty strings in $\finNaturalSequence$. Let $x = s_0\hat{\ }s_1\hat{\ }s_2 ...$ be the concatenation of the moves. Player 2 wins $G_A^*$ if and only if $x \in A$. 

(1) $A$ is comeager if and only if Player 2 has a winning strategy in $G^*_A$.

(2) There is some $s \in \finNaturalSequence$ so that $\bairespace\setminus A$ is comeager in $N_s$ if and only if Player 1 has a winning strategy in $G_A^*$.
\end{fact}

\Begin{fact}{unfolded banach-mazur games}
(Unfolded Banach-Mazur Game) Assume $\ZF$. Let $A \subseteq \bairespace$ and $B \subseteq \bairespace \times \bairespace$ be such that for all $x \in A$, there exists a $y \in \bairespace$ so that $B(x,y)$. Let $G^*_{A,B}$ be the following game
$$\begin{tikzpicture}
\node at (0,0) {$G^*_{A,B}$};

\node at (1,.5) {I};
\node at (2,.5) {$s_0$};
\node at (4,.5) {$s_2$};
\node at (6,.5) {$s_4$};

\node at (1,-.5) {II};
\node at (3,-.5) {$s_1,y(0)$};
\node at (5,-.5) {$s_3,y(1)$};
\node at (7,-.5) {$s_5,y(2)$};
\end{tikzpicture}$$
Player 1 plays nonempty strings in $\finNaturalSequence$. Player 2 plays nonempty strings in $\finNaturalSequence$ and an element of $\omega$. Let $x = s_0\hat{\ }s_1\hat{\ }s_2...$ be the concatenation of the strings played by both player. Let $y$ be the real produced by the extra elements of $\omega$ played by Player 2. Player 2 wins $G^*_{A,B}$ if and only if $x \in A$ and $(x,y) \in B$. 

If Player 2 wins $G_{A,B}^*$, then Player 2 wins $G^*_A$. If Player 1 wins $G_{A,B}^*$, then Player 1 wins $G^*_A$. 
\end{fact}

\Begin{corollary}{comeagerness and unfolded banach-mazur}
Assume $\ZF$. Let $A \subseteq \bairespace$ be a $\analytic$ set. Let $B \subseteq \bairespace\times\bairespace$ be a $\bPi_1^0$ set so that $A = \pi_1[B]$. $A$ is comeager if and only if Player 2 has a winning strategy in $G^*_{A,B}$.
\end{corollary}

\begin{proof}
Note that if $A = \pi_1[B]$, then the winning condition for Player 2 in $G^*_{A,B}$ is equivalent to merely $(x,y) \in B$. Thus in this case $G^*_{A,B}$ is a closed game. Thus the determinacy of such games hold under $\ZF$. Now the result follows from Fact \ref{banach-mazur games} and \ref{unfolded banach-mazur games}.
\end{proof}

If $\alpha < \omega_1$, one can define a topology on ${}^\omega\alpha$ be declaring the basic open sets to be sets of the form $N_s^\alpha = \{f \in {}^\omega\alpha : f \supseteq s\}$, where $s \in {}^{<\omega}\alpha$. Using this topology, one can define the notions of comeagerness and meagerness for ${}^\omega\alpha$. As $\alpha$ is countable, ${}^\omega\alpha$ with this topology is homeomorphic to $\bairespace$. Observe that the set $\mathrm{surj}_\alpha = \{f \in {}^\omega\alpha : f \text{ is a surjection}\}$ is a comeager subset of ${}^\omega\alpha$. If $A \subseteq \reals \times {}^\omega\alpha$, then one writes $(\forall^*_\alpha f)A(x,f)$ as an abbreviation for the statement that $A(x,f)$ holds for comeagerly many $f$ in ${}^\omega\alpha$.

\Begin{corollary}{complexity category quantifier}
Assume $\ZF$. Let $A \subseteq \bairespace \times \bairespace$ be a $\analytic$ set. Then $A_0(x) \Leftrightarrow (\forall_\omega^* y)A(x,y)$ is $\analytic$ and $A_1(x) \Leftrightarrow (\forall^*_\omega y)\neg A(x,y)$ is $\coanalytic$.
\end{corollary}

\begin{proof}
Let $B \subseteq \bairespace \times \bairespace \times \bairespace$ be such that $A = \{(x,y) : (\exists z)((x,y,z) \in B)\}$ where $B \in \bPi_1^0$. For each $x \in \bairespace$, let $B_x = \{(y,z) : (x,y,z) \in B\}$ and $A_x = \{y : (x,y) \in A\}$. Note that $B_x \in \bPi_1^0$ and $\pi_1[B_x] = A_x$. 

If $\sigma$ is a Player 1 strategy and $\tau$ is a Player 2 strategy for a game of the form $G^*_{C,D}$ for some appropriate $C$ and $D$, then let $x_{\sigma * \tau}$ be the real produced by the concatenation of the finite strings played by each player. Let $y_{\sigma * \tau}$ be the auxiliary sequence produced by the moves of $\tau$.

Note that $A_0(x)$ if and only if Player 2 has a winning strategy in $G^*_{A_x,B_x}$ if and only if $(\exists \tau)(\forall \sigma)((x_{\sigma * \tau}, y_{\sigma*\tau}) \in B_x)$ by Fact \ref{unfolded banach-mazur games}. Since $B_x \in \bPi_1^0$ and $\bPi_1^0$ is closed under $\forall^\reals$, the latter expression is $\analytic$. 

For $s \in \finNaturalSequence$, let $\Phi_s : N_s \rightarrow \bairespace$ be the canonical homeomorphism between $N_s$ and $\bairespace$. Let $A^s \subseteq \bairespace \times \bairespace$ be defined by $(x,y) \in A^s \Leftrightarrow (x,\Phi_s^{-1}(y)) \in A$. Let $B^s \subseteq \bairespace \times \bairespace \times \bairespace$ be defined by $(x,y,z) \in B^s$ if and only if $(x,\Phi_s^{-1}(y),z) \in B$. 

For the second statement, note that $\neg A_1(x)$ if and only if $\neg((\forall^*_\omega y)\neg A(x,y))$. Since $\coanalytic$ sets have Baire property, the latter holds if and only if there exists a $s \in \finBinarySequence$ so that $A_x$ is comeager in $N_s$. Now there exists an $s \in \finBinarySequence$ so that $A_x$ is comeager in $N_s$ if and only if there exists an $s \in \finBinarySequence$ so that $\Phi_s[A_x \cap N_s]$ is comeager if and only if there exists an $s \in \finBinarySequence$ so that Player 2 has a winning strategy in $G^*_{(A^s)_x,(B^s)_x}$. As before, the last expression can be checked to be $\analytic$. Hence $A_1$ is $\coanalytic$.
\end{proof}

The following is the (essentially trivial) version of the Kechris-Woodin generic coding function for $\omega_1$. (See \cite{Generic-Codes-for-Uncountable-Ordinals}.)

\Begin{fact}{generic coding for omega1}
There is a continuous function $G : {}^{\omega}\omega_1 \rightarrow \WO$ so that for all $\ell \in {}^\omega\omega_1$ such that $\ell(0) = \{\ell(n + 1) : n \in \omega\}$, $G(\ell) \in \WO$ and $\ot(G(\ell)) = \ell(0)$.
\end{fact}

\begin{proof}
For each $\ell$ as above, let $A_\ell = \{n \in \omega\setminus \{0\} : (\forall m)(\ell(n) = \ell(m) \Rightarrow n \leq m)\}$. Let $G(\ell)$ be the wellordering with domain $A_\ell$ so that $m <_{G(\ell)} n$ if and only if $\ell(m) < \ell(n)$. $G(\ell) \in \WO$ and $\ot(G(\ell)) = \ell(0)$. 

This function is continuous in the sense that for any $n$, $G(\ell)\upharpoonright n$ is determined by $\ell \upharpoonright m$ for some $m$.
\end{proof}

\Begin{fact}{solovay game for coding functions}
Assume $\ZF + \AD$. Let $f : \omega_1 \rightarrow \omega_1$. There is a Player 2 strategy so that for all $v \in \WO$, $\{(\alpha,\beta) : (\exists n)((\tau(v)_n)_0 \in \WO_\alpha \wedge (\tau(v)_n)_1 \in \WO_\beta)\} = f \upharpoonright \gamma$ for some $\gamma > \ot(v)$. 

Here $\tau(v)$ is $(v * \tau)_\mathrm{odd}$, i.e. the real produced by Player 2 when played against Player 1 playing the bits of $v$. Also $f \upharpoonright \gamma = \{(\alpha,\beta) \in f : \alpha < \gamma\}$. Recall that if $x \in \reals$, $x_n$ is defined by $x_n(k) = x(\langle n,k\rangle)$. This is the $n^\text{th}$ section of $x$. This result states that when given $v \in \WO$, the integer sections of $\tau(v)$ codes $f$ up to some ordinal $\gamma$ greater than $\ot(v)$. 
\end{fact}

\begin{proof}
Consider the game $S_f$ defined by
$$\begin{tikzpicture}
\node at (0,0) {$S_f$};

\node at (1,.5) {I};
\node at (2,.5) {$v(0)$};
\node at (4,.5) {$v(1)$};
\node at (6,.5) {$v(2)$};
\node at (8,.5) {$v(3)$};
\node at (11,.5) {$v$};

\node at (1,-.5) {II};
\node at (3,-.5) {$r(0)$};
\node at (5,-.5) {$r(1)$};
\node at (7,-.5) {$r(2)$};
\node at (9,-.5) {$r(3)$};
\node at (11,-.5) {$r$};
\end{tikzpicture}$$
Player 2 wins if and only if $v \in \WO$ implies that $\{(\alpha,\beta) : (\exists n)((r_n)_0 \in \WO_\alpha \wedge (r_n)_1 \in \WO_\beta)\} = f \upharpoonright \gamma$ for some $\gamma > \ot(v)$.

By essentially the same bounding argument of Fact \ref{solovay subset of omega1 constructible}, one has the Player 2 must have the winning strategy in this game.
\end{proof}

\Begin{remark}{solovay function coding remark}
In the above statement, it is very important that some section of $\tau(v)$ contains a code for the image of $f(\alpha)$. To search for $f(\alpha)$, which could be quite large in comparison to $\alpha$, would push the complexity of $\GC_{\alpha,\beta}$ beyond $\borel$. Instead, a particular winning strategy $\tau$ will take a code $v$ for $\alpha$ and ouput $\tau(v)$ which magically contains codes for $f(\alpha)$ among its integer sections, $\{(\tau(v)_n)_1 : n \in \omega\}$.
\end{remark}

\Begin{theorem}{kechris reasonableness proof}
(Kechris' proof) Assume $\ZF + \AD$. $\omega_1$ is $\omega_1$-reasonable.
\end{theorem}

\begin{proof}
One will define a good coding system for ${}^{\omega_1}\omega_1$ witnessing $\omega_1$-reasonableness of $\omega_1$. 

The associated pointclass is $\analytic$. 

Let $\add : \omega_1 \times {}^\omega\omega_1 \rightarrow {}^\omega\omega_1$ be defined by 
$$\add(\alpha,g)(n) = \begin{cases}
\alpha & \quad n = 0 \\
g(n - 1) & \quad n > 0
\end{cases}$$
$\add(\alpha,g)$ simply inserts $\alpha$ at the beginning of $g$.

In some fix recursive manner, every real $x$ codes a Player 2 strategy $\tau_x$.

Now fix $\alpha,\beta \in \omega_1$. Define the formula $\varphi^\alpha_{\beta}(x,z)$ to state that 
$$(\exists n)([(\tau_x(z)_n)_0 \in \WO_\alpha \wedge (\tau_x(z)_n)_1 \in \WO_\beta] \wedge (\forall m < n)((\tau_x(z)_m)_0 \notin \WO_\alpha)).$$ 
This formula defines a $\borel$ relation using elements of $\WO_\alpha$ and $\WO_\beta$ as parameters.

Let $\phi_\beta^\alpha(x)$ be the statement $(\forall^*_\alpha g)\varphi^\alpha_\beta(x,G(\add(\alpha,g)))$, where $G$ is the generic coding function of Fact \ref{generic coding for omega1}. Note that since $\mathrm{surj}_\alpha$ is comeager in ${}^\omega\alpha$, $(\forall^*_\alpha g)(G(\add(\alpha,g)) \in \WO_\alpha)$. 

Now fix a bijection $B : \omega \rightarrow \alpha$. Let $T : \bairespace \rightarrow {}^\omega\alpha$ be defined by $T(r) = B\circ r$. Define $\tilde G : \bairespace \rightarrow \WO$ by $\tilde G(r) = G(\add(\alpha,T(r)))$. Now $\tilde G$ is a continuous function with the property that for comeagerly many $r \in \bairespace$ (in the usual sense of comeager), $\tilde G(r) \in \WO_\alpha$.

$\phi^\alpha_\beta(x)$ is equivalent to $(\forall^*_\omega r)\varphi^\alpha_\beta(x,\tilde G(r))$. Since $\varphi^\alpha_\beta$ defines a $\borel$ relation (in parameters from $\WO_\alpha$ and $\WO_\beta$), the collection of $(x,r)$ satisfying $\varphi^\alpha_\beta(x,\tilde G(r))$ is $\borel$ using parameters from $\WO_\alpha$, $\WO_\beta$, and a parameter coding the continuous function $\tilde G$. Now the set defined by $\phi^\alpha_\beta(x)$ is $\borel$ by Corollary \ref{complexity category quantifier}. (Note this is not done uniformly in the $\alpha$ and $\beta$.)

If $x \in \reals$, define $\decode(x)(\alpha,\beta)$ if and only if $\phi_\beta^\alpha(x)$. 

Define $\GC_{\alpha,\beta} = \{x \in \reals : \phi_\beta^\alpha(x)\}$. $\GC_{\alpha,\beta}$ is $\borel$.

For any $f : \omega_1 \rightarrow \omega_1$, Fact \ref{solovay game for coding functions} states that there is some $x$ so that $\decode(x) = f$. 

Fix $\alpha < \omega_1$. For each $x \in \bairespace$, let $\psi_0(x,z)$ assert that $(\exists n)((\tau_x(z)_n)_0 \in \WO_\alpha)$. $\psi_0(x,z)$ is $\borel$ using any code of $\alpha$ as a parameter. Let $\psi_1(x,y,z)$ be the conjunction of the following statements

\noindent (1) $\psi_0(x,z) \wedge \psi_0(y,z)$.

\noindent (2) There exists $w_0,w_1 \in \reals$, $n_0,n_1 \in \omega$ so that 

(2a) $(\tau_x(z)_{n_0})_0 \in \WO_\alpha$, $(\tau_x(z)_{n_0})_1 = w_0$, and for all $m < n_0$, $(\tau_x(z)_m)_0 \notin \WO_\alpha$.

(2b) $(\tau_y(z)_{n_0})_0 \in \WO_\alpha$, $(\tau_y(z)_{n_1})_1 = w_1$, and for all $m < n_1$, $(\tau_y(z)_m)_0 \notin \WO_\alpha$.

(2c) $w_0 <_{\Sigma_1^1} w_1$, where $<_{\Sigma_1^1}$ is the $\Sigma_1^1$ relation witnessing that the ordertype function is a $\Pi_1^1$-norm.
\newline
\newline
Observe that $\psi_1(x,y,z)$ is $\analytic$.

Let $A \subseteq \GC_\alpha$ and $A$ is $\analytic$. Define a relation $\prec$ on $\bairespace$ by $x \prec y$ if and only if 
$$x \in A \wedge y \in A \wedge (\forall^*_\alpha g)\psi_1(x,y,G(\add(\alpha,g))).$$
By Corollary \ref{complexity category quantifier}, $\prec$ is a $\analytic$ relation on $\bairespace$. 

By checking the various definitions, one can see that for all $x,y \in A$, $\decode(x)(\alpha) < \decode(y)(\alpha)$ if and only if $x \prec y$.

Since $\prec$ is a $\analytic$ relation, it is an $\omega$-Suslin relation. By the Kunen-Martin theorem (\cite{Notes-on-the-Theory-of-Scales} Section 7), the rank must be less than $\omega_1$. Since $\omega_1$ is regular, there is some $\gamma < \omega_1$ so that $A \subseteq \bigcup_{\beta < \gamma} \GC_\beta$. 

(Some presentations of the Kunen-Martin theorem uses $\DC_\reals$, but with some care, one can remove the use of $\DC_\reals$. Another approach in the case of $\omega$-Suslin relations  is to observe that if $K$ denotes a tree witnessing a relation is $\omega$-Suslin, one can prove the Kunen-Martin theorem in $L[K] \models \AC$ and argue in $L[K]$ that the Kunen-Martin tree is wellfounded there and hence in the real world. Then one shows that this tree still works in the real world $V$. Alternatively, if one assumes Kechris's result that $L(\reals) \models \AD$ implies $L(\reals) \models \DC_\reals$, then one can absorb this problem into $L(\reals)$ and apply the Kunen-Martin Theorem in $L(\reals)$.)

It has been shown that $(\analytic, \decode,\GC_{\beta,\gamma} : \beta < \omega_1, \gamma < \omega_1)$ is a good coding system for ${}^{\omega_1}\omega_1$. 
\end{proof}

\Begin{corollary}{omega1 is strong partition}
(Martin) Assume $\ZF + \AD$. $\omega_1 \rightarrow_* (\omega_1)^{\omega_1}_2$. 
\end{corollary}

\Begin{definition}{club partition measure}
Let $\lambda \leq \omega_1$. Let $\mu_{\omega_1}^\lambda$ consists of those subsets $X$ of $[\omega_1]^{\lambda}$ which contain a set of the form $[C]^{\lambda}_*$ where $C$ is a club subset of $\omega_1$. Each $\mu_{\omega_1}^\lambda$ is a countably complete ultrafilter.
\end{definition}

\Begin{corollary}{absorbing function on omega1 into L(R)}
Assume $\ZF + \AD$. Let $\lambda \leq \omega_1$. $\prod_{[\omega_1]^{\lambda}}\omega_1 \slash \mu_{\omega_1}^\lambda = (\prod_{[\omega_1]^{\lambda}}\omega_1 \slash \mu_{\omega_1}^\lambda)^{L(\reals)}$.

Since $L(\reals) \models \DC$ by a result of Kechris \cite{The-Axiom-of-Determinacy-Implies-Dependent-Choice}, $\prod_{[\omega_1]^{\lambda}}\omega_1 \slash \mu_{\omega_1}^\lambda$ is an ordinal.
\end{corollary}

\begin{proof}
Note that the good coding system for ${}^\lambda\omega_1$ constructed above belong to $L(\reals)$. Now apply Theorem \ref{absorbing function on function spaces}.
\end{proof}

\Begin{remark}{final good coding remark}
This section presented two proofs of the strong partition property for $\omega_1$. The original proof of Martin uses indiscernibility. Many early results were proved using this indiscernibility idea. For example, Kunen showed in $\AD + \DC$ that $\delta_3^1$ is a weak partition cardinal and $\delta_3^1 = \aleph_{\omega + 1}$. 

The second proof of Kechris is perhaps the simplest proof using classical descriptive set theoretic ideas. The ideas of the Kechris-Woodin generic coding function is very useful in various different settings. 

However, no proof in the flavor of these two arguments is known to establish the strong partition property at $\delta_3^1$. There is another proof of the strong partition property for $\omega_1$ due to Jackson. It uses ideas such as the Kunen tree and an analysis of all the measures on $\omega_1$ due to Kunen. An exposition can be found in \cite{A-New-Proof-Strong-Partition-Relation}, \cite{Structural-Consequences-AD}, and \cite{A-Computation-of-delta15}. More generally, Jackson developed the theory of description. This theory produces the only known proof that $\delta_{3}^1$ and in fact all $\delta_{2n + 1}^1$ are strong partition cardinal under $\AD + \DC$. Jackson has also computed the identity of these cardinals. For example, $\delta_5^1 = \aleph_{{\omega^{\omega^\omega}} + 1}$ and in fact there is a general formula for computing $\delta_{2n + 1}^1$. See \cite{A-Computation-of-delta15} for more information. It should noted that Jackson's proof that $\omega_1$ is a strong partition cardinal, Kunen's result that $\delta_3^1$ is a weak partition cardinal, and all of Jackson's results mentioned above using description theory are proved using $\DC$. The original Martin proof and the Kechris proof of the strong partition property of $\omega_1$ are within $\ZF + \AD$.
\end{remark}

\section{Kunen Functions and Partition Properties at $\omega_2$}\label{kunen functions and partition properties at omega2}

\Begin{fact}{normal measure and closure points}
Let $\mu$ be a normal measure on a cardinal $\kappa$. Let $A \in \mu$, let $\enum_A : \omega_1 \rightarrow A$ be the increasing enumeration of $A$. Then $\{\alpha \in A : \enum_A(\alpha) = \alpha\} \in \mu$. 

Let $f : \kappa \rightarrow \kappa$ and $X \subseteq f[\kappa]$ so that $\{\alpha < \kappa : f(\alpha) \in X\} \in \mu$, then $\enum_X =_\mu f$. 
\end{fact}

\begin{proof}
Suppose not. Then $B = \{\alpha \in A : \alpha < \enum_A(\alpha)\} \in \mu$. Since $\enum_A$ is an order-preserving function, $B = \{\alpha \in A : \enum_A^{-1}(\alpha) < \alpha\}$. By normality, there is some $C \subseteq B$ and $\gamma < \kappa$ so that for all $\alpha \in C$, $\enum_A^{-1}(\alpha) = \gamma$. This contradicts the injectiveness of $\enum_A$.

For the second statement: Let $A = \{\alpha : f(\alpha) \in X\}$. Note that $f \circ \enum_A = \enum_X$. Let $B \subseteq A$ with $B \in \mu$ be such that $\enum_A(\alpha) = \alpha$ for all $\alpha \in B$. Then for all $\gamma \in B$, $\enum_X(\alpha) = f(\enum_A(\alpha)) = f(\alpha)$. Hence $\enum_X =_\mu f$.
\end{proof}

The next result states that an ultrapower of a strong partition cardinal by a measure is a cardinal as long as that ultrapower is well founded.

\Begin{fact}{wf ultrapower by measure is cardinal}
(Martin) Let $\kappa$ be a strong partition cardinal and $\mu$ be a measure on $\kappa$. If ${}^\kappa\kappa\slash\mu$ is wellfounded, then ${}^\kappa\kappa\slash \mu$ is a cardinal.

Assuming $\ZF + \AD + \DC_\reals$. If $\kappa < \Theta$ is a strong partition cardinal and $\mu$ is a measure on $\kappa$, then ${}^\kappa\kappa\slash\mu$ is a cardinal.
\end{fact}

\begin{proof}
Suppose there is an $h \in {}^\kappa\kappa\slash\mu$ and an injection $\Phi : {}^\kappa\kappa\slash\mu \rightarrow \prod_{\kappa}h \slash \mu$, where $\prod_{\kappa}h \slash \mu$ is the collection of $[f]_\mu$ such that $f <_\mu h$. Note that $\prod_{\kappa} h \slash \mu$ is the initial segment of ${}^\kappa\kappa\slash\mu$ determined by $[h]_\mu$.

Let $\mathcal{T} = (\kappa \times 2, \sqsubset)$ denote the lexicographic ordering on $\kappa\times 2$. If $F : \kappa \times 2 \rightarrow \kappa$ is a $\mathcal{T}$-increasing function, then let $f_0, f_1: \kappa \rightarrow \kappa$ be defined by $f_0(\alpha) = F(\alpha,0)$ and $f_1(\alpha) = F(\alpha,1)$. Note that $f_0$ and $f_1$ are increasing functions.

Define $P : [\kappa]^\mathcal{T} \rightarrow 2$ by 
$$P(F) =  \begin{cases}
0 & \quad \Phi([f_0]_\mu) < \Phi([f_1]_\mu) \\
1 & \quad \text{otherwise}
\end{cases}$$
Since $\mathrm{ot}(\mathcal{T}) = \kappa$, $\kappa \rightarrow (\kappa)^\kappa_2$ implies that there is an $E \subseteq \kappa$ with $|E| = \kappa$ which is homogeneous for this partition.

Suppose $P$ is homogeneous taking value $1$. For each $n \in \omega$, define $g_n : \kappa \rightarrow \kappa$ by $g_n(\alpha) = \enum_E(\omega\cdot\alpha + n)$. Note that if $m < n$, then $g_m <_\mu g_n$. For any $m<n$, define $G^{m,n} : \mathcal{T} \rightarrow \kappa$ by $G^{m,n}(\alpha,0) = g_m(\alpha)$ and $G^{m,n}(\alpha,1) = g_n(\alpha)$. Note that $G^{m,n}$ is $\mathcal{T}$-increasing and $G^{m,n} \in [\kappa]^{\mathcal{T}}$. Since $P(G^{m,n}) = 1$, one has that $\Phi([g_n]_\mu) = \Phi([G^{m,n}_1]_{\mu}) \leq \Phi([G^{m,n}_0]_\mu) = \Phi([g_m]_\mu)$. Since $\Phi$ is an injection, one must actually have $\Phi([g_n]_\mu) < \Phi([g_m]_\mu)$. Thus $\langle \Phi([g_k]_\mu) : k \in \omega\rangle$ is an infinite decreasing sequence in ${}^\kappa\kappa\slash\mu$. This contradicts the assumption that ${}^\kappa\kappa\slash\mu$ is wellfounded.

Thus $E$ must be homogeneous for $P$ taking value $0$.

Let $\ell \in {}^\kappa \kappa$ be such that $[\ell]_\mu > [h]_\mu$ and $\ell(\alpha) > 0$ for all $\alpha \in \kappa$. Let $U = \{(\alpha,\beta) \in \kappa \times \kappa : \beta < \ell(\alpha)\}$. Let $\mathcal{U}$ denote $(U,\prec)$ where $\prec$ denotes the lexicographic ordering of $U$. Again, $\mathrm{ot}(\mathcal{U}) = \kappa$. 

Let $F = \{\enum_E(\omega\cdot\alpha) : \alpha < \kappa)\}$. Note that $|F| = \kappa$. Pick some $K : \mathcal{U} \rightarrow F$ which is order-preserving. For any $f <_\mu \ell$, define $k_f : \kappa \rightarrow \kappa$ by
$$k_f(\alpha) = \begin{cases}
K(\alpha,f(\alpha)) & \quad f(\alpha) < \ell(\alpha) \\
K(\alpha,0) & \quad \text{otherwise}
\end{cases}.$$
Note that if $f <_\mu g <_\mu \ell$, then $k_f <_\mu k_g$ and that for all $f <_\mu \ell$, $k_f \in [F]^\kappa$.

Define $\Psi : \prod_{\kappa}\ell \slash \mu \rightarrow \prod_{\kappa} h \slash \mu$ as follows: Let $A < [\ell]_\mu$. Pick $f \in A$, i.e. a representative of $A$. Let $\Psi(A) = \Phi([k_f]_{\mu})$. One can check that $\Psi$ is well defined. 

The next claim is that $\Psi$ is order preserving: Suppose $A < B < [\ell]_\mu$. Let $f \in A$ and $g \in B$ be representatives for $A$ and $B$, respectively. Since $k_f <_\mu k_g$, the set $C = \{\alpha \in \kappa : k_g(\alpha) \leq k_f(\alpha)\} \notin \mu$.  Define $k_g' : \kappa \rightarrow E$ by $k_g'(\alpha) = k_g(\alpha)$ if $\alpha \notin C$ and $k_g'(\alpha)$ is the next element of $E$ greater than $k_f(\alpha)$ if $\alpha \in C$. The purpose of defining $F$ in the manner above was to ensure that between two successive elements of $F$, there are at least $\omega$ many point of $E$ in between. From this, one can verify that $k_g' \in [E]^\kappa$ and $k_f(\alpha) < k_g'(\alpha) < k_f(\alpha + 1)$ for all $\alpha < \kappa$. Note that $[k_g']_\mu = [k_g]_\mu$ since $k_g$ and $k_g'$ agree off of $C \notin \mu$.

Let $F^{k_f,k_g'} : \mathcal{T} \rightarrow \kappa$ be defined by
$$F^{k_f,k_g'}(\alpha,i) = \begin{cases}
k_f(\alpha) & \quad i = 0\\
k_g'(\alpha) & \quad i = 1
\end{cases}$$
Again using the main property of $F$, one can verify that $F^{k_f,k'_g} \in [\kappa]^{\mathcal{T}}$. Thus $P(F^{k_f,k_g'}) = 0$. Thus $\Psi(A) = \Phi([k_f]_\mu) < \Phi([k_g']_\mu) = \Phi([k_g]_\mu) = \Psi(B)$. 

It has been shown that $\Psi$ is order preserving. However, this is not possible since $[h]_\mu < [\ell]_\mu$. This completes the proof.
\end{proof}

\Begin{fact}{presliding lemma}
Let $\kappa$ be a measurable cardinal possessing a normal $\kappa$-complete measure $\mu$. Let $\epsilon < \kappa$. Let $\mathcal{T}^\epsilon = (\kappa \times \epsilon, \sqsubset)$ where $\sqsubset$ is the lexicographic ordering. For $F : \mathcal{T}^\epsilon \rightarrow \kappa$ which is order preserving and $\alpha < \epsilon$, let $F^\alpha \in [\kappa]^\kappa$ be defined by $F^\alpha(\gamma) = F(\gamma,\alpha)$. 

Let $\langle f_\alpha : \alpha < \epsilon\rangle$ be a sequence in $[\kappa]^\kappa$ such that for all $\alpha < \beta < \epsilon$, $f_\alpha <_\mu f_\beta$. Then there is an $F \in [\kappa]^{\mathcal{T}^\epsilon}$ so that for all $\alpha < \epsilon$, $F^\alpha =_\mu f_\alpha$.

Moreover, if $D \subseteq \kappa$ with $|D| = \kappa$ and each $f_\alpha \in [D]^\kappa$, then one can find $F \in [D]^{\mathcal{T}^\epsilon}$ with the above property.
\end{fact}

\begin{proof}
Suppose $F(\alpha',\beta')$ has been defined for all $(\alpha',\beta') \sqsubset (\alpha,\beta)$. Then let $F(\alpha,\beta)$ be the least element in the range of $f_\beta$ greater than $F(\alpha',\beta')$ for all $(\alpha',\beta') \sqsubset (\alpha,\beta)$.

The claim is that $F$ has the desired properties: This is proved by induction on $\beta < \epsilon$.

Consider $\beta = 0$. Suppose $A = \{\alpha < \kappa : F(\alpha,0) \neq f_0(\alpha)\} \in \mu$. This means for each $\alpha \in A$, there is some $\alpha' < \alpha$, and some $\beta < \epsilon$ so that $F(\alpha',\beta) \geq f_0(\alpha)$. Let $g : \kappa \rightarrow \kappa$ be defined by $g(\alpha)$ is the least such $\alpha'$ with the property above if $\alpha \in A$ and $0$ if $\alpha \notin A$. Note that $g(\alpha) < \alpha$ for all $\alpha \in A$. By the normality of $\mu$, there is some $A' \subseteq A$ with $A' \in \mu$ and some $\delta$ so that $g(\alpha) =\delta$ for all $\alpha \in A'$. This means for all $\alpha \in A' \in \mu$, $f_0(\alpha) \leq \sup\{F(\delta,\beta) : \beta < \epsilon\}$. By the $\kappa$-completeness of $\mu$, one may even find a $\beta^*$ and a $A'' \in \mu$ so that for all $\alpha \in A''$, $f_0(\alpha) \leq F(\delta,\beta^*)$. This is clearly impossible. This shows that $F^0 =_\mu f_0$. 

Suppose $\beta < \epsilon$ is such that for all $\beta' < \beta$, $F^{\beta'} = f_{\beta'}$. Suppose $\neg(F^\beta =_\mu f_\beta)$. Then $A = \{\alpha : F^\beta(\alpha) \neq f_\beta(\alpha)\} \in \mu$. Since $f_0(\alpha) = F^0(\alpha)$ for almost all $\alpha$, it can not be the case that for almost all $\alpha \in A$, there is some $\alpha' < \alpha$ and some $\beta' < \epsilon$ so that $F^{\beta'}(\alpha') \geq f_\beta(\alpha)$. Thus for almost all $\alpha \in A$, there is some $\beta' < \beta$ so that $F^{\beta'}(\alpha) \geq f_\beta(\alpha)$. By the $\kappa$-additivity of $\mu$, there is some $\beta^* < \beta$ so that for almost all $\alpha$, $F^{\beta^*}(\alpha) \geq f_\beta(\alpha)$. By the induction hypothesis, one has that $F^{\beta^*} =_\mu f_{\beta^*}$. One has shown that $f_\beta \leq_\mu f_{\beta^*}$ despite the fact that $\beta^* < \beta$. Contradiction. The result has been established.
\end{proof}

\Begin{fact}{ultrapower normal strong partition is regular}
(Martin) Let $\kappa$ be a strong partition cardinal and $\mu$ is a normal $\kappa$-complete measure on a cardinal $\kappa$. If ${}^\kappa\kappa \slash \mu$ is wellfounded, then ${}^\kappa\kappa\slash\mu$ is a regular cardinal.

Assuming $\ZF + \AD + \DC_\reals$, if $\kappa < \Theta$ is a strong partition cardinal and $\mu$ is a normal $\kappa$-complete measure on $\kappa$, then ${}^\kappa\kappa\slash\mu$ is a regular cardinal.
\end{fact}

\begin{proof}
By Fact \ref{wf ultrapower by measure is cardinal}, $\delta = \mathrm{ot}({}^\kappa\kappa\slash\mu)$ is a cardinal. Assume that ${}^\kappa\kappa \slash \mu$ is not a regular cardinal. There is some $\gamma < \delta$ and an increasing function $\Phi : \gamma \rightarrow \delta$ which is cofinal.

Let $\mathcal{T}^2$, $F^0$, $F^1$ come from the notation of Fact \ref{presliding lemma}. Define a partition $P : [\kappa]^{\mathcal{T}^2} \rightarrow 2$ by
$$P(F) = \begin{cases}
0 & \quad (\exists \beta < \gamma)([F^0]_\mu < \Phi(\beta) < [F^1]_\mu) \\
1 & \quad \text{otherwise}
\end{cases}$$
Let $D \subseteq \kappa$ with $|D| = \kappa$ be a homogeneous set for $P$.

(Case I) Assume $D$ is homogeneous for $P$ taking value $0$. 

By Fact \ref{normal measure and closure points}, every nonconstant function $f \in {}^\kappa D$ has some $g \in [D]^\kappa$ so that $f =_\mu g$. Thus one can show that $[D]^\kappa\slash \mu$ has ordertype $\delta$. 

For each $\alpha$, let $f_\alpha$ and $f_{\alpha + 1}$ be two elements of $[D]^\kappa$ which represents the elements of $[D]^\kappa\slash\mu$ of rank $\alpha$ and $\alpha + 1$ respectively. By Fact \ref{presliding lemma}, there is some $F : [D]^{\mathcal{T}^2} \rightarrow \kappa$ so that $F^0 = f_\alpha$ and $F^1 = F_{\alpha + 1}$. Then $P(F) = 0$. Hence let $\nu_\alpha$ be the least ordinal less than $\gamma$ so that $[f_\alpha] = [F^\alpha]_\mu < \Phi(\nu_\alpha) < [F^{\alpha + 1}]_\mu = [f_{\alpha + 1}]_\mu$. One can check that $\nu_\alpha$ depends only on $\alpha$ and not on the choice of $f_{\alpha}$ and $f_{\alpha + 1}$. Also if $\alpha \neq \alpha'$, then $\nu_{\alpha} < \nu_{\alpha'}$. This gives an order preserving map from $\delta$ into $\gamma < \delta$ which is impossible. 

(Case II) Assume $D$ is homogeneous for $P$ taking value $1$. 

As argued above, $[D]^\kappa\slash\mu$ has ordertype $\delta$ so $[D]^\kappa \slash \mu$ is a cofinal subset of $\kappa^\kappa\slash\mu$. Fix $f^* \in [D]^\kappa\slash\mu$. Since $\Phi$ is a cofinal map, there is some $\gamma^* < \gamma$ so that $[f^*]_\mu < \Phi(\gamma^*)$. Then there is some $f' \in [D]^\kappa$ so that $[f^*]_\mu < \Phi(\gamma^*) < [f']_{\mu}$. By Fact \ref{presliding lemma}, there is some $F \in [D]^{\mathcal{T}^2}$ so that $F^0 =_\mu f^*$ and $F^1 =_\mu f'$. Thus $P(F) = 0$. This is impossible since $D$ is homogeneous for $P$ taking value $1$. 

The failure of both cases would imply $P$ has no homogeneous subset of size $\kappa$ violating the assumption that $\kappa$ is a strong partition cardinal. This completes the proof.
\end{proof}

\Begin{definition}{kunen function definition}
Let $\kappa$ be a cardinal. Let $\mu$ be a $\kappa$-complete normal ultrafilter on $\kappa$. 

A function $f : \kappa \rightarrow \kappa$ is a block function if and only if $\{\alpha \in \kappa : f(\alpha) < |\alpha|^+\} \in \mu$.

Let $\Xi : \kappa \times \kappa \rightarrow \kappa$. For each $\alpha < \kappa$, let $\delta^{\Xi}_\alpha = \sup\{\Xi(\alpha,\beta) : \beta < \alpha\}$. Let $\Xi_\alpha : \alpha \rightarrow \delta^\Xi_\alpha$ be defined by $\Xi_\alpha(\beta) = \Xi(\alpha,\beta)$. 

$\Xi$ is a Kunen function for $f$ with respect to $\mu$ if and only if $K^\Xi_f = \{\alpha < \kappa : f(\alpha) \leq \delta^\Xi_\alpha \wedge \Xi_\alpha \text{ is a surjection}\} \in \mu$. (Here, when one says that $\Xi_\alpha$ is a surjection, one is considering $\Xi_\alpha$ as a function $\Xi_\alpha : \alpha \rightarrow \delta^\Xi_\alpha$.) $K^\Xi_f$ is the set of $\alpha$ on which $\Xi$ provides a bounding for $f$.

For $\beta < \kappa$, let $\Xi^\beta : \kappa \rightarrow \kappa$ be defined by $\Xi^\beta(\alpha) = \Xi(\alpha,\beta)$ where $\alpha > \beta$ and $0$ otherwise. $\Xi^\beta$ is a block function provided that $\{\alpha < \kappa : \Xi_\alpha \text{ is a surjection}\} \in \mu$.
\end{definition}

\Begin{fact}{selection of representative from kunen function}
Assume $\ZF$. Let $\mu$ be a normal measure on a cardinal $\kappa$. Suppose $f : \kappa \rightarrow \kappa$ is a block function which possesses a Kunen function $\Xi$ with respect to $\mu$. Suppose $G \in \prod_{\alpha \in \kappa} f(\alpha) \slash \mu$. Then there is a $\beta < \kappa$ so that $[\Xi^\beta]_\mu = G$
\end{fact}

\begin{proof}
Take any $g \in G$. Let $A = \{\alpha \in K^\Xi_f : g(\alpha) < f(\alpha)\}$. Note that for $\alpha \in A$, $g(\alpha) < f(\alpha) < \delta_\alpha^\Xi$. Define $\Phi : A \rightarrow \kappa$ by $\Phi(\alpha)$ is the least $\beta < \alpha$ so that $g(\alpha) = \Xi(\alpha,\beta)$ which exists since $\Xi_\alpha$ is a surjection onto $\delta_\alpha^\Xi$. Thus on $A$, $\Phi$ is a regressive function. By normality, there is some $\beta < \kappa$ so that $\Phi(\alpha) = \beta$ for $\mu$-almost all $\alpha$. Then $\Xi^\beta =_\mu g$. Hence $[\Xi^\beta]_\mu = G$. One can check this $\beta$ does not depend on the initial choice of $g \in G$. 
\end{proof}

\Begin{definition}{minimal code definition}
Let $\kappa$ be a cardinal. Let $\mu$ be a normal measure on $\kappa$. Let $h$ be a block function. Suppose $h$ possesses a Kunen function $\Xi$ with respect to $\mu$. An ordinal $\beta < \kappa$ is a minimal code (relative to $\Xi$) if and only if for all $\gamma < \beta$, $\neg(\Xi^\gamma =_\mu \Xi^\beta)$. Let $J^\Xi_h$ be the collection of $\beta$ which are minimal codes and $\Xi^\beta <_\mu h$. Define an ordering $\prec^\Xi_h$ on $J^\Xi_h$ by $\alpha \prec_h^\Xi \beta$ if and only if $\Xi^\alpha <_\mu \Xi^\beta$. By Fact \ref{selection of representative from kunen function}, for every $G < [h]_\mu$, there is a unique $\beta \in J^\Xi_h$ so that $\Xi^\beta \in G$. In this way, one says that $\beta$ is a minimal code for $G$ or for any $g \in G$ with respect to $\Xi$. 
\end{definition}

\Begin{fact}{ultrapower with kunen function is WO}
Let $\mu$ be a normal measure on a cardinal $\kappa$. Let $f : \kappa \rightarrow \kappa$ be a block function possessing a Kunen function $\Xi$ with respect to $\mu$. Then $\prod_{\alpha \in \kappa} f(\alpha) \slash \mu$, i.e. the initial segment of ${}^\kappa\slash\kappa \slash \mu$ determined by $[f]_\mu$, is a wellordering.

If every block function has a Kunen function, then $\prod_{\alpha < \kappa} |\alpha|^+ \slash \mu$ is wellfounded.

For each $F \in \prod_{\alpha < \kappa}|\alpha^+|\slash \mu$, $F < \kappa^+$. Thus $\prod_{\alpha < \kappa}|\alpha|^+ \slash \mu \leq \kappa^+$.
\end{fact}

\begin{proof}
Let $f$ be a block function possessing a Kunen function $\Xi$. Every $G \in \prod_{\alpha \in \kappa} f(\alpha) \slash \mu$ has a unique $\beta \in J^\Xi_f$ so that $[\Xi^\beta]_\mu = G$. There is a bijection of $J^\Xi_f$ with $\prod_{\alpha \in \kappa} f(\alpha) \slash \mu$ given by $\beta \mapsto [\Xi^\beta]_\mu$. This shows that $|\prod_{\alpha \in \kappa}f(\alpha)\slash\mu| \leq \kappa$.

Now suppose that $\prod_{\alpha \in \kappa}f(\alpha) \slash \mu$ is not wellfounded. Then $\prec^\Xi_f$ is an illfounded linear ordering on $J^\Xi_f \subseteq \kappa$. Let $\beta_0$ be the least (in the usual ordinal ordering) element of $J^\Xi_f$ so that the initial segment determined by $\beta_0$ in $(J^\Xi_f,\prec^\Xi_h)$ has no $\prec$-minimal element. Suppose $\beta_n$ has been defined, let $\beta_{n + 1}$ be the least ordinal in $J^\Xi_f$ which is $\prec^\Xi_f$-below $\beta_n$. This process defines a sequence of ordinals $\langle \beta_n : n \in \omega\rangle$ in $J^\Xi_f$. 

Then $\langle \Xi^{\beta_n} : n \in \omega\rangle$ is a sequence with the property that $[\Xi^{\beta_n}]_\mu$ is an infinite decreasing sequence in $\prod_{\alpha < \kappa} f(\alpha) \slash \mu$. Let $D_n = \{\alpha \in \kappa : \Xi^{\beta_{n + 1}}(\alpha) < \Xi^{\beta_n}(\alpha)\} \in \mu$. Then $\bigcap_{n \in \omega} D_n \in \mu$. Let $\xi \in \bigcap_{n \in \omega} D_n$. Then $\langle \Xi^{\beta_n}(\xi) : n \in \omega\rangle$ is an infinite decreasing sequence of ordinals. This is impossible.

Since $f$ was arbitrary, this shows that $\prod_{\alpha < \kappa}|\alpha|^+\slash\mu$ is wellfounded. By the first paragraph, each $F \in \prod_{\alpha < \kappa} |\alpha|^+ \slash \mu$ has cardinality less than or equal to $\kappa$. Thus $F < \kappa^+$. Thus $\prod_{\alpha < \kappa} |\alpha|^+ \slash \mu \leq \kappa^+$. 
\end{proof}

\Begin{definition}{minimal code identification}
Let $\mu$ be a normal measure on $\kappa$. Let $h : \kappa \rightarrow \kappa$ be a block function. Let $\Xi$ be a Kunen function for $h$ with respect to $\mu$. By Fact \ref{selection of representative from kunen function}, for each $G < [h]_\mu$, there is a minimal code $\beta \in J^\Xi_h$ so that $\Xi^\beta \in G$. Thus $(J_h^\Xi,\prec^\Xi_h)$ has the same ordertype as $[h]_\mu$. By Fact \ref{ultrapower with kunen function is WO}, $[h]_\mu$ is a wellordering. Let $\epsilon_h^\Xi \in \mathrm{ON}$ denote the ordertype of $([h]_{\mu},<)$ which is equal to the ordertype of $(J_h^\Xi,\prec^\Xi_h)$. Let $\pi^\Xi_h : \epsilon_h^\Xi \rightarrow (J_h^\Xi, \prec_h^\Xi)$ be the unique order-preserving isomorphism.
\end{definition}

Note that every function $f : \omega_1 \rightarrow \omega_1$ is (everywhere) a block function.

\Begin{theorem}{kunen tree and function}
(Kunen) Assume $\ZF + \AD$. Every function $f : \omega_1 \rightarrow \omega_1$ has a Kunen function with respect to the club measure $\mu$.
\end{theorem}

\begin{proof}
The Kunen function is derived from the Kunen tree which will be defined below. 

Recall that each $x \in \bairespace$ codes a relation $R_x$ defined in Definition \ref{WO definitions}. Also $\mathrm{WF}$ denotes the collection of $x$ so that $R_x$ is a wellfounded relation. Let $\pi_\mathrm{pair} : \omega \times \omega \rightarrow \omega$ denote a fixed recursive bijection.

Let $S \subseteq \omega \times \omega_1$ be the tree of partial rankings of relations on $\omega$ which is defined as follows: $(s,\bar{\alpha}) \in S$ if and only if $|s| = |\bar\alpha|$ and for all $i,j < |s|$, if $s(\pi_\mathrm{pair}(i,j)) = 1$, then $\bar{\alpha}(i) < \bar{\alpha}(j)$.

If $(x,f) \in [S] = \{(y,g) : (\forall n)(y \upharpoonright n, g\upharpoonright n) \in S)\}$, then $x \in \mathrm{WF}$ since $f$ is a ranking of $R_x$ into $\omega_1$. Conversely, if $x \in \mathrm{WF}$, then $R_x$ has a ranking $f$ with image in $\omega_1$. Thus $(x,f) \in [S]$. It has been shown that $\mathrm{WF} = \pi_1[[S]]$ is the projection of $[S]$ onto the first coordinate.

Let $T$ be a recursive tree on $\omega \times \omega$ so that $\pi_1[[T]] = \{a \in \bairespace : (\exists b)((a,b) \in [T])\} = \bairespace \setminus \mathrm{WF}$. The important observation is that $\bairespace \setminus \mathrm{WF}$ is a $\Sigma_1^1$ set which is $\mathbf{\Sigma}_1^1$-complete.

Let $\pi_\mathrm{seq} : {}^{<\omega}\omega \rightarrow \omega$ be a recursive bijection. For each $x \in \bairespace$, let $\sigma_x : {}^{<\omega}\omega \rightarrow \omega$ be a strategy in the usual integer game defined by $\sigma_x(s) = n$ if and $x(\pi_\mathrm{seq}(s)) = n$. Suppose $p \in {}^{<\omega}\omega$ and $s \in {}^{<\omega}\omega$, then let $\sigma_{p} * s$ denote the partial play where Player 1 uses the partial strategy $\sigma_{p}$ and Player 2 plays the bits of $s$ turn by turn. The game goes on for as long as $p$ codes a response to the partial play that is produced each turn.

Define a tree $K$ on $\omega \times \omega \times \omega_1 \times \omega \times \omega$ as follows: $(p,s,\bar{\alpha},t,u) \in K$ if and only if the conjunction of the following holds:

\noindent (1) $(s,\bar{\alpha}) \in S$.

\noindent (2) $(t,u) \in T$.

\noindent (3) $\sigma_{p} * s$ is a substring of $t$.

The meaning of $K$ becomes clear if one looks at what a path through $K$ represents: If $(x,y,f,v,w) \in [K]$, then $\sigma_x(y) = v$, $(y,f) \in [S]$, and $(v,w) \in [T]$. Therefore, $y \in \mathrm{WF}$ and $v \in \reals \setminus \mathrm{WF}$ since $S$ and $T$ are trees that project onto $\mathrm{WF}$ and $\reals \setminus \mathrm{WF}$.

The above defines the tree $K$. One now introduces an arbitrary function $f : \omega_1 \rightarrow \omega_1$. Consider the game $G_f$ defined as follows:
$$\begin{tikzpicture}
\node at (0,0) {$G_f$};
\node at (1,.5) {I};
\node at (2,.5) {$y(0)$};
\node at (4,.5) {$y(1)$};
\node at (6,.5) {$y(2)$};
\node at (8,.5) {$y(3)$};
\node at (11,.5) {$y$};

\node at (1,-.5) {II};
\node at (3,-.5) {$v(0)$};
\node at (5,-.5) {$v(1)$};
\node at (7,-.5) {$v(2)$};
\node at (9,-.5) {$v(3)$};
\node at (11,-.5) {$v$};
\end{tikzpicture}$$
Player 2 wins if and only if $y \in \WO \Rightarrow (v \in \mathrm{WF} \wedge \mathrm{rk}(T_v) > \sup\{f(\alpha) : \alpha \leq \mathrm{ot}(y))\}$. Here $T_v = \{u \in {}^{<\omega}\omega : (u,v\upharpoonright |u|) \in T\}$. Note that since $T$ projects onto $\bairespace \setminus \mathrm{WF}$, if $v \in \mathrm{WF}$, then $T_v$ is a wellfounded tree.

Player 2 must have the winning strategy: Suppose $\rho$ is a Player 1 winning strategy. $\rho[\bairespace] \subseteq \mathrm{WO}$ since otherwise Player 2 can win against $\rho$. Thus $\rho[\bairespace]$ is a $\analytic$ subset of $\WO$. By the bounding principle, let $\mu < \omega_1$ be such that $\ot(r) < \mu$ for all $r \in \rho[\bairespace]$. Let $\zeta = \sup\{f(\alpha) : \alpha < \mu\}$. $\zeta < \omega_1$ because $\omega_1$ is regular. Since the projection of $T$ is $\reals \setminus \mathrm{WF}$ which is a $\analytic$-complete set, one must have that $\sup\{\mathrm{rk}(T_r) : r \in \mathrm{WF}\} = \omega_1$. Choose some $v \in \mathrm{WF}$ so that $\mathrm{rk}(T_v) > \zeta$. If Player 2 plays $v$ against $\rho$, Player 2 will win. This contradicts the fact that $\rho$ is a Player 1 winning strategy.

Thus there is a winning strategy $\sigma$ for the game $G_f$. Let $x \in \bairespace$ be such that $\sigma_x = \sigma$. 

Let $K_x$ be the tree consisting of $(s,\bar\alpha,t,u)$ so that there is an $n \in \omega$ such that $|s| = |\bar\alpha| = |t| = |u| = n$ and $(x\upharpoonright n, s,\bar\alpha,t,u) \in K$. For each $\alpha < \omega_1$, let $K_x\upharpoonright \alpha$ denote the restrict of $K_x$ to ${}^{<\omega}(\omega\times\alpha\times\omega\times\omega)$. Note that $K_x\upharpoonright \alpha$ is wellfounded. To see this, suppose otherwise. This means there is some $y,v,w \in \bairespace$ and $f \in {}^{<\omega}\alpha$ so that $(x,y,f,v,w) \in [K]$. Thus $\sigma_x(y) = \sigma(y) = v$, $y \in \mathrm{WO}$, and $v \in \reals \setminus \mathrm{WF}$. However, $\sigma$ is a Player 2 winning strategy and $y \in \mathrm{WO}$, so one must have that $v \in \mathrm{WF}$. This is a contradiction.

Let $\alpha \geq \omega$. Suppose $y \in \mathrm{WO}$ with $\mathrm{ot}(y) = \alpha$. Then there is a ranking $f$ of $R_y$ using ordinals below $\alpha$. Let $v = \sigma_x(y) = \sigma(y)$. Since $\sigma$ is a Player 2 winning strategy, $\mathrm{rk}(T_v) > f(\ot(y))$. Note that $K_x\upharpoonright \alpha$ has a subtree which is isomorphic to $T_v$. In particular, this subtree is 
$$\hat{T}_v = \{(y\upharpoonright n, f\upharpoonright n, t,u) \in K_x : n \in \omega \wedge |t| = |u| = n \wedge (t,u) \in T_v\}.$$
This implies that $K_x\upharpoonright \alpha$ is a wellfounded tree of rank greater than $\mathrm{rk}(T_v) > f(\alpha)$.

Note that there is a club set of $\alpha$ so that $\alpha$ is closed under the G\"odel pairing function. Using this pairing function on $\alpha$, one can define uniformly a bijection $\pi_\alpha : \alpha \rightarrow {}^{<\omega}\alpha$. For all $\alpha$ closed under the G\"odel pairing function, define $\Xi(\alpha,\beta)$ to be the rank of $\pi_\alpha(\beta)$ in $K_x\upharpoonright\alpha$ whenever $\beta < \alpha$ and $\pi_\alpha(\beta) \in K_x\upharpoonright\alpha$. (Technically in the definition of a Kunen function, $\Xi$ should be a function on $\omega_1 \times \omega_1$; however, the value of $\Xi$ on $(\alpha,\beta)$ where $\beta \geq \alpha$ is not relevant in any applications.) If $\pi_\alpha(\beta)$ does not belong to $K_x\upharpoonright\alpha$ then let $\Xi(\alpha,\beta) = 0$.

$\Xi$ is a Kunen function for $f$. For any $\alpha \geq \omega$ which is closed under the G\"odel pairing function, $\Xi_\alpha$ is a surjection onto the $\mathrm{rk}(K_x\upharpoonright\alpha)$. Using the notation of Definition \ref{kunen function definition}, $\mathrm{rk}(K_x\upharpoonright \alpha) = \delta_\alpha^\Xi$. Also, it was shown above that $\mathrm{rk}(K_x\upharpoonright\alpha) > f(\alpha)$ and hence $\delta_\alpha^\Xi > f(\alpha)$. This verifies that $\Xi$ is a Kunen function for $f$ with respect to $\mu$.
\end{proof}

\Begin{remark}{kunen tree remarks}
The tree $K$ produced in Theorem \ref{kunen tree and function} is called the Kunen tree. 

There are some simplifications of $K$ that can be made. $K$ is a tree on $\omega\times\omega\times\omega_1\times\omega\times\omega$. One can merge the last four coordinates on $\omega\times\omega_1\times\omega\times\omega$ into $\omega_1$ to produce a tree on $\omega\times\omega_1$ with the same features. One can also Kleene-Brouwer order the various trees to produce linear orderings which can be used to produce the desired Kunen functions. See \cite{Structural-Consequences-AD} for more details on the Kunen tree.

In this survey, one will only use the existence of Kunen functions for functions $f : \omega_1 \rightarrow \omega_1$; thus, it is not important here how uniformly these Kunen functions are obtained. However, the proof shows that all Kunen functions are index uniformly by reals. For instance, there is a single tree $K$, the Kunen tree, so that for any $f : \omega_1 \rightarrow \omega_1$, there is a section of this tree by some strategy which can be used to produce the Kunen function for $f$. The Kunen tree $K$ is also $\Delta_1^1$ in the codes. See \cite{Structural-Consequences-AD} for the details and the precise meaning of these remarks.
\end{remark}

\Begin{corollary}{bound on omega1 ultrapower}
Assume $\ZF + \AD$. Let $\mu$ denote the club measure on $\omega_1$. $|\prod_{\omega_1}\omega_1 \slash \mu| \leq \omega_2$. 
\end{corollary}

\Begin{fact}{value of omega1 ultrapower}
(Martin) Assume $\ZF + \AD$. Let $\mu$ denote the club measure on $\omega_1$. Then $\prod_{\omega_1} \omega_1 \slash \mu = \omega_2$ and is a regular cardinal.
\end{fact}

\begin{proof}
Note that the club filter on $\omega_1$, $\mu$, is a normal measure on $\omega_1$. Therefore by Fact \ref{ultrapower normal strong partition is regular}, the ultrapower is a regular cardinal. Thus it must be greater than or equal $\omega_2$. Then Corollary \ref{bound on omega1 ultrapower} implies that the ultrapower is $\omega_2$.
\end{proof}

\Begin{definition}{kunen function and higher type}
Let $\mu$ be a normal measure on a cardinal $\kappa$. Let $h : \kappa \rightarrow \kappa$ be a function so that $h(\alpha) > 0$ $\mu$-almost everywhere. Let $T^h = \{(\alpha,\beta) \in \kappa\times\kappa : \beta < h(\alpha)\}$. Let $\mathcal{T}^h = (T^h,\sqsubset)$ where $\sqsubset$ is the lexicographic ordering. Note that $\ot(\mathcal{T}^h) = \kappa$. 

Suppose $F : \mathcal{T}^h \rightarrow \kappa$ is an order-preserving function. Let $g \in {}^\kappa\kappa$ be such that $g <_\mu h$. Let $A^g = \{\alpha : g(\alpha) < h(\alpha)\}$. Let $F^g \in {}^\kappa\kappa$ be defined by
$$F^g(\alpha) = \begin{cases}
F(\alpha,g(\alpha)) & \quad \alpha \in A^g \\
F(\alpha,0) & \quad \text{otherwise}
\end{cases}$$
Note that if $g_1 <_\mu g_2 <_\mu h$, then $F^{g_1} <_\mu F^{g_2}$. 

Fix a Kunen function $\Xi$ for $h$. Recall that $\epsilon_h^\Xi$ is the ordertype of the wellordering $(J_h^\Xi,\prec_h^\Xi)$ and $\pi^\Xi_h : \epsilon_h^\Xi \rightarrow (J_h^\Xi,\prec^\Xi_h)$ is the unique order isomorphism. If $\beta \in \epsilon_h^\Xi$, then let $F^{(\beta)} = F^{\Xi^{\pi^\Xi_h(\beta)}}$. Let $\funct(F) : \epsilon_h^\Xi \rightarrow \ON$, be defined by $\funct(F)(\alpha) = [F^{(\alpha)}]_\mu$. 
\end{definition}

$\mathcal{T}^h = (T^h, \sqsubset)$ is order isomorphic to $\kappa$. It is merely a reorganization of $\kappa$ into successive blocks of length $h(\alpha)$. If $F : \mathcal{T}^h \rightarrow \kappa$ is order preserving and $\beta < \epsilon^\Xi_h$, then $F^{(\beta)} : \kappa \rightarrow \kappa$ is defined by letting $F^{(\beta)}(\alpha)$ be the value of $F$ at $(\alpha,\Xi^{\pi^\Xi_h(\beta)}(\alpha))$, which can be construed to be the $\Xi^{\pi^\Xi_h(\beta)}(\alpha)^\text{th}$ element in the $\alpha^\text{th}$ block of length $h(\alpha)$.

\Begin{lemma}{higher type sliding lemma}
(Sliding lemma) Let $\mu$ be a normal measure on a cardinal $\kappa$. Let $h : \kappa \rightarrow \kappa$ be a block function possessing a Kunen function $\Xi$ with respect to $\mu$. 

Let $\langle g_\alpha : \alpha \in \epsilon_h^\Xi\rangle$ be an order preserving sequence elements of ${}^\kappa\kappa$ which are non-constant $\mu$-almost everywhere: order preserving means that if $\alpha < \beta < \epsilon_h^\Xi$, then $g_\alpha <_\mu g_\beta$. 

Then there exists an order-preserving $F : \mathcal{T}^h \rightarrow \kappa$ so that for all $\beta < \epsilon_h^\Xi$, $F^{(\beta)} =_\mu g_\beta$. 

Moreover, if $X \subseteq \kappa$ and for all $\beta < \epsilon_h^\Xi$, $\rang(g_\beta) \subseteq X$, then $\rang(F) \subseteq X$.
\end{lemma}

\begin{proof}
Fix some $X \subseteq \kappa$ so that $\rang(g_\beta) \subseteq X$ for all $\beta < \epsilon_h^\Xi$, $\rang(g_\beta) \subseteq X$. 

Define $F : \mathcal{T}^h \rightarrow \kappa$ by recursion as follows:

Let $(\alpha,\beta) \in T^h$. Suppose $F(\alpha',\beta')$ has been defined for all $(\alpha',\beta') \sqsubset (\alpha,\beta)$. (Recall $\mathcal{T}^h = (T^h,\sqsubset)$ where $\sqsubset$ is the lexicographic ordering.) 

(Case I) There is a $\gamma < \epsilon_h^\Xi$ so that $\beta = \Xi^{\pi_h^\Xi(\gamma)}(\alpha)$:

Let $\gamma$ be least with this property. Let $F(\alpha,\beta)$ be the least element in the range $g_\gamma$ which is greater than or equal to $\sup\{F(\alpha',\beta') : (\alpha',\beta') \sqsubset (\alpha,\beta)\}$. 

(Case II) There is no $\gamma < \epsilon_h^\Xi$ so that $\beta = \Xi^{\pi_h^\Xi(\gamma)}(\alpha)$:

Let $F(\alpha,\beta)$ be the $\omega^\text{th}$-element of $X$ above $\sup\{F(\alpha',\beta') : (\alpha',\beta') \sqsubset (\alpha,\beta)\}$. 

Claim: For all $\gamma < \epsilon_h^\Xi$, $F^{(\gamma)} =_\mu g_\gamma$. 

This will be proved by induction on $\epsilon_h^\Xi$.

It is clear $B = \{\alpha \in \kappa : \Xi^{\pi_h^\Xi(0)}(\alpha) = 0\} \in \mu$. For $\alpha \in B$, $F^{(0)} = F(\alpha,\Xi^{\pi_h^\Xi(0)}(\alpha)) \in \rang(g_{0})$. Thus on $B$, one must have that $F^{(0)}(\alpha) \geq g_{0}(\alpha)$. Suppose that $F^{(0)} >_\mu g_{0}$. Let $C = \{\alpha \in B : F^{(0)}(\alpha) > g_0(\alpha)\} \in \mu$. 

Since for $\alpha \in C$, $\Xi^{\pi_h^\Xi(0)}(\alpha) = 0$, there must be some $(\alpha',\beta') \sqsubset (\alpha,0)$ with $\alpha' < \alpha$ so that $F(\alpha',\beta') \geq g_{0}(\alpha)$. Define $\Phi : C \rightarrow \kappa$ by letting $\Phi(\alpha)$ be the least $\alpha' < \alpha$ so that there exists some $\beta'$ with $(\alpha',\beta') \sqsubset (\alpha,\beta)$ and $F(\alpha',\beta') \geq g_{0}(\alpha)$. Thus $\Phi$ is regressive on $C \in \mu$. There is some $\alpha'$ and a $D \in \mu$ so that $\Phi(\alpha) = \alpha'$ for all $\alpha \in D$. Thus for all $\alpha \in D$, $g_{0}(\alpha) \leq F(\alpha',\beta')$ for some $\beta' < h(\alpha')$. By $\kappa$-completeness, there is some $\bar\beta < h(\alpha')$ and an $E \subseteq D$ with $E \in \mu$ so that for all $\alpha \in E$, $g_{0}(\alpha) \leq F(\alpha',\bar\beta)$. This is impossible since $g_{0}$ is not constant $\mu$-almost everywhere. It has been shown that $F^{(0)} = g_{0}$. 

Suppose $\gamma < \epsilon_h^\Xi$ and that it has been shown that $F^{(\gamma')} =_{\mu} g_{\gamma'}$ for all $\gamma' < \gamma$. One will seek to show that $F^{(\gamma)} = g_\gamma$. 

Let $A = \{\alpha : (\forall \gamma' \in \gamma)(\Xi^{\pi^\Xi_h(\gamma)}(\alpha) \neq \Xi^{\pi^\Xi_h(\gamma')}(\alpha))\}$. If $A \notin \mu$, then by $\kappa$-completeness, there is some $\gamma' \in \gamma$ so that $\Xi^{\pi^\Xi_h(\gamma')} =_{\mu} \Xi^{\pi^\Xi_h(\gamma)}$. This is impossible since $\gamma' < \gamma$ implies $\Xi^{\pi_h^\Xi(\gamma')} <_\mu \Xi^{\pi^\Xi_h(\gamma)}$.

By induction, it has been shown that $F^{(0)} = g_0$. Let $B \subseteq A$ with $B \in \mu$ have the property that for all $\alpha \in B$, $\Xi^{\pi^\Xi_h(0)}(\alpha) = 0$, $F^{(0)}(\alpha) = g_{0}(\alpha)$, and $g_{\gamma}(\alpha) > g_{0}(\alpha)$. 

Suppose $F^{(\gamma)}$ is not equal to $g_\gamma$ for $\mu$-almost all $\alpha$. For $\alpha \in B$, $F^{(\gamma)}(\alpha) \in \rang(g_\gamma)$. Thus one must have that there is a set $C \subseteq B$ with $C \in \mu$ so that for all $\alpha \in C$, $F^{(\gamma)}(\alpha) > g_\gamma(\alpha)$. There must be some $(\alpha',\beta') \sqsubset (\alpha,\Xi^{\pi^\Xi_h(\gamma)}(\alpha))$ so that $F(\alpha',\beta') \geq g_\gamma(\alpha)$. However, $B$ was chosen so that for all $\alpha \in B$, $\Xi^{\pi_h^\Xi(0)}(\alpha) = 0$ and $F^{(0)}(\alpha) = g_{0}(\alpha) < g_{\gamma}(\alpha)$. Therefore, for $\alpha \in C$, $F^{(0)}(\alpha) < g_{\gamma}(\alpha) < F^{(\gamma)}(\alpha)$. Thus the least $(\alpha',\beta') \sqsubset (\alpha,\Xi^{\pi_h^\Xi(\gamma)}(\alpha))$ with $F(\alpha',\beta') \geq g_\gamma(\alpha)$ must have that $\alpha' = \alpha$. 

Therefore, for all $\alpha \in C$, there is some $\beta' < \Xi^{\pi_h^\Xi(\gamma)}(\alpha)$ so that $F(\alpha,\beta') \geq g_\gamma(\alpha)$. Let $\Phi(\alpha)$ be this $\beta'$. Then $\Phi <_\mu \Xi^{\pi_h^\Xi(\gamma)} < h$. Thus there is some $\bar{\gamma} < \gamma$ so that $\Phi =_\mu \Xi^{\pi_h^\Xi(\bar\gamma)}$. By induction, $F^{(\bar\gamma)} =_\mu g_{\bar{\gamma}}$. However by definition of $\Phi$, $F^\Phi \geq_\mu g_\gamma$. Also $F^\Phi =_\mu F^{(\bar{\gamma})} =_\mu g_{\bar{\gamma}}$. Thus $g_{\bar\gamma} \geq_\mu g_\gamma$. Since $\langle g_\delta : \delta \in \epsilon_h^\Xi\rangle$ is an increasing sequence and $\bar{\gamma} < \gamma$, one has that $g_{\bar{\gamma}} <_\mu g_{\gamma}$. This is a contradiction. This shows that $F^{\gamma} = g_\gamma$. 

The lemma has been proved.
\end{proof}

It will be helpful to have the correct type version of the sliding lemma. A sketch of the necessary modification will be given:

\Begin{lemma}{correct type sliding lemma}
(Correct-type sliding lemma) Let $\mu$ be a normal measure on a cardinal $\kappa$. Let $h : \kappa \rightarrow \kappa$ be a block function possessing a Kunen function $\Xi$ with respect to $\mu$. 

Let $\langle g_\alpha : \alpha < \epsilon^\Xi_h\rangle$ be an increasing sequence of functions from $\kappa$ to $\kappa$ of the correct type which are $\mu$-almost everywhere non-constant. Further, suppose there is a sequence $\langle G_\alpha : \alpha < \epsilon_h^\Xi\rangle$ where each $G_\alpha : \kappa \times \omega \rightarrow \kappa$ is a function witnessing that $g_\alpha$ has uniform cofinality $\omega$. Suppose $\langle [g_\alpha]_\mu : \alpha < \epsilon_h^\Xi\rangle$ is discontinuous everywhere.

Then there exists an order preserving function $F : \mathcal{T}^h \rightarrow \kappa$ which is of the correct type so that for all $\alpha < \epsilon_h^\Xi$, $F^{(\alpha)} =_\mu g_\alpha$. 

Moreover, if $C \subseteq \kappa$ is an $\omega$-club set such that for all $\beta \in \epsilon^\Xi_h$, $\rang(g_\beta) \subseteq C$, then $\rang(F) \subseteq C$. 
\end{lemma}

\begin{proof}
Fix some $C \subseteq \kappa$ be an $\omega$-club so that $\rang(g_\beta) \subseteq C$ for all $\beta < \epsilon_h^\Xi$. 

Define $F : \mathcal{T}^h \rightarrow \kappa$ by recursion as follows:

Let $(\alpha,\beta) \in T^h$. Suppose $F(\alpha',\beta')$ has been defined for all $(\alpha',\beta') \sqsubseteq (\alpha,\beta)$. 

(Case I) There is some $\gamma < \epsilon^\Xi_h$ so that $\beta = \Xi^{\pi_h^\Xi(\gamma)}(\alpha)$: 

Let $\gamma$ be least with this property. Let $F(\alpha,\beta)$ be the least element in the range of $g_\gamma$ which is greater than $\sup\{F(\alpha',\beta') : (\alpha',\beta') \sqsubset (\alpha,\beta)\}$. 

(Case II) There is no $\gamma < \epsilon^\Xi_h$ so that $\beta = \Xi^{\pi^\Xi_h(\gamma)}(\alpha)$: 

Find $\delta$ least so that $\enum_C(\delta) > \sup\{F(\alpha',\beta') : (\alpha',\beta') \sqsubset (\alpha,\beta)\}$. Let $F(\alpha,\beta) = \enum_C(\delta + \omega)$. 

Claim 1: $F : \mathcal{T}^h \rightarrow C$ is a function of the correct type. 

To prove this: Note that it is clear from the construction that $F$ is discontinuous everywhere. 

One will create a function $H : T^h \times \omega \rightarrow \kappa$ to witness the uniform cofinality of $F$. 

Let $(\alpha,\beta) \in T^h$. 

Suppose Case I had occurred at $(\alpha,\beta)$ with $\gamma < \epsilon_h^\Xi$ least so that $\beta = \Xi^{\pi_h^\Xi(\gamma)}(\alpha)$. Let $\delta < \kappa$ be least so that $g_\gamma(\delta) > \sup\{F(\alpha',\beta') : (\alpha',\beta') \sqsubset (\alpha,\beta)\}$. Now let $n$ be least so that $G_\gamma(\delta,n) > \sup\{F(\alpha',\beta') : (\alpha',\beta') \sqsubset (\alpha,\beta)\}$. Define $H((\alpha,\beta),m) = G_\gamma(\delta,n + m)$. 

Suppose Case II had occurred at $(\alpha,\beta)$. Let $H((\alpha,\beta),m)$ be the $m^\text{th}$ element of $C$ above $\sup\{F(\alpha',\beta') : (\alpha',\beta') \sqsubset (\alpha,\beta)\}$. 

This function $H$ witnesses that $F$ has uniform cofinality $\omega$. $F$ has the correct type.

Claim 2: For all $\gamma < \epsilon_h^\Xi$, $F^{(\gamma)} =_\mu g_\gamma$. 

To show Claim 2: This will be proved by induction on $\epsilon^\Xi_h$. One will indicate some of the necessary modification from the proof of Lemma \ref{higher type sliding lemma}.

Suppose that $\gamma < \epsilon^\Xi_h$ and that it has been shown that $F^{(\gamma')} =_\mu g_{\gamma'}$ for all $\gamma' < \gamma$. Let $A = \{\alpha : (\forall \gamma' \in \gamma)(\Xi^{\pi^\Xi_h(\gamma)}(\alpha) \neq \Xi^{\pi^\Xi_h(\gamma')}(\alpha))\}$. As in Lemma \ref{higher type sliding lemma}, $A \in \mu$. Also just as before, one can show that $B \in \mu$ where $B \subseteq A$ and has the property that for all $\alpha \in B$, $\Xi^{\pi^\Xi_h(0)}(\alpha) = 0$, $F^{(0)}(\alpha) = g_0(\alpha)$, and $g_\gamma(\alpha) > g_0(\alpha)$. 

Now suppose $F^{(\gamma)}$ is not equal to $g_\gamma$ for $\mu$-almost all $\alpha$. In the present situation, there are two ways this can happen:

(i) There is a $D \subseteq B$ with $D \in \mu$ so that for all $\alpha \in D$, there exists some $(\alpha',\beta') \sqsubset (\alpha,\Xi^{\pi^\Xi_h(\gamma)}(\alpha))$ so that $F(\alpha',\beta') \geq g_\gamma(\alpha)$. As argued in Lemma \ref{higher type sliding lemma}, this can not occur.

(ii) There is a $D \subseteq B$ with $D \in \mu$ so that for all $\alpha \in D$, $g_\gamma(\alpha) = \sup\{F(\alpha',\beta') : (\alpha',\beta') \sqsubset (\alpha,\Xi^{\pi^\Xi_h(\gamma)}(\alpha))\}$. 

For each $n \in \omega$, let $\Phi_n : D \rightarrow \omega_1$ be defined as follows: Since $\alpha \in D \subseteq B$, one has that $F^{(0)}(\alpha) = g_0(\alpha) < g_\gamma(\alpha)$. For each $n \in \omega$ and $\alpha \in D$, let $\Phi_n(\alpha)$ be the least $\beta' < \Xi^{\pi^\Xi_h(\gamma)}(\alpha)$ so that $F(\alpha,\beta') > G_\gamma(\alpha,n)$. 

By Fact \ref{selection of representative from kunen function}, there is some $\gamma_n < \gamma$ so that $\Phi_n =_\mu \Xi^{\gamma_n}$. By the induction hypothesis, $F^{\Xi^{\pi_h^\Xi(\gamma_n)}} = F^{(\gamma_n)} =_\mu g_{\gamma_n}$. By construction, $\sup\{[g_{\gamma_n}]_\mu : n \in \omega\} = [g_\gamma]_\mu$. However, by assumption, $\langle [g_\alpha]_{\mu} : \alpha < \epsilon_h^\Xi\rangle$ was discontinuous everywhere.

This complete the proof of the lemma.
\end{proof}

\Begin{fact}{uniqueness of the slide}
Let $\mu$ be a normal measure on a cardinal $\kappa$. Let $h : \kappa \rightarrow \kappa$ be a block function possessing a Kunen function $\Xi$ with respect to $\mu$. Suppose $F_0,F_1 \in [\omega_1]^{\mathcal{T}^h}$ have the property that $F_0^{(\beta)} =_\mu F_1^{(\beta)}$ for all $\beta < \epsilon^\Xi_h$. Then for $\mu$-almost all $\alpha$, $F_0(\alpha,\beta) = F_1(\alpha,\beta)$ for all $\beta < h(\alpha)$. 
\end{fact}

\begin{proof}
Suppose $A = \{\alpha : (\exists \beta < h(\alpha))(F_0(\alpha,\beta) \neq F_1(\alpha,\beta))\} \in \mu$. Define $g : \omega_1 \rightarrow \omega_1$ by $g(\alpha)$ is the least $\beta < h(\alpha)$ so that $F_0(\alpha,\beta) \neq F_1(\alpha,\beta)$. Since $g <_\mu h$, Fact \ref{selection of representative from kunen function} implies there is some $\gamma < \epsilon^\Xi_h$ so that $g =_\mu \Xi^{\pi_h^\Xi(\gamma)}$. Then $\neg(F_0^{(\gamma)} =_\mu F_1^{(\gamma)})$. This contradicts the assumptions.
\end{proof}

\Begin{fact}{representation alpha increasing ordinal sequence ultrapower}
Assume $\ZF + \AD$. Let $\mu$ denote the club measure on $\omega_1$. Let $\alpha < \omega_2$. Let $A \subseteq \omega_1$ with $|A| = \omega_1$. Let $B = [A]^{\omega_1} \slash \mu$, which is a set of cardinality $\omega_2$. Let $h : \omega_1 \rightarrow \omega_1$ have the property that $h(\beta) > 0$ for all $\beta < \omega_1$ and $[h]_\mu = \alpha$. Let $\Xi$ be a Kunen function for $h$ with respect to $\mu$.

For all $I \in [B]^\alpha$, there exists some $F: \mathcal{T}^h \rightarrow A$ which is order preserving and for all $\gamma < \alpha = \epsilon_h^\Xi$, $[F^{(\gamma)}]_\mu = I(\gamma)$.
\end{fact}

\begin{proof}
Since $\omega_2 = {}^{\omega_1}\omega_1\slash \mu$ is regular by Fact \ref{ultrapower normal strong partition is regular}, $\rang(I) < \omega_2$. Let $k \in [\omega_1]^{\omega_1}$ be such that $\rang(I) < [k]_\mu$. Let $\Xi'$ be a Kunen function for $k$ with respect to $\mu$. 

Let $J = \{\beta \in \epsilon_k^{\Xi'} : (\exists \gamma < \alpha)(I(\gamma) = [{\Xi'}^{\pi_{k}^{\Xi'}(\beta)}]_\mu)\}$. $J \subseteq \epsilon_k^{\Xi'}$ and $\ot(J) = \alpha = \epsilon_h^\Xi$. Let $\rho : \epsilon^\Xi_h \rightarrow J$ be the unique order isomorphism. Define $g'_\gamma = {\Xi'}^{\pi_k^{\Xi'}(\rho(\gamma))}$. Since $[g'_\gamma]_\mu = I(\gamma)$ and $I(\gamma) \in B = [A]^{\omega_1} \slash \mu$, the set $B_\gamma = \{\alpha \in \omega_1 : g'_\gamma(\alpha) \in A\} \in \mu$. By Fact \ref{normal measure and closure points}, $g'_\gamma =_\mu g'_\gamma \circ \enum_{B_\gamma}$.  Let $g_\gamma = g'_\gamma \circ \enum_{B_\gamma}$. $[g_\gamma]_\mu = I(\gamma)$ and $\rang(g_\gamma) \subseteq A$.  The result now follows from Lemma \ref{higher type sliding lemma}.
\end{proof}

\Begin{theorem}{omega2 is a weak partition cardinal}
(Martin-Paris) Assume $\ZF + \AD$. Let $\mu$ be the club measure on $\omega_1$. Then for all $\alpha < \omega_2$, the partition relation $\omega_2 \rightarrow (\omega_2)^\alpha_2$ holds. That is, $\omega_2$ is a weak partition cardinal.
\end{theorem}

\begin{proof}
By Fact \ref{value of omega1 ultrapower}, $\omega_2$ is isomorphic to ${}^{\omega_1}\omega_1 \slash \mu$. One will identify $\omega_2$ with ${}^{\omega_1}\omega_1 \slash \mu$.

Let $\alpha < \omega_2$ and $P : [\omega_2]^{\alpha} \rightarrow 2$ be a partition.

Let $h : \omega_1 \rightarrow \omega_1$ be such that $[h]_\mu = \alpha$ and fix a Kunen function $\Xi$ for $h$. Define a partition $Q : [\omega_1]^{\mathcal{T}^h} \rightarrow 2$ by $Q(F) = P(\funct(F)))$, where $\funct(F)$ is an $\alpha$-sequence of ordinals in $\omega_2$ defined in Definition \ref{kunen function and higher type} (relative to $\Xi$).

Since $\mathcal{T}^h$ is order isomorphic to $\omega_1$, $\omega_1 \rightarrow (\omega_1)^{\omega_1}_2$ implies that there is a $A \subseteq \omega_1$ with $|A| = \omega_1$ which is homogeneous for $Q$. Without loss of generality, suppose $A$ is homogeneous for $Q$ taking value $0$.

Note that $[A]^{\omega_1} \slash \mu$ has cardinality ${}^{\omega_1}\omega_1\slash \mu = \omega_2$. Let $B$ denote $[A]^{\omega_1}\slash \mu$. 

Let $I \in [B]^\alpha$. By Fact \ref{representation alpha increasing ordinal sequence ultrapower}, let $F : \mathcal{T}^h \rightarrow A$ be such that $[F^{(\gamma)}]_\mu = I(\gamma)$ for all $\gamma < \alpha$. Since $F \in [A]^{\mathcal{T}^h}$, one has $0 = Q(F) = P(\funct(F)) = P(I)$ since $I = \funct(F)$. $B$ is homogeneous for $P$ taking value $0$. The proof is complete.
\end{proof}

By Fact \ref{equivalence of partition properties}, a suitable form of the ordinary partition property will imply the appropriate correct type version of the partition property which uses a club set as its version of the homogeneous set. In particular, the ordinary weak partition property for $\omega_2$ proved in Theorem \ref{omega2 is a weak partition cardinal} implies the correct type club version of the weak partition property, that is $\omega_2 \rightarrow_* (\omega_2)^\alpha_2$ for all $\alpha < \omega_2$.

\Begin{corollary}{normal measure on omega2}
Assume $\ZF + \AD$. $W^{\omega_2}_\omega$ and $W^{\omega_2}_{\omega_1}$ are the only two $\omega_2$-complete normal ultrafilters on $\omega_2$. 
\end{corollary}

\begin{proof}
By Fact \ref{characterize normal measure by partition}.
\end{proof}

\Begin{fact}{representation of club subset of omega2}
Assume $\ZF + \AD$. Let $\mu$ be the club measure on $\omega_1$. If $C \subseteq \omega_1$ is a club subset of $\omega_1$, then $[C]^{\omega_1} \slash \mu$ is a club subset of $\omega_2$. 

If $D \subseteq \omega_2$ is club, then there is a club $C \subseteq \omega_1$ so that $[C]^{\omega_1} \slash \mu \subseteq D$. 
\end{fact}

\begin{proof}
Let $\epsilon < \omega_2$. Suppose $\langle \nu_\gamma : \gamma < \epsilon\rangle$ is an increasing sequence in $[C]^{\omega_1} \slash \mu$. Let $h : \omega_1 \rightarrow \omega_1$ with $h(\alpha) > 0$ for all $\alpha \in \omega_1$ be such that $\epsilon = [h]_\mu$. Let $\Xi$ be a Kunen function for $h$ with respect to $\mu$.

By Fact \ref{representation alpha increasing ordinal sequence ultrapower}, there is a function $F : \mathcal{T}^h \rightarrow C$ which is order preserving so that $[F^{(\gamma)}]_\mu = \nu_\gamma$.

Let $\ell : \omega_1 \rightarrow \omega_1$ be defined by $\ell(\alpha) = \sup\{F(\alpha,\beta) : \beta < h(\alpha)\}$. Since $\rang(F) \subseteq C$, one has that $\ell(\alpha) \in C$ for each $\alpha \in \omega_1$. 

Suppose $j : \omega_1 \rightarrow \omega_1$ is such that $j <_\mu \ell$. Let $B = \{\alpha \in \omega : j(\alpha) < \ell(\alpha)\} \in \mu$. For $\alpha \in B$, let $p(\alpha)$ be the least $\beta < h(\alpha)$ so that $j(\alpha) < F(\alpha,\beta)$. For $\alpha \notin B$, let $p(\alpha) = 0$. Since $p <_\mu h$, Fact \ref{selection of representative from kunen function} implies there $\gamma < \alpha$ so that $p =_\mu \Xi^{\pi^\Xi_h(\gamma)}$. Then $j <_\mu F^{(\gamma)}$. Thus $[j]_\mu < \nu_{\gamma}$. This establishes that $[\ell]_\mu$ is the limit of $\langle \nu_\gamma : \gamma < \epsilon\rangle$. Since $\ell \in [C]^{\omega_1}$, this shows the supremum belongs to $[C]^{\omega_1}\slash\mu$. This shows that $[C]^{\omega_1} \slash \mu$ is closed.  It is easy to see that $[C]^{\omega_1} \slash \mu$ is unbounded. 

Now suppose $D \subseteq \omega_2$ is club. Let $\mathcal{T}^2 = (\omega_1 \times 2, \sqsubset)$ where $\sqsubset$ is the lexicographic ordering. If $F : \mathcal{T}^2 \rightarrow \omega_1$ is an increasing function, then let $F_0,F_1 \in [\omega_1]^{\omega_1}$ be defined by $F_i(\alpha) = F(\alpha,i)$. Note that $F_0 <_\mu F_1$. 

Define $P : [\omega_1]^{\mathcal{T}^2} \rightarrow 2$ by $P(F) = 0$ if and only if $(\exists \alpha \in D)([F_0]_\mu < \alpha < [F_1]_\mu)$. By $\omega_1 \rightarrow_* (\omega_1)^{\omega_1}_2$ (the correct type strong partition property), there is a club $C \subseteq \omega_1$ which is homogeneous for $P$ (for all functions $F \in [C]^{\mathcal{T}^2}_*$ of correct type).

Suppose $C$ is homogeneous for $P$ taking value $1$. Fix some $\alpha_0 \in [C]^{\omega_1}_* \slash \mu$. Let $f_0 \in [C]^{\omega_1}_*$ with $[f_0]_\mu = \alpha_0$. Pick $\beta < \omega_2$ with $\alpha_0 < \beta$. Since $[C]^{\omega_1}_* \slash \mu$ is cofinal through $\omega_2$, pick some $\alpha_1 \in [C]^{\omega_1}_* \slash \mu$ with $\beta < \alpha_1$. Let $f_1 \in [C]^{\omega_1}_*$ be such that $[f_1]_\mu = \alpha_1$. By an argument as in Fact \ref{presliding lemma} with some additional care to maintain the correct type, there is a function $F \in [C]^{\mathcal{T}^2}_*$ so that $[F_0]_\mu = \alpha_0$ and $[F_1]_\mu = \alpha_1 > \beta$. Since $P(F) = 1$, there is no element of $D$ between $\alpha_0$ and $\alpha_1$. In particular, there are no elements of $D$ between $\alpha_0$ and $\beta$. However since $\beta$ is arbitrary, this implies $D \subseteq \alpha_0$. This contradicts that $D$ is unbounded.

This shows that $C$ is homogeneous for $P$ taking value $0$. Let $\tilde C = \{\alpha \in C : \enum_C(\alpha) = \alpha\}$ which is a club subset of $\omega_1$. Let $f \in [\tilde C]^{\omega_1}$. Suppose $h <_\mu f$. For each $\alpha < \omega_1$, let $\gamma_\alpha$ be least so $h(\alpha) < \enum_C(\gamma_\alpha)$. Let $f_0(\alpha) = \enum_C(\gamma_\alpha + \omega)$ and $f_1 = \enum_C(\gamma_\alpha + \omega + \omega)$. Since $f(\alpha) \in \tilde C$, $f_0(\alpha) < f_1(\alpha) <  f(\alpha)$. Let $g_0(\alpha, n) = \enum_C(\gamma_\alpha + n)$ and $g_1(\alpha, n) = \enum_C(\gamma_\alpha + \omega + n)$. Note that $h <_\mu f_0 <_\mu f_1 <_\mu f$, $f_0, f_1 \in [C]^{\omega_1}_*$, and $g_0,g_1$ witnesses that $f_0,f_1$ are functions of the correct type. As above, one can find some $F \in [C]^{\mathcal{T}^2}_*$ so that $F_0 =_\mu f_0$ and $F_1 =_\mu f_1$. Since $P(F) = 0$, there is some $\alpha \in D$ so that $[h]_\mu < [f_0]_\mu < \alpha < [f_1]_\mu < [f]_\mu$. Since $D$ is a club and $h<_\mu f$ was arbitrary, this shows that $[f]_\mu \in D$. 

It has been shown that for all $D \subseteq \omega_1$, there is some club $\tilde C \subseteq \omega_1$ so that $[\tilde C]^{\omega_1}\slash \mu \subseteq D$. 
\end{proof}

\section{Failure of the Strong Partition Property at $\omega_2$}\label{failure SP at omega2}

\Begin{fact}{representation of omega cofinal points in omega2}
Assume $\ZF + \AD$. Let $\mu$ be the club measure on $\omega_1$. Suppose $f \in {}^{\omega_1}\omega_1$ is a function of uniform cofinality $\omega$. Then $\mathrm{cof}([f]_\mu) = \omega$.
\end{fact}

\begin{proof}
Suppose $f$ has uniform cofinality $\omega$. Let $g : \omega_1 \times \omega \rightarrow \omega_1$ be such that for all $\alpha < \omega_1$, $f(\alpha) = \sup\{g(\alpha, n) : n \in \omega\}$. Let $f_n : \omega_1 \rightarrow \omega_1$ be defined by $f_n(\alpha) = g(\alpha,n)$. Then $m < n$ implies $f_m <_\mu f_n$. Suppose $h <_\mu f$. Let $A = \{\alpha : h(\alpha) < f(\alpha)\}$. Define $k(\alpha)$ to be the least $n \in \omega$ so that $h(\alpha) < f_n(\alpha)$ whenever $\alpha \in A$. By the countable additivity of $\mu$, there is some $n^*$ so that $k(\alpha)=n^*$ for $\mu$-almost all $\alpha \in \kappa$. Then $h <_\mu f_{n^*}$. Thus $[f]_\mu = \sup\{[f_n]_\mu : n \in \omega\}$. So $\mathrm{cof}([f]_\mu) = \omega$.
\end{proof}

\Begin{fact}{representation of omega club subset omega2}
Assume $\ZF + \AD$. Let $\mu$ be the club measure on $\omega_1$. Let $C \subseteq \omega_1$ be a club. Then $[C]^{\omega_1}_* \slash \mu$ is an $\omega$-club subset of $\omega_2$. Moreover, for every $\omega$-club $D \subseteq \omega_2$, there is a club $C \subseteq \omega_1$ so that $[C]^{\omega_1}_* \slash \mu \subseteq D$. 
\end{fact}

\begin{proof}
It is clear that $[C]^{\omega_1}_* \slash \mu$ is unbounded. 

Suppose $\langle \nu_n : n \in \omega\rangle$ is an increasing $\omega$-sequence in $[C]^{\omega_1}_* \slash \mu$. By $\AC_\reals^\omega$, let $\langle g_n : n \in \omega\rangle$ be such that $\nu_n = [g_n]_\mu$ and each $g_n$ is of the correct type. Using $\AC_\omega^\reals$, one can also select a sequence $\langle k_n : n \in \omega\rangle$ so that $k_n$ witnesses that $g_n$ has uniform cofinality $\omega$. By Fact \ref{normal measure and closure points}, one may assume $g_n : \omega_1 \rightarrow C$. Using Lemma \ref{correct type sliding lemma}, let $F : \mathcal{T}^\omega \rightarrow C$ be an order preserving map of the correct type so that for all $n \in \omega$, $F^n =_\mu g_n$. Let $g : \omega_1 \rightarrow C$ be defined by $g(\alpha) = \sup\{F^n(\alpha) : n \in \omega\}$. Note that $g$ is of the correct type. As before, one can check that $[g]_\mu = \sup\{\nu_n : n \in \omega\}$ and $[g]_\nu \in [C]^{\omega_1}_* \slash \mu$. This shows that $[C]^{\omega_1}_* \slash \mu$ is an $\omega$-club subset of $\omega_2$.

Let $D \subseteq \omega_2$ be an $\omega$-club subset of $\omega_2$. Define $P : [\omega_1]^{\omega_1}_* \rightarrow 2$ by
$$P(f) = \begin{cases}
0 & \quad [f]_\mu \notin D \\
1 & \quad [f]_\mu \in D 
\end{cases}$$

By the correct type partition property $\omega_1 \rightarrow_* (\omega_1)^{\omega_1}_*$, there is a club $C \subseteq \omega_1$ which is homogeneous for $P$ (in the correct type sense). Since it was shown that $[C]^{\omega_1}_* \slash \mu$ is an $\omega$-club, $([C]^{\omega_1}_* \slash \mu) \cap D \neq \emptyset$. Thus $C$ must be homogeneous for $D$ taking value $1$. It has been shown that every $\omega$-club $D \subseteq \omega_2$ contains an $\omega$-club of the form $[C]^{\omega_1}_* \slash \mu$.
\end{proof}

\Begin{fact}{correct type function through club and its witness}
Assume $\ZF + \AD$. Let $\mu$ denote the club measure on $\omega_1$. Let $C \subseteq \omega_1$ be club. Let $B = [C]^{\omega_1}_* \slash \mu$ which is an $\omega$-club subset of $\omega_2$. 

Let $\epsilon < \omega_2$. Let $h : \omega_1 \rightarrow \omega_1$ with $h(\alpha) > 0$ for all $\alpha < \omega_1$ and $[h]_\mu = \epsilon$. Let $\Xi$ be a Kunen function for $h$.

Let $\mathcal{F} \in [B]^{\epsilon}_*$ (be of correct type). Let $\mathcal{G} : [\omega_2]^{\epsilon} \times \omega$ be a function witnessing that $\mathcal{F}$ has uniform cofinality $\omega$ with the property that for all $\alpha_0 < \alpha_1 < \epsilon$ and $m,n \in \omega$, $\mathcal{G}(\alpha_0,m) < \mathcal{G}(\alpha_1,n)$.

Then there is a sequence $\langle G_n : n \in \omega \rangle$  with each $G_n: \mathcal{T}^{h} \rightarrow \omega_1$ so that $\mathcal{G}(\alpha,n) = [G_n^{(\alpha)}]_\mu$.

There is an $F \in [C]^{\mathcal{T}^{h}}_*$ so that for all $\alpha < \epsilon$, $[F^{(\alpha)}]_\mu = \mathcal{F}(\alpha)$.
\end{fact}

\begin{proof}
For each $n$, let $\mathcal{G}_n : \epsilon \rightarrow \omega_1$ be defined by $\mathcal{G}_n(\alpha) = \mathcal{G}(\alpha,n)$. Using Fact \ref{representation alpha increasing ordinal sequence ultrapower} on $h$, $\Xi$, and $\mathcal{G}_n$, one obtains $G_n$.

Fix $\gamma < \epsilon$. Define $g'_{\gamma,n} : \omega_1 \rightarrow \omega_1$ by defined by
$$g'_{\gamma,n}(\alpha) = \begin{cases}
G_n^{(\gamma)}(\alpha) & \quad \Xi^{\pi_h^\Xi(\gamma)}(\alpha) < h(\alpha)\\
0 & \quad \text{otherwise}
\end{cases}$$
Let $f'_\gamma : \omega_1 \rightarrow \omega_1$ be defined by $f'_\gamma(\alpha) = \sup\{g'_{\gamma,n}(\alpha) : n \in \omega\}$. Note that $\mathcal{F}(\gamma) = [f'_\gamma]_{\mu}$.

Claim 1: The set $D_\gamma$ of $\alpha < \omega_1$ so that $f'_\gamma$ is discontinuous at $\alpha$ belongs to $\mu$.

Suppose $\omega_1 \setminus D_\gamma \in \mu$ and hence contains a club $E_\gamma$. Since $\mathcal{F}(\gamma) \in [C]^{\omega_1}_* \slash \mu$, there is some $\bar{f}_\gamma \in [C]^{\omega_1}_*$ (of the correct type) so that $f'_\gamma =_\mu \bar{f}_\gamma$. Let $J_\gamma = \{\alpha : \bar{f}_\gamma(\alpha) = f'_\gamma(\alpha)\} \in \mu$. Let $K_\gamma \subseteq J_\gamma$ be a club. Then $E_\gamma \cap K_\gamma$ is a club subset of $\omega_1$. Let $\langle\lambda_n : n \in \omega\rangle$ be an increasing sequence in $E_\gamma \cap K_\gamma$. Then $\lambda = \sup\{\lambda_n : n \in \omega\rangle \in E_\gamma \cap K_\gamma$ since $E_\gamma \cap K_\gamma$ is a club. Since $\lambda_n$ and $\lambda$ belong to $K_\gamma \subseteq J_\gamma$, $\bar{f}_\gamma(\lambda_n) = f'_\gamma(\lambda_n)$ and $\bar{f}_\gamma(\lambda) = f'_\gamma(\lambda)$. Since $\lambda \in E_\gamma$, one has that $f'_\gamma$ is continuous at $\lambda$. However, the previous statement implies that $\bar{f}_\gamma$ is also continuous at $\lambda$. But $\bar{f}_\gamma$ is discontinuous everywhere since it is of the correct type. This proves the claim.

Since $[f'_\gamma]_\mu = \mathcal{F}(\gamma) \in [C]^{\omega_1}_*\slash \mu$, there is some $\bar{f}_\gamma \in [C]^{\omega_1}_*$ so that $f'_\gamma =_\mu \bar{f}_\gamma$. Thus $P_\gamma = \{\alpha \in \omega_1: f_\gamma'(\alpha) \in C\} \in \mu$. Then $Q_\gamma = P_\gamma \cap D_\gamma \in \mu$. Let $f_\gamma = f'_\gamma \circ \enum_{Q_\gamma}$. By Fact \ref{normal measure and closure points}, $f'_\gamma =_\mu f_\gamma$. Thus $[f_\gamma]_\mu = \mathcal{F}(\gamma)$. Note that $f_\gamma$ is discontinuous everywhere. Define $g_\gamma : \omega_1 \times \omega \rightarrow \omega_1$ by $g_\gamma(\alpha,n) = g'_{\gamma,n}(\enum_{Q_\gamma}(\alpha))$. $g_\gamma$ witnesses that $f_\gamma$ has uniform cofinality $\omega$. Thus $f_\gamma$ is a function of the correct type.

By the uniformity of the construction, one has now produced a sequence $\langle f_\gamma : \gamma < \epsilon\rangle$ of functions of the correct type and $\langle g_\gamma : \gamma < \epsilon\rangle$ so that $[f_\gamma]_\mu = \mathcal{F}(\gamma)$ and $g_\gamma$ witnesses the uniform cofinality of $f_\gamma$. Now apply Lemma \ref{correct type sliding lemma} to $h$, $\Xi$, $\langle f_\gamma : \gamma < \epsilon\rangle$, and $\langle g_\gamma : \gamma < \epsilon\rangle$ to obtained the desired function $F \in [C]^{\mathcal{T}^{h}}_*$. 
\end{proof}

\Begin{fact}{uniform cofinality id implies cofinality is omega1}
Assume $\ZF + \AD$. Let $\mu$ be the club measure on $\omega_1$. Suppose $f$ has uniform cofinality $\mathsf{id}$. Then $[f]_\mu$ (as an ordinal in $\omega_2$) has cofinality $\omega_1$.
\end{fact}

\begin{proof}
Let $T^\mathsf{id} = \{(\alpha,\beta) : \beta < \alpha\}$. Let $g : T^\mathsf{id} \rightarrow \omega_1$ witness that $f$ has uniform cofinality $\mathsf{id}$. For each $\beta < \omega_1$, let 
$$f_\beta(\gamma) = \begin{cases}
g(\gamma,\beta) & \quad \gamma > \beta \\
0 & \quad \text{otherwise}
\end{cases}$$

One can show that $\langle [f_\beta]_\mu : \beta < \omega_1\rangle$ is an increasing sequence cofinal below $[f]_\mu$.
\end{proof}

\Begin{fact}{uniform cofinality omega or id}
Assume $\ZF + \AD$. Let $\mu$ denote the club measure on $\omega_1$. Let $f : \omega_1 \rightarrow \omega_1$ be a function so that $f(\alpha)$ is a limit ordinal for $\mu$-almost all $\alpha$. Then either $f$ has uniform cofinality $\omega$ or uniform cofinality $\mathsf{id}$ $\mu$-almost everywhere but not both.
\end{fact}

\begin{proof}
Let $\Xi$ be a Kunen function for $f$. Let $A = \{\alpha : f(\alpha) \in \mathrm{Lim}\} \in \mu$. By Definition \ref{kunen function definition}, for all $\alpha \in K^\Xi_f$, $f(\alpha) < \delta_\alpha^\Xi = \sup\{\Xi(\alpha,\beta) : \beta < \alpha\}$. 

Fix $\alpha \in K^\Xi_f \cap A$, define $\epsilon_0^\alpha$ be the least $\epsilon < \alpha$ so that $\Xi(\alpha,\epsilon) < f(\alpha)$. Suppose $\beta < \alpha$ and $\epsilon^\alpha_{\beta'}$ has been defined for all $\beta' < \beta < \alpha$. If $\sup\{\Xi_\alpha(\epsilon^\alpha_{\beta'}) : \beta' < \beta\} = f(\alpha)$, then let $\alpha^*$ denote this $\beta$ and declare the construction to have terminated. Otherwise, let $\epsilon^\alpha_\beta$ be the least $\epsilon < \alpha$ so that $\sup\{\Xi_\alpha(\epsilon^\alpha_{\beta'}) : \beta' < \beta\} < \Xi_\alpha(\epsilon) < f(\alpha)$. This defines a sequence $\langle \epsilon^\alpha_\beta : \beta < \alpha^*\rangle$ where $\alpha^* \leq \alpha$. Note that $\langle \epsilon^\alpha_\beta : \beta < \alpha^*\rangle$ is an increasing sequence and $\sup\{\Xi_\alpha(\epsilon_\beta^\alpha) : \beta < \alpha^*\} = f(\alpha)$.

Consider the function $\Phi : K_f^\Xi \cap A \rightarrow \omega_1$ by $\Phi(\alpha) = \alpha^*$. Let $B = \{\alpha \in A \cap K_f^\Xi : \Phi(\alpha) = \alpha^* < \alpha\}$. 

(Case I) $B \in \mu$. 

Then $\Phi$ is a regressive function. Since $\mu$ is normal, there is a $C \subseteq B$ and $\delta < \omega_1$ so that $C \in \mu$ and $\Phi(\alpha) = \delta$ for all $\alpha \in C$. Let $\phi : \omega \rightarrow \delta$ be a cofinal sequence. Let $g : \omega_1 \times \omega \rightarrow \omega_1$ be defined by 
$$g(\alpha,n) = \begin{cases}
n & \quad \alpha \notin C \\
\Xi_\alpha(\epsilon^\alpha_{\phi(n)}) & \quad \alpha \in C
\end{cases}$$

The function $g$ witnesses that $f$ is $\mu$-almost everywhere a function of uniform cofinality $\omega$.

(Case II) $B \notin \mu$.

Let $C = (K_f^\Xi \cap A)\setminus B$. For all $\alpha \in C$, $\Phi(\alpha) = \alpha^* = \alpha$. Let $T^\mathsf{id} = \{(\alpha,\beta) : \alpha < \omega_1 \wedge \beta < \alpha\}$. Define $g : T^\mathsf{id} \rightarrow \omega_1$ by
$$g(\alpha,\beta) = \begin{cases}
\beta & \quad \alpha \notin C \\
\Xi_\alpha(\epsilon^\alpha_\beta) & \quad \alpha \in C
\end{cases}$$
Then $g$ witnesses that $f$ has uniform cofinality $\mathsf{id}$ $\mu$-almost everywhere.

Fact \ref{representation of omega cofinal points in omega2} implies that if $f$ has uniform cofinality $\omega$, then $[f]_\mu$ has cofinality $\omega$. Fact \ref{uniform cofinality id implies cofinality is omega1} implies that if $f$ has uniform cofinality $\mathsf{id}$, then $[f]_\mu$ has cofinality $\omega_1$. Since $\omega_1$ is regular under $\AD$, $f$ cannot have both uniform cofinality.
\end{proof}

In the previous result, it is shown that $f$ cannot have uniform cofinality $\omega$ and $\mathsf{id}$ by showing these uniform cofinalities correspond to different cofinalities of $[f]_\mu$. This is however not the case in general. The following argument avoids considering the function in the ultrapower. This argument can be generalized to show functions from $f : [\omega_1]^n \rightarrow \omega_1$ can only have one uniform cofinality. 

\Begin{fact}{not both uniform cofinality on omega1}
Assume $\ZF + \AD$. Every function $f : \omega_1 \rightarrow \omega_1$ has uniform cofinality $\omega$ or $\mathsf{id}$ $\mu$-almost everywhere but not both.
\end{fact}

\begin{proof}
This is already proved in Fact \ref{uniform cofinality omega or id}, so one will just give an argument that no function can have both uniform cofinality which does not involve the ultrapower.

Suppose that $f : \omega_1 \rightarrow \omega_1$ $\mu$-almost everywhere has both uniform cofinality $\omega$ and $\mathsf{id}$. Let $g_\omega : \omega_1 \times \omega \rightarrow \omega_1$ witness that $f$ has uniform cofinality $\omega$ $\mu$-almost everywhere. Let $g_{\mathsf{id}} : \omega_1 \times \omega_1 \rightarrow \omega_1$ witness that $f$ has uniform cofinality $\mathsf{id}$ $\mu$-almost everywhere.

Let $A_\omega \in \mu$ be the set of $\alpha < \omega_1$ so that $f(\alpha) = \sup\{g_\omega(\alpha,n) : n \in \omega\}$. Let $A_{\mathsf{id}} \in \mu$ be the set of $\alpha < \omega_1$ so that $f(\alpha) = \sup\{g_\mathsf{id}(\alpha,\beta) : \beta < \alpha\}$. Let $A = A_\omega \cap A_{\mathsf{id}}$ which also belongs to $\mu$.

For each $\beta \in A$ and $\alpha < \beta$, let $n_{\alpha,\beta}$ be least $n \in \omega$ so that $g_\omega(\beta,n) > g_{\mathsf{id}}(\beta,\alpha)$ which exists since $\beta \in A$. Consider the function $\Phi : [A]^2 \rightarrow \omega$ defined by $\Phi(\alpha,\beta) = n_{\alpha,\beta}$. By the weak partition property and the countable additivity of $\mu$, one has that there is some $B \subseteq A$ with $B \in \mu$ and  $n \in \omega$ so that for all $(\alpha,\beta) \in [B]^2$, $n_{\alpha,\beta} = n$. 

Now let $\beta \in B$ be a limit point of $B$, that is $\sup B\cap \beta = \beta$. For each $\alpha < \beta$ with $\alpha \in B$, one has that $n_{\alpha,\beta} = n$. Then $f(\beta) = \sup\{g_{\mathsf{id}}(\beta,\alpha) : \alpha \in \beta\} \leq g_\omega(\beta,n) < f(\beta)$. Contradiction.
\end{proof}

Next one will consider the ultrapower, $\mathrm{ult}(V,\mu)$, where $\mu$ is the club measure on $\omega_1$. One needs some care when working with this structure as one may not have \Los's theorem without $\AC$. It will be shown later that \Los's theorem fails for this ultrapower and in fact $\mathrm{ult}(V,\mu)$ is not a model $\ZF$.

\Begin{fact}{representation of a subset of omega2}
Assume $\ZF + \AD$. Let $\mu$ denote the club measure on $\omega_1$. Suppose $f : \omega_1 \rightarrow V$ is such that $\mathrm{ult}(V,\mu) \models [f]_\mu \subseteq \omega_2$, then there is an $f' : \omega_1 \rightarrow \powerset{\omega_1}$ so that $[f']_\mu = [f]_\mu$.
\end{fact}

\begin{proof}
By Corollary \ref{value of omega1 ultrapower}, one knows that $\omega_2 = \prod_{\omega_1} \omega_1 \slash \mu$. Thus in $\mathrm{ult}(V,\mu)$, each $\zeta < \omega_2$ is represented by some $h : \omega_1 \rightarrow \omega_1$. 

Now let $f'(\alpha) = f(\alpha)\cap \omega_1$. The claim is that $[f']_\mu = [f]_\mu$.

Suppose $h : \omega_1 \rightarrow V$ is a function so that $\mathrm{ult}(V,\mu) \models [h]_\mu \in [f]_\mu$. Since $\mathrm{ult}(V,\mu) \models [f]_\mu \subseteq \omega_2$, one must have that $A = \{\alpha \in \omega_1 : h(\alpha) \in f(\alpha) \cap \omega_1\} \in \mu$. Let $h' : \omega_1 \rightarrow \omega_1$ be defined by $h'(\alpha) = h(\alpha)$ if $\alpha \in A$ and $h'(\alpha) = \emptyset$ otherwise. Note that $[h]_\mu = [h']_\mu$ and $[h']_\mu \in [f']_\mu$. This shows that $[f]_\mu \subseteq [f']_\mu$. It is clear that $[f']_\mu \subseteq [f]_\mu$. Thus $[f]_\mu = [f']_\mu$.
\end{proof}

\Begin{fact}{uniform witness implies enumeration correct type}
Assume $\ZF + \AD$ and $\mu$ is the club measure on $\omega_1$. Suppose $h : \omega_1 \rightarrow \powerset{\omega_1}$ has the property that for all $\alpha < \omega_1$, $|h(\alpha)| = \omega_1$ and $\enum_{h(\alpha)}$ is a function of correct type. Suppose there is a function $G$ so that for each $\alpha < \omega_1$, $G(\alpha) : \omega_1\times\omega\rightarrow\omega_1$ is a witness to $\enum_{h(\alpha)}$ having uniform cofinality $\omega$, then $\enum_{[h]_\mu} : \omega_2 \rightarrow \omega_2$ is a function of correct type.
\end{fact}

\begin{proof}
For each $\xi < \omega_1$, let $g_\xi : \omega_1 \times \omega \rightarrow \omega_1$ be defined by $g_\xi = G(\xi)$. Hence $g_\xi$ is a witnesses to $\enum_{h(\alpha)}$ having the correct type.

For each $\xi < \omega_2$, let $p_\xi : \omega_1 \rightarrow \omega_1$ have the property that for all $\alpha < \omega_1$, $p_\xi(\alpha) \in h(\alpha)$ and $[p_\xi]_\mu$ represents the $\xi^\text{th}$-element of $[h]_\mu$. (Note that one does not have a uniform procedure for finding $p_\xi$ as $\xi$ ranges over ordinals below $\omega_2$.) Define $k_{\xi,n} : \omega_1 \rightarrow \omega_1$ by $k_{\xi,n}(\alpha) = g_\xi(\enum_{h(\xi)}^{-1}(p_\xi(\alpha)),n)$. Let $\delta_{\xi,n} = [k_{\xi,n}]_\mu$. Note that if $p'_\xi =_\mu p_\xi$ and $k'_{\xi,n}$ was defined in the same manner as $k$ using $p'_\xi$ instead of $p_\xi$, then $k'_{\xi,n} =_\mu k_{\xi,n}$. This shows that $\delta_{\xi,n}$ is well defined independent of the choice of $p_\xi$. 

Define $r : \omega_2 \times \omega \rightarrow \omega_2$ by $r(\xi,n) = \delta_{\xi,n}$. One can check that $r$ witnesses that $\enum_{[h]_\mu}$ has uniform cofinality $\omega$. 

Let $\zeta \in [h]_\mu$. Let $\ell : \omega_1 \rightarrow \omega_1$ be such that $[\ell]_\mu = \zeta$. One may assume that for all $\alpha$, $\ell(\alpha) \in h(\alpha)$. Let $\iota : \omega_1 \rightarrow \omega_1$ be defined by $\iota(\alpha) = \sup(h(\alpha) \cap \ell(\alpha))$. Note that $\iota(\alpha) < \ell(\alpha)$ since $\enum_{h(\alpha)}$ was assumed to be discontinuous everywhere. One can check that every element of $[h]_\mu$ which is below $[\ell]_\mu$ is below $[\iota]_\mu < [\ell]_\mu$. Thus $\enum_{[h]_\mu}$ is discontinuous everywhere.
\end{proof}

\Begin{fact}{enumeration of the elements of h}
Suppose $h : \omega_1 \rightarrow \powerset{\omega_1}$ with $|h(\alpha)| = \omega_1$ for all $\alpha < \omega_1$. Let $\xi < \omega_2$ and $\ell : \omega_1 \rightarrow \omega_1$ be such that $[\ell]_\mu = \xi$. $\enum_{[h]_\mu}(\xi)$ is represented by the function $g : \omega_1 \rightarrow \omega_1$ defined by $g(\alpha) = \enum_{h(\alpha)}(\ell(\alpha))$. 
\end{fact}

\begin{proof}
For each $\ell : \omega_1 \rightarrow \omega_1$, let $g_\ell : \omega_1 \rightarrow \omega_1$ be defined by $g_\ell(\alpha) = \enum_{h(\alpha)}(\ell(\alpha))$. Note that if $\ell =_\mu \ell'$, then $g_\ell =_\mu g_{\ell'}$. Note also that $[g_\ell]_\mu \in [h]_\mu$. 

For each $\xi < \omega_1$, let $\gamma_\xi \in \omega_2$ be defined by $\gamma_\xi = [g_\ell]_\mu$ where $\ell$ is any function so that $[\ell]_\mu = \xi$. This is well defined by the previous paragraph. $\langle \gamma_\xi : \xi < \omega_2\rangle$ is an increasing sequence through $[h]_\mu$. 

Let $k : \omega_1 \rightarrow \omega_1$ be such that $[k]_\mu \in [h]_\mu$. Note $A = \{\alpha \in \omega_1 : k(\alpha) \in h(\alpha)\} \in \mu$. By modifying $k$ off $A$, one will assume without loss of generality that $k(\alpha) \in h(\alpha)$ for all $\alpha < \omega_1$. Let $\ell(\alpha) = \enum_{h(\alpha)}^{-1}(k(\alpha))$. Then $k =_\mu g_{\ell}$. 

This shows that $\enum_{[h]_\mu}(\xi) = \gamma_\xi$. This completes the proof.
\end{proof}

\Begin{fact}{omega1 size subset of omega2 in ultrapower}
Assume $\ZF + \AD$ and let $\mu$ be the club measure on $\omega_1$. Suppose $E \subseteq \omega_2$ is such that $|E| = \omega_1$. Then $E \in \mathrm{ult}(V,\mu)$.
\end{fact}

\begin{proof}
Let $\zeta = \ot(E)$. Let $h : \omega_1 \rightarrow \omega_1$ be such that $[h]_\mu = \zeta$ and $h(\alpha) > 0$ for all $\alpha < \omega_1$. Let $\Xi$ be a Kunen function for $h$. Using Fact \ref{representation alpha increasing ordinal sequence ultrapower}, there is an $F : \mathcal{T}^h \rightarrow \omega_1$ which is increasing so that for all $\xi < \zeta$, $[F^{(\xi)}]_\mu = \enum_E(\xi)$. 

Let $g : \omega_1 \rightarrow \powerset{\omega_1}$ be defined by $g(\alpha) = \{F(\alpha,\beta) : \beta < h(\alpha)\}$. The claim is that $[g]_\mu = E$. 

Suppose $p : \omega_1 \rightarrow \omega_1$ is such that $[p]_\mu \in [g]_\mu$. Note $A = \{\alpha : p(\alpha) \in g(\alpha)\} \in \mu$. By modifying $p$ off $A$, one may assume that for all $\alpha$, $p(\alpha) \in g(\alpha)$. Let $f(\alpha)$ be the $\beta < h(\alpha)$ so that $F(\alpha,\beta) = p(\alpha)$. Then $p = F^f = F^{(\xi)} = \enum_E(\xi)$ for some $\xi < \zeta$. Then $[g]_\mu \subseteq E$. It is straightforward to see that $E \subseteq [g]_\mu$. 
\end{proof}

\Begin{fact}{adding countable set below the sup}
Assume $\ZF + \AD$. Suppose $D \subseteq \omega_2$ is such that there is some $g : \omega_1 \rightarrow \powerset{\omega_1}$ so that $[g]_\mu = D$. Suppose $E \subseteq \omega_2$ is such that $\sup E < \min D$, then $D \cup E \in \mathrm{ult}(V,\mu)$.
\end{fact}

\begin{proof}
By Fact \ref{omega1 size subset of omega2 in ultrapower}, $E \in \mathrm{ult}(V,\mu)$. Thus there is a $g' : \omega_1 \rightarrow \powerset{\omega_1}$ so that $[g']_\mu = E$. Let $g'' : \omega_1 \rightarrow \powerset{\omega_1}$ be defined by $g''(\alpha) = g(\alpha) \cup g'(\alpha)$. One can check that $[g'']_\mu = D \cup E$.
\end{proof}

\Begin{fact}{subset of omega2 not in ultrapower}
Assume $\ZF + \AD$ and let $\mu$ be the club measure on $\omega_1$. Let $D = \{\alpha < \omega_2 : \mathrm{cof}(\alpha) = \omega_1\}$. Then $D \notin \mathrm{ult}(V,\mu)$.
\end{fact}

\begin{proof}
Suppose $D \in \mathrm{ult}(V,\mu)$. By Fact \ref{representation of a subset of omega2}, there is some $h : \omega_1 \rightarrow \powerset{\omega_1}$ so that $D = [h]_\mu$. 

(Case I) $A = \{\alpha < \omega_1 : h(\alpha) \in \mu\} \in \mu$.

By Fact \ref{omega1 club selection}, there is a sequence $\langle C_\alpha : \alpha \in A\rangle$ of club subsets of $\omega_1$ so that $C_\alpha \subseteq h(\alpha)$ for all $\alpha \in A$. 

Define $\ell : \omega_1 \rightarrow \omega_1$ by $\ell(\alpha) = \enum_{C_{\enum_A(\alpha)}}(\omega)$. Let $\ell' : \omega_1 \rightarrow \omega_1$ be defined by 
$$\ell'(\alpha) = \begin{cases}
\enum_{C_\alpha}(\omega) & \quad \alpha \in A \\
\min(h(\alpha)) & \quad \alpha \notin A
\end{cases}$$
Clearly $[\ell']_\mu \in [h]_\mu$. Note that $\ell = \ell' \circ \enum_A$. By Fact \ref{normal measure and closure points}, $\ell =_\mu \ell'$. Define $g : \omega_1 \times \omega \rightarrow \omega$ by $g(\alpha,n) = \enum_{C_{\enum_A(\alpha)}}(n)$. The function $g$ witnesses that $\ell$ has uniform cofinality $\omega$. Thus $\mathrm{cof}([\ell]_\mu) = \omega$ by Fact \ref{representation of omega cofinal points in omega2}. Since $[\ell]_\mu \in D = \{\alpha < \omega_2 : \mathrm{cof}(\alpha) = \omega_1\}$, one has a contradiction.

(Case II) $A = \{\alpha < \omega_1 : h(\alpha) \in \mu\} \notin \mu$. 

Then $B = \omega_1 \setminus A = \{\alpha < \omega_1 : \omega_1 \setminus h(\alpha) \in \mu\} \in \mu$. By Fact \ref{omega1 club selection}, there is a sequence $\langle C_\alpha : \alpha \in B\rangle$ of club subsets of $\omega_1$ so that $C_\alpha \subseteq \omega_1 \setminus h(\alpha)$. Define $\ell : \omega_1 \rightarrow \omega_1$ by $\ell(\alpha) = \enum_{C_{\enum_B(\alpha)}}(\alpha)$. Let $\ell' : \omega_1 \rightarrow \omega_1$ be defined by
$$\ell'(\alpha) = \begin{cases}
\enum_{C_\alpha}(\alpha) & \quad \alpha \in B \\
\min h(\alpha) & \quad \alpha \notin B
\end{cases}$$
Note that $[\ell']_\mu \notin [h]_\mu$. Observe $\ell = \ell \circ \enum_B$. By Fact \ref{normal measure and closure points}, $\ell =_\mu \ell'$. Let $g: \mathcal{T}^{\mathsf{id}} \rightarrow \omega_1$ be defined by $g(\alpha,\beta) = \enum_{C_{\enum_{B}(\alpha)}}(\beta)$ if $\beta < \alpha$. Then $g$ witnesses that $\ell$ has uniform cofinality $\mathrm{id}$. Therefore, $\mathrm{cof}([\ell]_\mu) = \omega_1$ by Fact \ref{uniform cofinality id implies cofinality is omega1}. However, $[\ell]_\mu \notin D = \{\alpha < \omega_2 : \mathrm{cof}(\alpha) = \omega_1\}$ yields a contradiction. 

The proof is complete.
\end{proof}

One can now show that $\mu$ does not satisfy \Los's theorem and in fact the ultrapower does not satisfy the $\ZF$ axioms.

\Begin{fact}{los theorem fails for club measure}
Assume $\ZF + \AD + V = L(\reals)$. Let $\mu$ denote the club measure on $\omega_1$. Then $\mathrm{ult}(L(\reals),\mu)$ is not a model of $\ZF$. Thus \Los's theorem fails for $\mu$.
\end{fact}

\begin{proof}
Note that $L(\reals) \models \AD$ implies $L(\reals) \models \DC_\reals$ and hence $L(\reals) \models \DC$ by a result of Kechris \cite{The-Axiom-of-Determinacy-Implies-Dependent-Choice}. Thus $\mathrm{ult}(L(\reals),\mu)$ may be considered a transitive inner model of $L(\reals)$. One can check that $\reals \subseteq \mathrm{ult}(L(\reals),\mu)$. If $\mathrm{ult}(L(\reals),\mu)$ is an inner model of $\ZF$ containing all the reals of $L(\reals)$, then one must have that $\mathrm{ult}(L(\reals),\mu) = L(\reals)$. This is impossible since Fact \ref{subset of omega2 not in ultrapower} asserts that $\mathrm{ult}(L(\reals),\mu)$ is missing a subset of $\omega_2$ which belongs to $L(\reals)$. 
\end{proof}

\Begin{theorem}{example of failure of strong partition at omega2}
(Jackson) Assume $\ZF + \AD$. Let $\mu$ denote the club measure on $\omega_1$. Define a partition $P : [\omega_2]^{\omega_2}_* \rightarrow 2$ by
$$P(f) = 0 \Leftrightarrow \rang(f) \in \mathrm{ult}(V,\mu).$$

Then there is no club $D \subseteq \omega_2$ and no $i \in 2$ so that $P(f) = i$ for all $f \in [D]^{\omega_2}_*$. 
\end{theorem}

\begin{proof}
(Case I) Suppose there is a club $D \subseteq \omega_2$ so that $P(f) = 1$ for all $f \in [D]^{\omega_2}_*$. By Fact \ref{representation of club subset of omega2}, there is a club $C \subseteq \omega_1$ so that $[C]^{\omega_1}\slash \mu \subseteq D$.

Let $A = \{\alpha : (\exists \gamma)(\alpha = \enum_{C}(\gamma + \omega)\}$. Note that $\enum_A : \omega_1 \rightarrow \omega_1$ is a function of the correct type. Let $g : \omega_1 \times \omega \rightarrow \omega$ witness that $\enum_A$ has uniform cofinality $\omega$. For each $\xi < \omega_1$, let $A_\xi = \{\alpha \in A : \alpha \geq \xi\}$. Let $g_\xi : \omega_1 \times \omega \rightarrow \omega$ be defined by $g_\xi(\alpha,n) = g(\enum_A^{-1}(\enum_{A_\xi}(\alpha)),n)$. Then for each $\xi < \omega_1$, $g_\xi$ witnesses that $\enum_{A_\xi}$ has uniform cofinality $\omega$. Let $h : \omega_1 \rightarrow \powerset{\omega_1}$ be defined by $h(\xi) = A_{\xi}$. By Fact \ref{uniform witness implies enumeration correct type}, $\enum_{[h]_\mu} : \omega_2 \rightarrow \omega_2$ is a function of the correct type. Let $f = \enum_{[h]_\mu}$. One can check that $\rang(f) = [h]_\mu \subseteq [C]^{\omega_1} \slash \mu$. Thus $f \in [D]^{\omega_2}_*$. However since $\rang(f) = [h]_\mu \in \mathrm{ult}(V,\mu)$, one must have that $P(f) = 0$. Contradiction.

(Case II) Suppose there is a club $D \subseteq \omega_2$ so that $P(f) = 0$ for all $f \in [D]^{\omega_2}_*$.

Fix any $f \in [D]^{\omega_2}_*$. Since $P(f) = 0$, there is some $h : \omega_1 \rightarrow \powerset{\omega_1}$ so that $\rang(f) = [h]_\mu$ by Fact \ref{representation of a subset of omega2}.

Now let $A \subseteq \omega_2$ be arbitrary. $\enum_{f[A]} : \omega_2 \rightarrow D$ is a function of the correct type since $f : \omega_2 \rightarrow D$ is a function of the correct type. Then $P(\enum_{f[A]}) = 0$. Thus $f[A] \in \mathrm{ult}(V,\mu)$. There is some $g : \omega_1 \rightarrow \powerset{\omega_1}$ so that $[g]_\mu = f[A]$. Let $k : \omega_1 \rightarrow \powerset{\omega_1}$ be defined by $k(\alpha) = \enum_{h(\alpha)}^{-1}[g(\alpha)]$.

The claim is that $[k]_\mu = A$.

To show $[k]_\mu \subseteq A$: Suppose $\xi \in [k]_\mu$. There is a representative $\ell : \omega_1 \rightarrow \omega_1$ of $\xi$ so that for all $\mu$-almost all $\alpha$, $\ell(\alpha) \in k(\alpha)$. There is a function $p : \omega_1 \rightarrow \omega_1$ so that for $\mu$-almost all $\alpha$, $p(\alpha) \in g(\alpha)$ and $\ell(\alpha) = \enum_{h(\alpha)}^{-1}(p(\alpha))$. In particular $[p]_\mu \in [g]_\mu = f[A]$. By Fact \ref{enumeration of the elements of h}, the $\xi^\text{th}$ element of $[h]_\mu$ is represented by the function $q : \omega_1 \rightarrow \omega_1$ defined by $q(\alpha) = \enum_{h(\alpha)}(\ell(\alpha))$. However $\mu$-almost everywhere, $q(\alpha) = \enum_{h(\alpha)}(\enum_{h(\alpha)}^{-1}(p(\alpha))) = p(\alpha)$. Thus $q =_\mu p$. It has been shown that the $\xi^\text{th}$ element of $[h]_\mu = \rang(f)$ belongs to $[g]_\mu = f[A]$. Thus $\xi \in A$. 

To show that $A \subseteq [k]_\mu$: Let $\xi \in A$. Let $\ell : \omega_1 \rightarrow \omega_1$ be such that $[\ell]_\mu = \xi$. By Fact \ref{enumeration of the elements of h}, the $\xi^\text{th}$-element of $[h]$ is represented by the function $q(\alpha) = \enum_{h(\alpha)}(\ell(\alpha))$. Thus $[q]_\mu \in f[A] = [g]_\mu$. So $\mu$-almost everywhere, $q(\alpha) = \enum_{h(\alpha)}(\ell(\alpha)) \in g(\alpha)$. Thus for $\mu$-almost all $\alpha$, $\ell(\alpha) \in \enum_{h(\alpha)}^{-1}[g(\alpha)] = k(\alpha)$. It has been shown that $\xi = [\ell]_\mu \in [k]_\mu$. 

The claim has been shown. Since $A$ was arbitrary, this implies that everywhere subset of $\omega_2$ belongs to $\mathrm{ult}(V,\mu)$. This contradicts Fact \ref{subset of omega2 not in ultrapower}.

It has been shown that the partition $P$ has no club set which is homogeneous.
\end{proof}

\Begin{corollary}{martin paris failure SP omega2}
(Martin and Paris) Assume $\ZF + \AD$. The partition relation $\omega_2 \rightarrow (\omega_2)^{\omega_2}_2$ does not hold. Thus $\omega_2$ is a weak partition cardinal which is not a strong partition cardinal.
\end{corollary}

\Begin{remark}{failure SP on omega2 remark}
The example of Jackson from Theorem \ref{example of failure of strong partition at omega2} gives an explicit example of a partition $P : [\omega_2]^{\omega_2}_* \rightarrow 2$ which has no homogeneous club subset.

Martin and Paris original argument roughly shows that if $\omega_2 \rightarrow (\omega_2)^{\omega_2}_2$ holds, then $\omega_3$ would satisfy $\omega_3 \rightarrow (\omega_3)^{\alpha}_2$ for all $\alpha < \omega_1$. Fact \ref{partition implies regularity regularity} implies that $\omega_3$ must be regular. However, it can be shown that $\omega_3$ is a singular cardinal of cofinality $\omega_2$. For more information on the result of Martin and Paris, see \cite{Infinitary-Combinatorics-and-the-Axiom-of-Determinateness} and specifically Lemma 5.19. Also see Section 13 of \cite{AD-and-Projective-Ordinals}.
\end{remark}

\Begin{corollary}{bounded addition failure SP on omega2}
Let $\sigma \in [\omega_2]^{<\omega_2}_*$. Define $P_\sigma : [\omega_2 \setminus (\sup(\sigma) + \omega)]^{\omega_2}_* \rightarrow 2$ by $P(\sigma\hat{\ }f)$, where $P$ is the partition from Theorem \ref{example of failure of strong partition at omega2}. $P_\sigma$ also does not have have a club homogeneous set.
\end{corollary}

\begin{proof}
Essentially the same argument as Theorem \ref{example of failure of strong partition at omega2} with the assistance of Fact \ref{adding countable set below the sup}.
\end{proof}

\section{$L(\reals)$ as a Symmetric Collapse Extension of $\HOD$}\label{L(R) symmetric collapse extension HOD}

\Begin{definition}{infinity borel codes}
Let $S \subseteq \ON$ be a set of ordinals. Let $\varphi$ be a formula of set theory. The pair $(S,\varphi)$ is called an $\infty$-Borel code. Let $n \in \omega$. Define $\mathfrak{B}^n_{(S,\varphi)} = \{x \in \reals^n : L[S,x] \models \varphi(S,x)\}$.

Let $A \subseteq \reals^n$. $(S,\varphi)$ is said to be an $\infty$-Borel code for $A$ if and only if $A = \mathfrak{B}^n_{(S,\varphi)}$.
\end{definition}

If a set $A \subseteq \reals^n$ has an $\infty$-Borel code, then $A$ has a very absolute definition. That is, in ordered to determine membership of $x$ in $A$, one needs only to ask whether $\varphi(S,x)$ holds in $L[S,x]$, which is the minimal model of $\ZFC$ containing $x$ and the code set $S$.

\Begin{definition}{ordinal definability}
Let $P$ be a set. Recall a set $A$ is ordinal definable from $P$ if and only if there is a formula $\varphi$ of set theory, a finite tuple of ordinals $\bar{\alpha}$, and a finite tuple $\bar{p}$ of elements of $P$ so that $A = \{x : \varphi(x,\bar{\alpha},\bar{p})\}$. 

Using the reflection theorem, one can show that a set $A$ is ordinal definable if and only if there is some $\xi \in \ON$, tuple of ordinals $\bar{\alpha}$, tuple $\bar{p}$ from $P$, and formula $\varphi$ so that all these objects belong to $V_\xi$ and $A = \{x \in V_\xi : V_\xi \models \varphi(x,\bar{\alpha},\bar{p})\}$.

This shows the collection $\OD_P$ of all sets which are ordinal definable from $P$ forms a first order class. If $P$ has an $\OD_P$ wellordering, then $\OD_P$ has a wellordering which is definable with parameter from ordinals and $P$. Thus there is a bijection of $\OD_P$ with $\ON$.

Let $\HOD_P$ denote the subclass of $\OD_P$ which is hereditarily $\OD_P$. That is, $\HOD_P$ consists of those $x \in \OD_P$ such that $\mathrm{tc}(\{x\}) \subseteq \OD_P$, where $\mathrm{tc}$ refers to the transitive closure.

As a matter of convention, if $S \subseteq \ON$ is a set of ordinals, one will often write $\OD_S$ and $\HOD_S$ for $\OD_{\{S\}}$ and $\HOD_{\{S\}}$. 
\end{definition}

\Begin{definition}{Vopenka forcing}
Let $n \in \omega$ and $S \subseteq \ON$ be a set of ordinals. Let ${}_n\bbO_S$ denote the forcing of nonempty ordinal definable in $S$ subsets of $\reals^n$. Let the ordering be $\leq_{{}_n\bbO_S} = \subseteq$. The largest element is $1_{{}_n\bbO_S} = \reals^n$. 

Since there is a definable (in $S$) bijection of the class $\OD_S$ with the ordinal $\ON$, one can identify ${}_n\bbO_S$ as a set of ordinals in $\HOD_S$. In this way, ${}_n\bbO_S \in \HOD_S$. ${}_n\bbO_S$ is called the $n$-dimensional $S$-Vop\v{e}nka forcing. If $n = 1$, ${}_1\bbO_S$ will be denoted simply $\bbO_S$.
\end{definition}

\Begin{definition}{infinity borel code forcing}
Let $n \in \omega$ and $S \subseteq \ON$ be a set of ordinals. Let ${}_n\bbA_S$ denote the forcing of nonempty subsets of $\reals^n$ which possess $\OD_S$ $\infty$-Borel codes. ${}_n\bbA_S$ is ordered by $\leq_{{}_n\bbA_S} = \subseteq$. It has a largest element $1_{{}_n\bbO_S} = \reals^n$.

Since ${}_n\bbA_S \subseteq {}_n\bbO_S$, one can consider ${}_n\bbA_S$ to be a forcing of $\HOD_S$. ${}_n\bbA_S$ will be called the $n$-dimensional $\OD_S$ $\infty$-Borel code forcing.
\end{definition}

\Begin{remark}{coding the A forcing into the ordinals}
One can be more specific about how ${}_n\bbA_S$ is coded as a set of ordinals. One can identify ${}_n\bbA_S$ with a (set sized) collection of pairs of $(S',\varphi)$, where $S'$ is a $\OD_S$ set of ordinals and $\varphi$ is a formula. Using the canonical global wellordering of $\HOD_S$, let $\langle (S_\alpha, \varphi_\alpha) : \alpha < \delta\rangle$, for some ordinal $\delta$, be an enumeration of $\infty$-Borel codes that include at least one code for each element of ${}_n\bbA_S$. Fix $\langle \phi_n : n \in \omega\rangle$ to be a coding of formulas of set theory by natural numbers. Using the G\"odel pairing function, let $K = \{(\alpha,\beta,n) : \beta \in S_\alpha \wedge \phi_n = \varphi_\alpha\}$. One will often identify ${}_n\bbA_S$ with this set of ordinals $K$. If this is done, then from ${}_n\bbA_S$, one can obtain uniformly $\infty$-Borel codes for each condition in ${}_n\bbA_S$. 
\end{remark}

In nearly every regard, the Vop\v{e}nka forcing is a more practical forcing than $\bbA_S$. It will shown that in $\ZF + \AD + V = L(\reals)$, $\bbO_S$ and $\bbA_S$ are identical. To establish this, one will prove a structural theorem about $L(\reals)$ due to Woodin that involves the forcing $\bbA_S$. The presentation of Woodin's result that $L(\reals)$ is a symmetric collapse extension of its $\HOD$ follows closely \cite{Proper-Forcing-and-Absoluteness-Comm}. Of particular importance to the study of cardinals and combinatorics in $L(\reals) \models \AD$ will be existence of an ultimate $\infty$-Borel code which follows from the proof.

For simplicity, $S = \emptyset$ in the following results. The result can be appropriately relativized.

The main benefit of $\bbA$ over $\bbO$ is the following result:

\Begin{fact}{HOD bbO extension is same as extension by real}
Let $x \in \reals^n$. Then there is a generic filter $G_x \subseteq {}_n\bbA$ which is ${}_n\bbA$-generic over $\HOD$ so that $\HOD[G_x] = \HOD[x]$, where $\HOD[x]$ refers to the smallest transitive model of $\ZF$ extending $\HOD$ and containing $x$.
\end{fact}

\begin{proof}
For simplicity, let $n = 1$. Let $x \in \reals$. Let $G_x = \{p \in \bbA : x \in p\}$.

Using the convention of Remark \ref{coding the A forcing into the ordinals}, from $\bbA \in \HOD$, one can obtain an enumeration, $\langle (S_p,\varphi_p) : p \in \bbA\rangle$ so that $(S_p,\varphi_p)$ is an $\OD$ $\infty$-Borel for the condition $p$.

First to show that $G_x$ is an $\bbA$-generic filter over $\HOD$: Let $A \subseteq \bbA$ be a maximal antichain which belongs to $\HOD$ and is hence $\OD$. Considering $\bbA$ as the set $K$ defined in Remark \ref{coding the A forcing into the ordinals} and using the fact that $A$ is $\OD$, one can find a formula $\varphi$ so that $(\bbA,\varphi)$ is an $\OD$-Borel code for $\bigcup A$. Therefore, $\bigcup A = \reals$ since otherwise $\reals \setminus \bigcup A$ would be a nonempty set with an $\OD$ $\infty$-Borel code. Then $\reals\setminus\bigcup A$ is a condition of $\bbA$ which is incompatible with every element of $A$. This contradicts $A$ being a maximal antichain. Thus $x \in \bigcup A$. There is some $a \in A$ so that $x \in a$. Thus $a \in G_x$. It has been shown that $G_x \cap A \neq \emptyset$. $G_x$ is $\bbA$-generic over $\HOD$.

Thinking of $\reals = \cantorspace$, let $b_n = \{x \in \reals : x(n) = 1\}$. Note that $b_n \neq \emptyset$ and $b_n$ clearly has an $\OD$ $\infty$-Borel code. Thus $b_n \in \bbA$. Note that $b_n \in G_x$ if and only if $x(n) = 1$. Thus $x \in \HOD[G_x]$.

Note that $p \in G_x$ if and only if $x \in p$ if and only if $L[S_p,x] \models \varphi_p(S_p,x)$. Note that $V \models L[S_p,x] \models \varphi_p(S_p,x)$ if and only if $\HOD[x] \models L[S_p,x]\models \varphi_p(S_p,x)$. (This is an application of the important absoluteness property of $\infty$-Borel codes.) Thus $G_x$ can be defined in $\HOD[x]$ as the set of $p$ such that $L[S_p,x] \models \varphi_p(S_p,x)$. This shows $G_x \in \HOD[x]$. 

It has been shown that $\HOD[x] = \HOD[G_x]$.
\end{proof}

\Begin{definition}{canonical name for a real}
Suppose $n \in \omega$. For each $s \in {}^n\omega$, let $b_s = \{x \in \reals^n : x(s) = 1\}$.

Let $\dot x_\mathrm{gen}^n = \{(\check s,b_{s}) : s \in {}^n\omega\}$. In light of the argument of Fact \ref{HOD bbO extension is same as extension by real}, if $x \in \reals^n$, then $\dot x_\mathrm{gen}^n[G_x] = x$.
\end{definition}

\Begin{definition}{forcing projection definition}
Let $\bbP$ and $\bbQ$ be two forcings. A surjective map $\pi : \bbQ \rightarrow \bbP$ is a forcing projection if and only if the following holds:

\noindent (1) For all $q_0,q_1 \in \bbQ$, $q_0 \leq_\bbQ q_1$ implies $\pi(q_0) \leq_\bbP \pi(q_1)$ and $\pi(1_\bbQ) = 1_\bbP$.

\noindent (2) For all $q \in \bbQ$ and $p' \in \bbP$ such that $p' \leq_\bbP \pi(q)$, there exists some $q' \leq_\bbQ q$ so that $\pi(q') = p'$. 

Let $\pi : \bbQ \rightarrow \bbP$ be a projection. If $D \subseteq \bbP$ is dense, then $\pi^{-1}[D]$ is dense in $\bbQ$. This implies that if $G \subseteq \bbQ$ is a $\bbQ$-generic over $V$, then $\pi[G]$ is a $\bbP$-generic over $V$.

Let $G \subseteq \bbP$ is $\bbP$-generic over $V$, then let $\bbQ \slash G = \{q \in \bbQ : \pi(q) \in G\}$. Let $\leq_{\bbQ \slash G} = \leq_\bbQ \upharpoonright \bbQ \slash G$.

Let $\dot G \in V^\bbP$ denote the $\bbP$-name for the generic filter. Let $\bbQ \slash \dot G$ denote the $\bbP$-name so that $1_\bbP \forces_\bbP \bbQ \slash \dot G = \{p \in \check \bbQ : \check\pi(p) \in \dot G\}$. 

It can be checked that $\bbQ$ embeds densely into the iteration $\bbP * (\bbQ \slash \dot G)$. Therefore, these two forcings are equivalent as forcings. Also if $G$ is $\bbP$-generic over $V$ then $(\bbQ \slash \dot G)[G] = \bbQ \slash G$. 

Suppose $H$ is $\bbQ$-generic over $V$. Let $G = \pi[H]$ be the associated $\bbP$-generic filter. One can check that $H$ is $(\bbQ \slash G)$-generic over $V[G]$, $G * H$ is $(\bbP * \bbQ \slash \dot G)$-generic over $V$ and $V[G][H] = V[G * H] = V[H]$. 
\end{definition}

\Begin{definition}{turing notions}
For the moment, consider $\reals = \cantorspace = \powerset{\omega}$. If $x,y \in \reals$, then one writes $x \leq_T y$ if and only if there is a Turing machine (taking oracle input) so that $x$ can be computed from this Turing machine when given $y$ as its oracle. Since Turing programs can be coded by natural numbers, for any $x \in \reals$, there are only countably many $y \in \reals$ so that $y \leq_T x$. A Turing program is also absolute between models of $\ZF$ with the same $\omega$.

Define $x =_T y$ if and only if $x \leq_T y$ and $y \leq_T x$.

Let $\degrees = \powerset{\omega} \slash =_T$ denote the collection of $=_T$ equivalence classes. An element of $\degrees$ is called a Turing degree. If $x \in \reals$, then one uses the notation $[x]_T$ rather than $[x]_{=_T}$ to denote the Turing degree of $x$. As observed above, each Turing degree contains only countable many reals.

If $X,Y \in \degrees$, one defines $X \leq Y$ if and only if there exist $x \in X$ and $y \in Y$ so that $x \leq_T y$. (One can check this is a well defined relation.)

A Turing cone of reals with base $x$ is the collection $C_x = \{y \in \reals : x \leq_T y\}$. A Turing cone of degrees with base $X$ is the collection $C_X = \{Y \in \degrees : X \leq Y\}$.

Let $\mu_\degrees \subseteq \powerset{\degrees}$ consists of those subsets of $\degrees$ which contain a Turing cone of degrees. $\mu_\degrees$ is a filter.
\end{definition}

\Begin{fact}{martin measure}
(Martin) Assume $\ZF + \AD$. $\mu_\degrees$ is a countably complete ultrafilter.
\end{fact}

\begin{proof}
Let $K \subseteq \degrees$. Let $\tilde K = \{x \in \reals : x \in K\}$. Note that $\tilde K$ is a $=_T$-invariant subset of $\reals$.

Consider the usual game $G_{\tilde K}$ 
$$\begin{tikzpicture}
\node at (0,0) {$G_{\tilde K}$};
\node at (1,.5) {I};
\node at (2,.5) {$a(0)$};
\node at (4,.5) {$a(1)$};
\node at (6,.5) {$a(2)$};
\node at (8,.5) {$a(3)$};
\node at (11,.5) {$a$};

\node at (1,-.5) {II};
\node at (3,-.5) {$b(0)$};
\node at (5,-.5) {$b(1)$};
\node at (7,-.5) {$b(2)$};
\node at (9,-.5) {$b(3)$};
\node at (11,-.5) {$b$};
\end{tikzpicture}$$
where Player 1 wins if and only if $a \oplus b \in \tilde K$, where $a \oplus b$ is defined by $a\oplus b(2n) = a(n)$ and $a \oplus b(2n + 1) = b(n)$.

By $\AD$, one of the two players has a winning strategy. Suppose Player 1 has a winning strategy $\sigma$.

Let $Z = [\sigma]_T$ be the Turing degree of $\sigma$. The claim is that $C_Z \subseteq K$:  

Let $Y \in \degrees$ be such that $Z \leq_T Y$. Pick any $y \in Y$. Thus $\sigma \leq_T y$. Let $\sigma * y$ denote the result of the play where Player 1 uses $\sigma$ and Player 2 plays the bits of $y$ each turn. Since $\sigma$ is a Player 1 winnings strategy, $\sigma * y \in \tilde K$.  Since $\sigma \leq_T y$, $\sigma * y \leq_T y$ and clearly $y \leq_T \sigma * y$. Thus $\sigma * y =_T y$. Since $\tilde K$ is $=_T$-invariant, $y \in \tilde K$. Thus $[y]_T = Y \in K$. 

This shows that if Player 1 has a winning strategy in $G_{\tilde K}$, then $K \in \mu_\degrees$. A similar argument shows that if Player 2 has a winning strategy then $\degrees\setminus K \in \mu_\degrees$. Thus $\mu_\degrees$ is an ultrafilter.

Next to show countable completeness: Suppose $\langle K_n : n \in \omega\rangle$ is a sequence in $\mu_\degrees$. By $\AC_\omega^\reals$, let $x_n$ be such that the cone above $X_n = [x_n]_T$ is contained in $K_n$. Let $x = \bigoplus x_n = \{(n,m) : m \in x_n\}$ and $X = [x]_T$. Then the cone above $X$ is contained inside of $\bigcap_{n \in \omega} K_n$. Thus $\bigcap_{n \in \omega} K_n \in \mu_{\degrees}$. 
\end{proof}

\Begin{remark}{turing degree and canonical wellorderings}
Let $S$ be a set. Note that if $x =_T y$, then $L[S,x] = L[S,y]$. Therefore, if $X \in \degrees$, one will often write $L[S,X]$ to denote $L[S,x]$ for any $x \in X$.

The canonical constructibility wellordering is based on the hierarchy $\{L_\alpha[S,x] : \alpha \in \ON\}$. Even if $x =_T y$, the levels $L_\alpha[S,x]$ and $L_\alpha[S,y]$ can differ. Thus, the canonical constructibility wellordering on $L[S,x]$ is not invariant under $=_T$.

If $V$ is a model of $\ZF$ (in the language $\dot \in$), then the canonical wellordering $<^{\HOD_S^V}$ of $\HOD^V_S$ depends only on $V$. Since $x =_T y$ implies that $L[S,x] = L[S,y]$, one has that $\HOD_S^{L[S,x]} = \HOD_S^{L[S,y]}$ and $<^{\HOD_S^{L[S,x]}} = <^{\HOD_S^{L[S,y]}}$. (Note that although the constructibility hierarchy of $L[S,x]$ is naturally formulated in the language $\{\dot \in, \dot E_0, \dot E_1\}$, where $\dot E_0$ and $\dot E_1$ are unary predicate symbols meant to interpret $S$ and $x$, when constructing $\HOD_S^{L[S,x]}$, $L[S,x]$ is considered as merely a $\{\dot \in\}$-structure.) Thus if $X$ is a Turing degree, one will often write $\HOD_S^{L[S,X]}$ and $\leq^{\HOD_S^{L[S,X]}}$ to refer to $\HOD^{L[S,x]}_S$ and $<^{\HOD_S^{L[S,x]}}$ for any $x \in X$.

The invariance of the canonical wellordering of $\HOD$ allows one to take ultrapowers of local $\HOD$'s (models of the form $\HOD_S^{L[S,X]}$) by the Martin's measure. This is a very powerful technique as indicated by the following results.
\end{remark}

The existence of the canonical wellorderings of the local $\HOD_S^{L[S,X]}$ shows that the resulting ultraproduct satisfies \Los's theorem:

\Begin{fact}{Los theorem for local hod}
Assume $\ZF + \AD$. Let $S \subseteq \ON$ be a set of ordinals. $\prod_{X \in \degrees}\HOD_S^{L[S,X]}\slash \mu_\degrees$ satisfies the \Los's theorem: Let $\varphi$ be a formula. Let $f_0,...,f_{n - 1}$ be functions on $\degrees$ with the property that for all $X \in \degrees$, $f_i(X) \in \HOD_S^{L[S,X]}$. Then 
$$\prod_{X \in \degrees}\HOD_S^{L[S,X]} \slash \mu_\degrees \models \varphi([f_0]_{\mu_\degrees},...,[f_{n - 1}]_{\mu_\degrees}) \Leftrightarrow \{X \in \degrees : \HOD_{S}^{L[S,X]} \models \varphi(f_0(X),...,f_{n - 1}(X))\} \in \mu_\degrees.$$
\end{fact}

\begin{proof}
Only the existential quantification case requires a choice-like principle. One will give a sketch:

Let $\CM$ denote this ultraproduct. Let $\varphi$ be a formula and assume inductively one has already shown the result for $\varphi$. Suppose
$$K = \{X \in \degrees : \HOD^{L[S,X]}_S \models (\exists v)\varphi(v,f_0(X),...,f_{n - 1}(X))\} \in \mu_\degrees.$$
Let $g$ be defined on $K$ by letting $g(X)$ be the $\HOD_S^{L[S,X]}$-least $v$ so that $\HOD_S^{L[S,X]}\models \varphi(v,f_0(X),...,f_{n - 1}(X))$. Then using the induction hypothesis, one can show that $\CM \models \varphi([g]_{\mu_\degrees},[f]_{\mu_\degrees}, ..., [f_{n - 1}]_{\mu_\degrees})$. Thus $\CM \models (\exists v)\varphi(v,[f]_{\mu_\degrees},...,[f_{n - 1}]_{\mu_\degrees})$.
\end{proof}

There is no claim that $\prod_{X \in \degrees}\HOD_{S}^{L[S,X]} \slash \mu_\degrees$ is a wellfounded model. This is true assuming $\DC$. It is open whether $\AD$ implies this.

\Begin{fact}{projection infinty borel set is infinty borel}
(Woodin) Assume $\ZF + \AD + \DC_\reals$. Suppose $A \subseteq \reals^2$ has an $\infty$-Borel code $(S,\varphi)$, then $B(x) = (\exists^\reals y)A(x,y)$ has an $\infty$-Borel code which is $\OD_S$. 
\end{fact}

\begin{proof}
Note that $L(S,\reals) \models \ZF + \AD + \DC$. Work in $L(S,\reals)$. 

Let $\CM = \prod_{X \in \degrees} \HOD_S^{L[S,X]} \slash \mu_\degrees$. $\CM$ is an $S$-definable class which is wellfounded by $\DC$. Implicitly, it will be assumed that all objects in this ultrapower have been Mostowski collapsed. By Fact \ref{Los theorem for local hod}, $\CM$ satisfies \Los's theorem. Define $\Phi_{\bbA_S^\infty}$ on $\degrees$ by $\Phi_{\bbA_S^\infty}(X) = \bbA_S^{L[S,X]}$. (Recall Remark \ref{coding the A forcing into the ordinals} concerning the convention on $\bbA_S$.) Let $\bbA^\infty_S = [\Phi_{A^\infty_S}]_{\mu_\degrees}$. By \Los's theorem, $\bbA_S^\infty$ is a forcing poset. Let $\lambda_X = |\bbA_S|^{\HOD_S^{L[S,X]}}$. Let $\Phi_\lambda$ be a function on $\degrees$ defined by $\Phi_\lambda(X) = \lambda_X$. Let $\lambda = [\Phi_{\lambda}]_{\mu_\degrees}$. By \Los's theorem, $\CM \models \lambda = |\bbA^\infty_S|$. Let $\Phi_{S^\infty}$ be a function of $\degrees$ defined by $\Phi_{S^\infty}(X) = S$. Let $S^\infty = [\Phi_{S^\infty}]_{\mu_\degrees}$.

Claim: For all $a \in \reals$, 
$$a \in B \Leftrightarrow L[S^\infty,\bbA_S^\infty,a] \models 1_{\mathrm{Coll}(\omega,\lambda)} \forces (\exists b)L[S^\infty,a,b] \models \varphi(S^\infty,a,b).$$

To prove the claim: First observe that for all $X \in \degrees$, $\lambda_X$ is countable in $L(S,\reals)$ which is a model of $\AD$. To see this: Note that $\bbR^{L[S,X]}$ is countable since it is a wellorderable collection of reals; thus, there is a bijection of $\omega$ with $\reals^{L[S,X]}$. $\bbA_S^{L[S,X]}$ is a collection of subsets of $\bbR^{L[S,X]}$ in $L[S,X]$. Identifying $\reals^{L[S,X]}$ with $\omega$, this collection $\bbA_S^{L[S,X]}$ can be identified as a wellorderable collection of $\reals$, as well. Thus $\lambda_X = |\bbA_S|^{L[S,X]}$ is countable in $L(S,\reals)$. By the same reasoning, $(2^{\lambda_X})^{L[S,X]}$ is countable in $L(S,\reals)$.

$(\Leftarrow)$ For all $X \in \degrees$ so that $a \in X$, one can define $\HOD_S^{L[S,X]}[a]$ as in Fact \ref{HOD bbO extension is same as extension by real}. Let 
$$\CM[a] = \prod_{X \in \degrees} \HOD^{L[S,X]}_S[a] \slash \mu_\degrees.$$
Assume that 
$$V \models L[S^\infty,\bbA_S^\infty,a] \models 1_{\mathrm{Coll}(\omega,\lambda)} \forces (\exists b)L[S^\infty,a,b] \models \varphi(S^\infty,a,b).$$
Thus
$$\CM[a] \models L[S^\infty,\bbA_S^\infty,a] \models 1_{\mathrm{Coll}(\omega,\lambda)} \forces (\exists b)L[S^\infty,a,b] \models \varphi(S^\infty,a,b).$$
By the idea of the proof of \Los's theorem (Fact \ref{Los theorem for local hod}), for $\mu_\degrees$-almost all $X \in \degrees$ with $a \in X$, 
$$L[S,\bbA_S^{L[S,X]},a] \models 1_{\mathrm{Coll}(\omega,\lambda_X)}\forces (\exists b)L[S,a,b]\models\varphi(S,a,b).$$
Fix such an $X$. Since $(2^{\lambda_X})^{L[S,X]}$ is countable in $L(S,\reals)$, there is a $g \in L(S,\reals)$ so that $g \subseteq \mathrm{Coll}(\omega,\lambda_X)$ is $\mathrm{Coll}(\omega,\lambda_X)$-generic over $\HOD_{S}^{L[S,X]}[a]$. Then $g$ is also generic over $L[S,\bbA_S^{L[S,X]},a]$. By the forcing theorem, $L[S,\bbA_S^{L[S,X]},a][g] \models (\exists b)L[S,a,b] \models \varphi(S,a,b)$. Pick some $b \in L[S,\bbA_S^{L[S,X]},a][g]$ witnessing the existential. Then one has in particular that $L[S,a,b] \models \varphi(S,a,b)$. Since $(S,\varphi)$ is the $\infty$-Borel code for $A$, one has that $(a,b) \in A$. Thus $a \in B$. 

$(\Rightarrow)$ Suppose that $a \in B$. Let $b \in \reals$ be such that $(a,b) \in A$. Let $X$ be a Turing degree such that $[a \oplus b] \leq X$. By Fact \ref{HOD bbO extension is same as extension by real}, let $G_{a \oplus b}$ be the $\bbA_{S}^{L[S,X]}$-generic over $\HOD_S^{L[S,X]}$-filter derived from $a \oplus b$. Note that $G_{a\oplus b}$ is also $\bbA_S^{L[S,X]}$-generic over $L[S,\bbA_S^{L[S,X]}]$. Using the convention from Remark \ref{coding the A forcing into the ordinals} and an argument similar to Fact \ref{HOD bbO extension is same as extension by real}, one can recover $a \oplus b$ from $G_{a,b}$ in $L[S,\bbA_S^{L[S,X]}][G_{a,b}]$. Thus
$$L[S,\bbA_S^{L[S,X]}][G_{a\oplus b}] \models (\exists y)L[S,a,y]\models\varphi(S,a,y).$$
Since $|\bbA_S^{L[S,X]}|$ is a forcing of size $\lambda_X$, one has that $\bbA_S^{L[S,X]}$ regularly embeds into $\mathrm{Coll}(\omega,\lambda_X)$ by \cite{Set-Theory} Corollary 26.8. There is some $g \subseteq \mathrm{Coll}(\omega,\lambda_X)$ which is generic over $L[S,\bbA_S^{L[S,X]}]$ so that $G_{a\oplus b} \in L[S,\bbA_S^{L[S,X]}][g]$. Thus
$$L[S,\bbA_S^{L[S,X]}][g] \models (\exists y)L[S,a,y]\models\varphi(S,a,y).$$
Note that $a \in L[S,\bbA_S^{L[S,X]}][g]$. By one of the main properties of the collapse forcing (\cite{Set-Theory} Corollary 26.10), there is some $h \subseteq \mathrm{Coll}(\omega,\lambda_X)$ which is $\mathrm{Coll}(\omega,\lambda_X)$-generic over $L[S,\bbA_S^{L[S,X]}][a]$ so that $L[S,\bbA_S^{L[S,X]}][a][h] = L[S,\bbA_S^{L[S,X]}][g]$. Thus
$$L[S,\bbA_S^{L[S,X]}][a][h] \models (\exists y)L[S,a,y]\models\varphi(S,a,y).$$
There is some $p \in \mathrm{Coll}(\omega,\lambda_X)$ which forces the inner statement. By the homogeneity of $\mathrm{Coll}(\omega,\lambda_X)$, $1_{\mathrm{Coll}(\omega,\lambda_X)}$ forces this statement:
$$L[S,\bbA_S^{L[S,X]}][a] \models 1_{\mathrm{Coll}(\omega,\lambda_X)} \forces (\exists y)L[S,a,y]\models\varphi(S,a,y).$$
In particular, 
$$\HOD_{S}^{L[S,X]}[a] \models L[S,\bbA_S^{L[S,X]},a] \models 1_{\mathrm{Coll}(\omega,\lambda_X)} \forces (\exists y)L[S,a,y]\models\varphi(S,a,y).$$
By \Los's theorem (Fact \ref{Los theorem for local hod}), one has
$$\CM[a] \models L[S^\infty,\bbA_S^\infty,a] \models 1_{\mathrm{Coll}(\omega,\lambda)} \forces (\exists y)L[S^\infty,a,y]\models\varphi(S^\infty,a,y).$$
So in particular, 
$$L[S^\infty,\bbA_S^\infty,a] \models 1_{\mathrm{Coll}(\omega,\lambda)} \forces (\exists y)L[S^\infty,a,y]\models\varphi(S^\infty,a,y).$$
This completes the proof of the claim.

Let $J$ be a set of ordinals that codes $S^\infty$ and $\bbA_S^\infty$ in some fixed way. Let $\psi(J,a)$ be the formula that asserts 
$$1_{\mathrm{Coll}(\omega,\lambda)} \forces (\exists y)L[S^\infty,a,y] \models \varphi(S^\infty,a,y).$$
$(J,\psi)$ is an $\infty$-Borel code for $B$. 
\end{proof}

\Begin{fact}{projection of infty borel forcing}
(Woodin) Assume $\ZF + \AD + \DC$. Let $1 \leq m \leq n$. Let $\pi_{n,m} : \reals^n \rightarrow \reals^m$ be defined by $\pi_{n,m}(s) = s\upharpoonright m$. Define $\pi_{n,m} : {}_n\bbA \rightarrow {}_m\bbA$ by $\pi_{n,m}(b) = \pi_{n,m}[b]$, where the latter $\pi_{n,m}$ refers to the earlier function $\pi_{n,m} : \reals^n \rightarrow \reals^m$. Then $\pi_{n,m}$ is a forcing projection.
\end{fact}

\begin{proof}
By Fact \ref{projection infinty borel set is infinty borel}, one can show that $\pi_{n,m}(p) \in {}_m\bbA$ for all $p \in {}_n\bbA$. 

Now suppose that $q \in {}_n\bbA$. Suppose $p' \in {}_m\bbA$ is such that $p' \leq_{{}_m\bbA} \pi_{n,m}(q)$. Let $q' = \{x \in \reals^n : x\upharpoonright n \in \pi_{n,m}(q)\}$. Since $\pi_{n,m}(q) \in {}_m\bbA$ (that is, it has an $\OD$ $\infty$-Borel code), $q'$ has an $\OD$ $\infty$-Borel code and hence belongs to ${}_n\bbA$. Clearly, $\pi_{n,m}(q') = p'$. This establishes that $\pi_{n,m}$ is a forcing projection.
\end{proof}

Now one will work in $L(\reals) \models \AD$. Kechris (\cite{The-Axiom-of-Determinacy-Implies-Dependent-Choice}) showed that if $\AD$ holds, then $L(\reals) \models \DC_\reals$. Thus, for the following results, the background theory $\ZF + \AD$ is sufficient.

\Begin{definition}{more names in the symmetric extension}
Using the projections from Fact \ref{projection of infty borel forcing}, let ${}_\omega\bbA$ be the finite support direct limit of $\langle {}_n\bbA : n \in \omega \setminus \{0\}\rangle$. That is, ${}_{\omega}\bbA$ is the collection of $p : (\omega \setminus \{0\}) \rightarrow \bigcup_{n \in \omega}{}_n\bbA$ so that for all $m \leq n$, $\pi_{n,m}(p(n)) = p(m)$ and there exists a $N \in \omega$ so that for all $k \geq N$, $p(k) = p(N) \times \reals^{k - N}$. The least such $N$ is denoted $\dim(p)$, the dimension of $p$. For $n \in \omega$ and $p \in {}_\omega\bbA$, let $\pi_{\omega,n}(p) = p(n)$. Each $\pi_{\omega,n} : {}_\omega\bbA \rightarrow {}_n\bbA$ is a forcing projection.

Since each ${}_n\bbA \in \HOD$ is identified as a set of ordinals (see Remark \ref{coding the A forcing into the ordinals}) and the projection maps are in $\HOD$, ${}_\omega\bbA$ belongs to $\HOD$ and may be identified as a set of ordinals having the property expressed in Remark \ref{coding the A forcing into the ordinals}.

Let $m \leq n$, let $\tau_m$ be a ${}_n\bbA$-name for the last real of the $\dot x_\mathrm{gen}^m$, which is a name for the generic $m$-tuple coming from a ${}_n\bbA$-generic filter. (Technically, $\tau_m$ is different for each $n$, but the projection can be used to interpret it in suitable $n$'s.)

Let $\symreals$ be a ${}_\omega\bbA$-name so that $1_{{}_\omega\bbA}\forces \symreals = \{\tau_n : n \in \omega\setminus\{0\}\}$. 
\end{definition}

Observe that from Fact \ref{HOD bbO extension is same as extension by real}, every $z \in \reals^n$ induces a ${}_n\bbA$-generic filter $G_z$ over $\HOD$ so that $\HOD[z] = \HOD[G_z]$. Note that $\dot x_\mathrm{gen}[G_z] = z$ and $\tau_m[G_z] = z(m)$ for all $m < n$. 

\Begin{theorem}{coll generic over L(R) gives a omegabbA generic}
(Woodin) Assume $\ZF + \AD + \mathsf{V = L(\reals)}$. Suppose $g \subseteq \mathrm{Coll}(\omega,\reals)$-generic over $L(\reals)$. From $g$, one can derive a ${}_\omega\bbA$-generic over $\HOD^{L(\reals)}$ filter $G_g$.
\end{theorem}

\begin{proof}
Let $g \subseteq \mathrm{Coll}(\omega,\reals)$ be a generic over $\HOD^{L(\reals)}$. Let $G_g \subseteq {}_\omega\bbA$ be the collection of condition $p \in {}_\omega\bbA$ so that $g \upharpoonright \dim(p) \in p(\dim(p))$. 

The claim is that $G_g$ is ${}_\omega\bbA$-generic over $\HOD^{L(\reals)}$. 

Suppose $D \subseteq {}_\omega\bbA$ belongs to $\HOD^{L(\reals)}$ and is dense. Let $\tilde D \subseteq \mathrm{Coll}(\omega,\reals)$ be the collection of $s \in \mathrm{Coll}(\omega,\reals)$ so that there is some $p \in {}_\omega\bbA$ with $\dim(p) = |s|$, $s \in p(|s|)$, and $p \in D$. (Note that $s$ as a condition of $\mathrm{Coll}(\omega,\reals)$ is a finite tuple of reals.)

One will show that $\tilde D$ is dense in $\mathrm{Coll}(\omega,\reals)$. Let $s \in \mathrm{Coll}(\omega,\reals)$. Let $n = |s|$. Let 
$$E = \{p \in {}_n\bbA : (\exists q \in {}_\omega\bbA)(\dim(q) \geq n \wedge q \in D \wedge \pi_{\omega,n}(q) = p)\}.$$
First, one will show $E$ is dense in ${}_n\bbA$. Let $r \in {}_n\bbA$. Since $D$ is dense in ${}_\omega\bbA$, there is some $q \in {}_\omega\bbA$ with $\dim(q) \geq n$ so that $q \leq_{{}_\omega\bbA} r$ and $q \in D$. Then $p = \pi_{\omega,n}(q)$ belongs to $E$ and $p \leq_{{}_n\bbA} r$. Thus $E$ is dense in ${}_n\bbA$. 

Let $G_s^n$ be the ${}_n\bbA$-generic over $\HOD^{L(\reals)}$ filter derived from $s$. By genericity, pick some $p \in G_s^n \cap E$. Thus there is some $q \in D$ so that $\pi_{\omega,n}(q) = p$. Let $s' \supseteq s$ be such that $s' \in q$. This means that $s' \leq_{\mathrm{Coll}(\omega,\reals)} s$ and $s' \in \tilde D$. It has been shown that $\tilde D$ is dense in $\mathrm{Coll}(\omega,\reals)$. 

Now since $g \subseteq \mathrm{Coll}(\omega,\reals)$ is $\mathrm{Coll}(\omega,\reals)$-generic over $L(\reals)$, $g \cap \tilde D \neq \emptyset$. There is some $n \in \omega$ so that $g \upharpoonright n \in \tilde D$. By definition, there is some $p \in D$ so that $g\upharpoonright n \in p(n)$. Thus $p \in G_g \cap D$. 

This shows that $G_g$ is ${}_\omega\bbA$-generic over $\HOD$. 
\end{proof}

\Begin{fact}{symreals are the reals}
(Woodin) Assume $\ZF + \AD + \mathsf{V = L(\reals)}$. $\HOD^{L(\reals)} \models 1_{{}_\omega\bbA} \forces$ ``the reals of $L(\symreals)$ is $\symreals$''.
\end{fact}

\begin{proof}
Let $p \in {}_\omega\bbA$ be some condition and let $n = \dim(p)$. Let $s \in \reals^n$ with $s \in p(n)$. Consider $s$ as a condition of $\mathrm{Coll}(\omega,\reals)$. Let $g \subseteq \mathrm{Coll}(\omega)$ be $\mathrm{Coll}(\omega,\reals)$-generic over $L(\reals)$ containing $s$. An easy density argument shows that if $g$ is considered as a function from $\omega$ to $\reals$, $g$ must be a surjection onto $\reals^{L(\reals)}$. Let $G_g$ be the ${}_\omega\bbA$-generic filter over $\HOD^{L(\reals)}$ derived from $g$. Then $\HOD^{L(\reals)}[G_g] \models$ ``the reals of $L(\symreals[G_g])$ is $\symreals[G_g]$'' since $L(\symreals[G_g]) = L(\reals)$. Thus there is some $q \leq_{{}_\omega\bbA} p$ so that $\HOD^{L(\reals)} \models q \forces_{{}_\omega\bbA}$ ``the reals of $L(\symreals)$ is $\symreals$''. Since $p$ was an arbitrary condition, one has that $1_{{}_\omega\bbA}$ forces this same statement. 
\end{proof}

\Begin{theorem}{all fact L(symreal) forced in the same way}
(Woodin) Assume $\ZF + \AD + \mathsf{V + L(\reals)}$. Let $s \in \reals^n$. Suppose $G^n_s$ is the ${}_n\bbA$-generic filter over $\HOD^{L(\reals)}$ derived from $s$. Let $\varphi$ be a formula. Suppose $z \in \HOD^{L(\reals)}[G^n_s]$. Then 
$$\HOD^{L(\reals)}[G_s^n] \models 1_{{}_\omega\bbA \slash G_s^n} \forces_{{}_\omega\bbA \slash G^n_s} L(\symreals) \models \varphi(\check z,\dot x_\mathrm{gen}^n)$$
or
$$\HOD^{L(\reals)}[G_s^n] \models 1_{{}_\omega\bbA \slash G_s^n} \forces_{{}_\omega\bbA \slash G_s^n} L(\symreals) \models \neg\varphi(\check z,\dot x_\mathrm{gen}^n)$$
\end{theorem}

\begin{proof}
In $L(\reals)$, one either has $L(\reals) \models \varphi(z,s)$ or $L(\reals) \models \neg\varphi(z,s)$. The claim is that which ever case occurs, this is the side that is forced over $\HOD^{L(\reals)}[G_s^n]$. 

So without loss generality, suppose $L(\reals) \models \varphi(z,s)$. 

Take any $q \in {}_\omega\bbA \slash G_s^n$ with $\dim(q) \geq n$. Let $m = \dim(q)$. This means that $\pi_{\omega,n}(q) \in G_s^n$. Let $r \in q$ be a $\reals^m$ sequence extending $s$. Let $g \subseteq \mathrm{Coll}(\omega,\reals)$ be a $\mathrm{Coll}(\omega,\reals)$-generic over $L(\reals)$ filter so that $r \subseteq g$. Let $G_g$ be the ${}_\omega\bbA$-generic over $\HOD^{L(\reals)}$ filter derived from $g$. Note that $G_g$ is a ${}_\omega\bbA \slash G_s^n$-generic filter over $\HOD^{L(\reals)}[G_s^n]$ (see the basic properties of projection from Definition \ref{forcing projection definition}) and that $q \in G_g$. Using Fact \ref{symreals are the reals} and the fact that $\HOD[G^n_s][G_g] = \HOD[G_g]$ (also see Definition \ref{forcing projection definition}),
$$\HOD[G_s^n][G_g] \models L(\symreals[G_g])\models \varphi(z, \dot x_\mathrm{gen}^n[G_s^n])$$
since $L(\reals) \models \varphi(z,s)$. Hence there is some $q' \leq_{{}_\omega\bbA \slash G_s^n} q$ so that
$$\HOD[G_s^n] \models q' \forces_{{}_\omega\bbA\slash G_s^n} L(\symreals) \models \varphi(\check z,\dot x_\mathrm{gen}^n).$$
Since $q \in {}_\omega\bbA\slash G_s^n$ was arbitrary, $1_{{}_\omega\bbA\slash G_s^n}$ forces this statement. This completes the proof.
\end{proof}

\Begin{theorem}{ultimate infinity borel code}
(Woodin) Assume $\ZF + \AD + \mathsf{V = L(\reals)}$. Let $s \in \reals^n$, $z \in L[{}_\omega\bbA,s]$, and $\varphi$ be a formula. Then $L(\reals) \models \varphi(s,z)$ if and only if 
$$L[{}_\omega\bbA,s] \models 1_{{}_\omega\bbA \slash G_s^n} \forces_{{}_\omega\bbA \slash G_s^n} L(\symreals) \models \varphi(\dot x_\mathrm{gen}^n, \check z)$$
\end{theorem}

\begin{proof}
Recall ${}_\omega\bbA \in \HOD^{L(\reals)}$. Using the convention from Remark \ref{coding the A forcing into the ordinals}, one may assume that from a $p \in {}_\omega\bbA$, one can obtain an $\infty$-Borel code for this condition. From this, one can see that given $s$ and ${}_\omega\bbA$, one can reconstruct $G_s^n$. Also the names $\dot x_\mathrm{gen}^n$ and $\symreals$ belong to $L[{}_\omega\bbA]$.

$(\Rightarrow)$ Using the idea from Fact \ref{HOD bbO extension is same as extension by real}, one can show that $L[{}_\omega\bbA,s] = L[{}_\omega\bbA][G_s^n]$. Take any $p \in {}_\omega\bbA\slash G^n_s$. In particular, $p \in {}_\omega\bbA$. Pick some $r \in \mathrm{Coll}(\omega,\reals)$ with $r \supseteq s$ and $r \in p(|r|)$. Let $g \subseteq \mathrm{Coll}(\omega,\reals)$ be $\mathrm{Coll}(\omega,\reals)$-generic over $L(\reals)$ and $r \subseteq g$. Using ${}_\omega\bbA$, one can reconstruct $G_g$ which is a ${}_\omega\bbA$-generic filter over $\HOD^{L(\reals)}$. Thus $G_g$ is $({}_\omega\bbA \slash G^n_s)$-generic over $\HOD^{L(\reals)}[G^n_s]$ and hence generic over $L[{}_\omega\bbA,s]$. Since $L(\reals) \models \varphi(s,z)$, $L[{}_\omega\bbA,s][G_g] \models L(\symreals[G_g]) \models \varphi(s,\check z)$ because $L(\symreals[G_g]) = L(\reals)$. Pick a $q \leq_{{}_\omega\bbA} p$ so that 
$$L[{}_\omega\bbA,s] \models q \forces_{{}_\omega\bbA \slash G_s^n} L(\symreals)\models \varphi(s,\check z).$$
Since $p$ was arbitrary, $1_{{}_\omega\bbA\slash G_s^n}$ forces this statement. 

$(\Leftarrow)$ Let $g \supseteq s$ be a $\mathrm{Coll}(\omega,\reals)$-generic filter over $L(\reals)$. Let $G_g$ be the ${}_\omega\bbA$-generic filter of $\HOD^{L(\reals)}$ derived from $g$. It is $({}_\omega\bbA \slash G_s^n)$-generic over $\HOD^{L(\reals)}[G^n_s]$ and hence generic over $L[{}_\omega\bbA,s]$. Then one has
$$L[{}_\omega\bbA,s][G_g] \models L(\symreals[G_g]) \models \varphi(s,z).$$
However, $L(\symreals[G_g]) = L(\reals)$. Thus $L(\reals) \models \varphi(s,z)$.
\end{proof}

\Begin{corollary}{infity borel codes in L(R)}
Assume $\ZF + \AD + \mathsf{V = L(\reals)}$. Every set of reals has an $\infty$-Borel code. Moreover, if $A$ is $\OD_{s}$ where $s$ is a finite tuple of reals, then $A$ has an $\infty$-Borel code in $L[{}_{\omega}\bbA,s] \subseteq \HOD_s$. 
\end{corollary}

\begin{proof}
Suppose $A \subseteq \reals$. In $L(\reals)$, every set is definable from ordinals and some reals. Let $s \in \reals^n$, $\bar{\alpha}$ be a tuple of ordinals, and $\varphi$ be a formula so that $r \in A \Leftrightarrow \varphi(s,r,\bar{\alpha})$. By Theorem \ref{ultimate infinity borel code}, one has that $r \in A$ if and only if
$$L(\reals) \models L[{}_\omega\bbA,s,r] \models 1_{{}_\omega\bbA\slash G_{s\hat{\ }r}^{n + 1}} \forces L(\symreals) \models \varphi(s,r,\bar{\alpha}).$$
Let $S$ be a set of ordinals coding $({}_\omega\bbA,s,\bar{\alpha})$. Let $\psi$ be the formula expressing the inner forcing statement. One has that $(S,\psi)$ is an $\infty$-Borel code for $A$. 
\end{proof}

Note that ${}_\omega\bbA$ can be considered an ultimate $\infty$-Borel code in the sense that for any $A \subseteq \reals$, there is some tuple of reals, tuple of ordinals, and a formula so that together with ${}_\omega\bbA$ they form an $\infty$-Borel code for $A$. This is very useful in diagonalization arguments.

\Begin{corollary}{HOD and L direct limit vopenka}
Assume $\ZF + \AD + \mathsf{V = L(\reals)}$. Then $\HOD^{L(\reals)} = L[{}_\omega\bbA]$. 
\end{corollary}

\begin{proof}
Since ${}_\omega\bbA \in \HOD^{L(\reals)}$, $L[{}_\omega\bbA]\subseteq \HOD^{L(\reals)}$. $\HOD^{L(\reals)} \subseteq L[{}_\omega\bbA]$ by Theorem \ref{ultimate infinity borel code}.
\end{proof}

\Begin{corollary}{vopenka and infty borel forcing are the same}
Assume $\ZF + \AD + \mathsf{V = L(\reals)}$. $\bbA = \bbO$, that is the forcing with nonempty subsets of $\reals$ possessing $\OD$ $\infty$-Borel codes is the same as the Vop\'{e}nka forcing.
\end{corollary}

\begin{proof}
This follows from Corollary \ref{infity borel codes in L(R)}.
\end{proof}

Therefore, in the following, one will use the more practical Vop\'enka forcing $\bbO$ rather than $\bbA$.

\section{The Vop\'enka Forcing}\label{vopenka forcing}

\Begin{fact}{restate vopenka theorem}
(\cite{Dichotomy-for-Definable-Universe} Theorem 2.4.) Let $M$ be a transitive inner model of $\ZF$. Let $S \in M$ be a set of ordinals. Suppose $K$ is an $\OD^{M}_S$ set of ordinals and $\varphi$ is a formula. Let $N$ be some transitive inner model with $\HOD_S^M \subseteq N$. Suppose $p = \{x \in \reals^M : L[K,x] \models \varphi(K,x)\}$ is a condition of $\bbO^M_S$ (i.e. is nonempty). Then $N \models p \forces_{\bbO_S^M} L[\check K,\dot x_\mathrm{gen}] \models \varphi(\check K, \dot x_\mathrm{gen})$. 
\end{fact}

\begin{proof}
Suppose not. Then there is some $q' \leq_{\bbO^M_S} p$ so that $N \models q' \forces_{\bbO^M_S} L[\check K,\dot x_\mathrm{gen}] \models \neg\varphi(\check K,\dot x_\mathrm{gen})$. Because every $\bbO_S^M$-generic filter over $N$ is also generic over $\HOD_S^M$, one can find some $q \leq_{\bbO_S^M} q'$ so that $\HOD_S^M \models q \forces_{\bbO_S^M} L[\check K,\dot x_\mathrm{gen}] \models \neg\varphi(\check K,\dot x_\mathrm{gen})$. Since $q \neq \emptyset$, let $y \in q$. Let $G_y$ be the $\bbO_S^M$-generic filter over $\HOD_S^M$ derived from $y$. Note that $q \in G_y$. Thus $\HOD_S^M[G_y] \models L[K,y] \models \neg \varphi(K,y)$. Hence $L[K,y] \models \neg\varphi(K,y)$. This contradicts $q \subseteq p$. 
\end{proof}

\Begin{fact}{individual vopenka generic}
Let $M$ be an inner model of $\ZF$. Let $S \in M$ be a set of ordinals. Let $N$ be an inner model of $\ZF$ such that $N \supseteq \HOD_S^M$. Let $n \geq 1$ be a natural number. Suppose $(g_0,...,g_{n - 1})$ is an ${}_n\bbO_S^M$-generic reals over $N$. Then each $g_0$,...,$g_{n - 1}$ is $\bbO^M_S$-generic over $N$.
\end{fact}

\begin{proof}
Here $g$ is a $\bbO_S^M$-generic real over $N$ if and only if there is a filter $G$ which $\bbO_S^M$-generic over $N$ so that $g = \dot x_\mathrm{gen}[G]$. For simplicity assume $n = 2$. 

Let $\pi_1 : \reals^2 \rightarrow \reals$ be the projection onto the first coordinate. For each $p \in {}_2\bbO_S^M$, $\pi_1[p] \in \bbO_S^M$. 

Let $(g_0,g_1)$ be ${}_2\bbO_S^M$-generic over $N$. Let $G_{(g_0,g_1)}$ be the ${}_2\bbO^M_S$-generic filter over $N$ which adds $(g_0,g_1)$. Let $G = \{\pi_1(p) : p \in G_{(g_0,g_1)}\}$. $G$ is a filter on $\bbO_S^M$. 

Suppose $D \subseteq \bbO_S^M$ belongs to $N$ and is dense open. Let $D' = \{p \in {}_2\bbO_S^M : \pi_1(p) \in D\}$. Let $r \in {}_2\bbO_S^M$. Since $D$ is dense, there is some $r' \leq_{\bbO_S^M} \pi_1(r)$ with $r' \in D$. Let $s = (r' \times \reals) \cap r$. Note $s \in {}_2\bbO_S^M$, $\pi_1(s) = r' \in D$, and $s \leq_{{}_2\bbO_S^M} r$. Hence $s \in D'$. This shows that $D'$ is dense in ${}_2\bbO_S^M$. 

By genericity, there is some $r \in D' \cap G_{(g_0,g_1)}$. Then $\pi_1(r) \in D \cap G$. This shows that $G$ is $\bbO^M_S$-generic over $N$. Since $g_0$ is the real added by $G$, $g_0$ is a $\bbO_S^M$-generic real over $N$.
\end{proof}

First, one will give a proof of Woodin's countable section uniformization.  This follows a presentation of this result in \cite{Ramsey-Ultrafiler-and-Countable-to-One-Uniformation}. It is not known whether countable section uniformization is provable in $\AD$ alone.

\Begin{lemma}{selection lemma for HOD}
Assume $\ZF + \AD$. Let $R \subseteq \reals \times \reals$ be a relation. Suppose there is a set of ordinals $S$ with the property that for all $x \in \reals$ so that $R_x = \{y \in \reals : R(x,y)\} \neq \emptyset$ (i.e. $x \in \dom(R)$), for $\mu_\degrees$-almost all Turing degrees $Z$, $R_x \cap \HOD_{S,x}^{L[S,x,Z]} \neq \emptyset$. Then $R$ has a uniformization.
\end{lemma}

\begin{proof}
Think of $\reals = \cantorspace$. Let $x \in \dom(R)$. Let $E^x = \{Z \in \degrees: R_x \cap \HOD_{S,x}^{L[S,x,Z]} \neq \emptyset\}$. By the assumption, $E^x \in \mu_\degrees$. 

For each $n \in \omega$ and $i \in 2$, let $E^x_{n,i}$ be the collection of $Z \in E^x$ so that $z(n) = i$, where $z$ is the $\HOD_{S,x}^{L[S,x,Z]}$-least element of $R_x \cap \HOD_{S,x}^{L[S,x,Z]}$ according to the canonical wellordering of $\HOD_{S,x}^{L[S,x,Z]}$. Note that $E^x = E^x_{n,0} \cup E^x_{n,1}$. Since $\mu_\degrees$ is an ultrafilter and $E^x \in \mu_\degrees$, for each $n \in \omega$, there is a unique $a^x_n \in 2$ so that $E^x_{n,a_n^x} \in \mu_\degrees$. 

Define a function $\Phi : \dom(R) \rightarrow \cantorspace$ by $\Phi(x)(n) = a_n^x$. The claim is that $\Phi$ is a uniformization for $R$.

Suppose $x \in \dom(R)$. Since $\mu_\degrees$ is a countably complete ultrafilter on $\degrees$, $\bigcap_{n \in \omega} E^x_{n,a^x_n} \in \mu_\degrees$. Pick any $Z \in \bigcap_{n \in \omega} E^x_{n,a^x_n}$. Then $\Phi(x)$ is the $\HOD_{S,x}^{L[S,x,Z]}$-least element of $R_x \cap \HOD_{S,x}^{L[S,x,Z]}$. In particular, $\Phi(x) \in R_x$. $\Phi$ is a uniformization.
\end{proof}

\Begin{theorem}{woodin countable section uniformization}
(Woodin's countable section uniformization) Assume $\ZF + \AD$ and that all sets of reals have an $\infty$-Borel code. Let $R \subseteq \reals \times \reals$ be such that for all $x \in \dom(R)$, $R_x = \{y \in \reals : R(x,y)\}$ is countable. Then $R$ has a uniformization function, that is some function $F : \dom(R) \rightarrow \reals$ so that for all $x \in \dom(R)$, $R(x,F(x))$.  

In particular, countable section uniformization holds in $L(\reals) \models \AD$.
\end{theorem}

\begin{proof}
Let $(S,\varphi)$ be an $\infty$-Borel code for $R$, that is $\mathfrak{B}^2_{(S,\varphi)} = R$.

Let $x \in \dom(R)$. Thus there is some $y^*$ such that $R(x,y^*)$. Let $Z \in \degrees$ be such that $Z \geq [y^*]_T$. In the model $L[S,x,Z]$, define 
$$p = \{y \in \reals : L[S,x,y] \models \varphi(S,x,y)\}.$$
This condition is clearly $\OD_{S,x}^{L[S,x,Z]}$. Note $[y^*]_T \leq Z$ implies that $y^* \in L[S,x,Z]$. Since in the real world $V$, $R(x,y^*)$ holds, one has that $V \models L[S,x,y^*] \models \varphi(S,x,y^*)$. Thus $L[S,x,Z] \models L[S,x,y^*] \models \varphi(S,x,y^*)$. So $y^* \in p$. Hence $p \neq \emptyset$. It has been shown that $p \in \bbO_{S,x}^{L[S,x,Z]}$.

Claim 1: There is a dense below $p$ set of conditions in $\bbO_{S,x}^{L[S,x,Z]}$ which forces a value for $\dot x_\mathrm{gen}$.

To prove Claim 1: Suppose otherwise. As argued in the proof of Fact \ref{projection infinty borel set is infinty borel}, since $\HOD_{S,x}^{L[S,x,Z]}\models \AC$, $\mathscr{P}(\bbO_{S,x}^{L[S,x,Z]})^{\HOD_{S,x}^{L[S,x,Z]}}$ is countable in the real world $V$. Let $\langle D_n : n \in \omega\rangle$ enumerate (in the real world) all the dense open subsets of $\bbO_{S,x}^{L[S,x,Z]}$ that belong to $\HOD_{S,x}^{L[S,x,Z]}$.

Let $p_\emptyset$ be any condition below $p$ that meets $D_0$. Let $m_\emptyset = 0$. Suppose for some $\sigma \in \finBinarySequence$, $p_\sigma$ and $m_\sigma$ have been defined. First find some $p' \leq p_\sigma$ so that $p' \in D_{|\sigma| + 1}$. Since $p'$ can not determine $\dot x_\mathrm{gen}$, there is some $N \geq m_\sigma$ and two conditions $q_0$ and $q_1$ below $p'$ so that 

\noindent (1) $q_0$ and $q_1$ determine $\dot x_\mathrm{gen} \upharpoonright N + 1$ and

\noindent (2) For $i \in 2$, $q_i \forces \dot x_\mathrm{gen}(\check N) = \check i$. (That is, $q_i$ forces the generic real to take value $i$ at $N$.)

Let $q_0$ and $q_1$ be the $\HOD_{S,x}^{L[S,x,Z]}$-least pair with the above property. Now let $p_{\sigma\hat{\ }i} = q_i$ and $m_{\sigma\hat{\ }i} = N + 1$. Observe that $p_{\sigma\hat{\ }i} \in D_{n + 1}$ (as $D_{n + 1}$ is dense open). This completes the construction of a sequence $\langle p_{\sigma} : \sigma \in \finBinarySequence\rangle$ of conditions in $\bbO^{L[S,x,Z]}_{S,x}$.

For each $f \in \cantorspace$, let $G_f$ be the $\leq_{\bbO_{S,x}^{L[S,x,Z]}}$ upward closure of the set $\{p_{f\upharpoonright n} : n \in \omega\}$. By construction, $G_f$ is $\bbO_{S,x}^{L[S,x,Z]}$-generic over $\HOD^{L[S,x,Z]}_{S,x}$. Also by construction, if $f \neq g$, then $\dot x_\mathrm{gen}[G_f] \neq \dot x_\mathrm{gen}[G_g]$.

Hence $A = \{\dot x_\mathrm{gen}[G_f] : f \in \cantorspace\}$ is an uncountable set of reals.

Subclaim 1.1: $A \subseteq R_x$.

To see Subclaim 1.1: Note that for all $f \in \cantorspace$, $p \in G_f$. Note that the condition $p$ takes the form specified in Fact \ref{restate vopenka theorem}. Hence the fact asserts that
$$\HOD_{S,x}^{L[S,x,Z]} \models p \forces_{\bbO_{S,x}^{L[S,x,Z]}} L[S,x,\dot x_\mathrm{gen}] \models \varphi(\check S,\check x,\dot x_\mathrm{gen}).$$
Therefore,
$$\HOD_{S,x}^{L[S,x,Z]}[G_f] \models L[S,x,\dot x_\mathrm{gen}[G_f]] \models \varphi(S,x,\dot x_\mathrm{gen}[G_f]).$$
So in particular, 
$$V \models L[S,x,\dot x_\mathrm{gen}[G_f]] \models \varphi(S,x,\dot x_\mathrm{gen}[G_f]).$$
Since $R = \mathfrak{B}^2_{(S,\varphi)}$, $R(x,\dot x_\mathrm{gen}[G_f])$ for each $f \in \cantorspace$. This shows that $A \subseteq R_x$ which completes the proof of the subclaim.

However, $A$ is uncountable and $R_x$ was assumed to be countable. Contradiction. This completes the proof of Claim 1.

Let $D \subseteq \bbO_{S,x}^{L[S,x,Z]}$ be the dense below $p$ set that witnesses Claim 1. Take any $q \in D$ with $q \leq p$. Define $t \in \cantorspace$ by $t(n) = i \Leftrightarrow q \forces \dot x_\mathrm{gen}(\check n) = \check i)$. It is clear that $t \in \HOD_{S,x}^{L[S,x,Z]}$. Again since $\mathscr{P}(\bbO_{S,x}^{L[S,x,Z]})^{\HOD_{S,x}^{L[S,x,Z]}}$ is countable in $V$, let $G$ be any $\bbO_{S,x}^{L[S,x,Z]}$-generic filter over $\HOD^{L[S,x,Z]}_{S,x}$ with $q \in G$ (and hence $p \in G$). By the forcing theorem, $\dot x_\mathrm{gen}[G] = t$. Using Fact \ref{restate vopenka theorem} as above, one has $\HOD_{S,x}^{L[S,x,Z]}[G]\models L[S,x,t] \models \varphi(S,x,t)$. Hence $V \models L[S,x,t] \models \varphi(S,x,t)$. Thus $R(x,t)$. So $t \in R_x \cap \HOD_{S,x}^{L[S,x,Z]}$. 

It has been shown that $R_x \cap \HOD_{S,x}^{L[S,x,Z]} \neq \emptyset$ for any $Z \in \degrees$ with $Z \geq_T [y^*]_T$. Lemma \ref{selection lemma for HOD} implies that $R$ has a uniformization.
\end{proof}

In $\ZF$: The Silver's dichotomy \cite{Counting-the-Number-of-Equivalence-Classes} states that for every $\coanalytic$ equivalence $E$ on $\reals$, $\reals \slash E$ injects into $\omega$ or $\reals$ injects into $\reals \slash E$. Burgess \cite{Effective-Enumeration-of-Classes-in-a} showed that if $E$ is a $\analytic$ equivalence relation, $\reals \slash E$ injects into $\omega$, $\reals \slash E$ is in bijection with $\omega_1$, or $\reals$ inject into $\reals \slash E$. Note that both results imply that either $\reals \slash E$ is wellorderable or $\reals$ injects into $\reals \slash E$. 

Woodin's perfect set dichotomy is an extension of this property to all equivalence relations on $\reals$ assuming $\AD$, all sets of reals have $\infty$-Borel codes, and the ultrapower, $\prod_{X \in D} \omega_1 \slash \mu_\degrees$, is wellfounded. The wellfoundedness of the ultrapower is certainly a consequence of $\DC$.

The proof of Woodin's perfect set dichotomy presented below follows the outline of Hjorth's \cite{Dichotomy-for-Definable-Universe} generalization of the \cite{Glimm-Effros-Dichotomy-for-Borel-Equivalence-Relations} $E_0$-dichotomy in $L(\reals) \models \AD$ and Harrington's proof of the Silver's dichotomy. Harrington's proof uses the Gandy-Harrington forcing of nonempty $\Sigma_1^1$ definable subsets of $\reals$. The Vop\'enka forcing is simply the $\OD$ version of the Gandy-Harrington forcing. 

The argument presented below appears in \cite{Cardinality-Wellordered-Disjoint-Unions-Quotients-Smooth} where the uniformity of this proof is needed to make further conclusions about wellordered disjoint unions of quotients of ``smooth'' equivalence relations with countable classes.

\Begin{theorem}{woodin perfect set dichotomy}
(Woodin's perfect set dichotomy) Assume $\ZF + \AD + \mathsf{V = L(\reals)}$. Let $E$ be an equivalence relation on $\reals$. Then either

\indent (i) $\reals \slash E$ is wellorderable.

\indent (ii) $\reals$ inject into $\reals \slash E$.
\end{theorem}

\begin{proof}
Let $E$ be an equivalence relation on $\reals$. An $E$-component is a nonempty set $A$ so that for all $x,y \in A$, $x \ E \ y$. That is, an $E$-component is simply a nonempty subset of an $E$-class.

By Corollary \ref{infity borel codes in L(R)}, every set of reals has an $\infty$-Borel code. Let $(S,\varphi)$ be an $\infty$-Borel code for $E$; that is, $E = \mathfrak{B}^2_{(S,\varphi)}$. Throughout this argument, $E$ will always be considered as the set defined by the $\infty$-Borel code $(S,\varphi)$.

(Case I) For all $X \in \degrees$, for all $a \in \reals^{L[S,X]}$, there is some $\OD_S^{L[S,X]}$ $E$-component containing $a$.

In other words, for all degrees $X \in \degrees$, in the local model $L[S,X]$, every real belongs to an $\OD_S$ $E$-component.

For each $\alpha \in \prod_{X \in \degrees} \omega_1 \slash \mu_\degrees$, let $f : \degrees \rightarrow \omega_1$ be such that $f$ is a representative for $\alpha$. Define $A_\alpha \subseteq \reals$ as follows: $a \in A_\alpha$ if and only if for $\mu_\degrees$-almost all $X$, $a$ belongs to the $f(X)^\text{th}$ $\OD_S^{L[S,X]}$ $E$-component according to the canonical wellordering of $\HOD_S^{L[S,X]}$. (If there is no $F(X)^\text{th}$ $\OD_S^{L[S,X]}$ $E$-component, then one just let this set be $\emptyset$.) One can check that $A_\alpha$ is well defined, independent of the choice of $f$ representing $\alpha$.

$A_\alpha$ is an $E$-component: To see this, let $a,b \in A_\alpha$. Thus there is a set $K \in \mu_\degrees$ so that for all $X \in K$, $[a\oplus b]_T \leq X$ and $a$ and $b$ both belong to the $f(X)^\text{th}$ $\OD_S$ $E$-component of $\HOD_{S}^{L[S,X]}$. Thus, $L[S,X] \models a \ E \ b$. Since $E$ is always defined by it $\infty$-Borel code, one has, $L[S,X] \models L[S,a,b] \models \varphi(S,a,b)$. Thus $V \models L[S,a,b] \models \varphi(S,a,b)$. Therefore, in $V$, $a \ E \ b$. It has been shown that $A_\alpha$ is an $E$-component. 

For every $a \in \reals$, there is some $\alpha \in \prod_{X \in \degrees}\omega_1 \slash \mu_\degrees$ such that $a \in A_\alpha$: To see this. Define $f: \degrees \rightarrow \omega_1$ by letting $f(X)$ be the least $\beta$ so that $a$ belongs to the $\beta^\text{th}$ $\OD_S^{L[S,X]}$ $E$-component of $L[S,X]$ according to the canonical wellordering of $\HOD_{S}^{L[S,X]}$. This $\beta$ exists due to the Case I assumption. Let $\alpha = [f]_{\mu_\degrees}$. Then $a \in A_\alpha$. 

Since $L(\reals) \models \DC$, $\prod_{X \in \degrees} \omega_1 \slash \mu_\degrees$ is a wellordered set. Thus $\langle A_\alpha : \alpha \in \prod_{X \in \degrees} \omega_1 \slash \mu_\degrees\rangle$ is a wellordered sequence of $E$-components with the property that every real belongs to some $A_\alpha$. One can now wellorder $\reals \slash E$ as follows: For two $E$-classes $u,v \in \reals \slash E$, $u \prec v$ if and only if the least $\alpha$ so that $A_\alpha \subseteq u$ is less than the least $\alpha$ so that $A_\alpha \subseteq v$. Thus $\prec$ wellorders $\reals \slash E$. 

(Case II) There exists an $X \in \degrees$ and an $a \in \reals^{L[S,X]}$ so that there is no $\OD_S^{L[S,X]}$ $E$-component containing $a$. 

In other word, there is a particular local model $L[S,X]$ so that within this model, there is a real $a$ which does not belong to any $\OD_S$ $E$-component.

Fix such a degree $X \in \degrees$. One will always work in this local model $L[S,X]$. 

In $L[S,X]$, define $u$ be the collection of reals of $L[S,X]$ that do not belong to any $\OD_S^{L[S,X]}$ $E$-component. Note that $u \neq \emptyset$ by the Case II assumption and $u$ is $\OD_S^{L[S,X]}$. Thus $u$ is condition of $\bbO_S^{L[S,X]}$. 

Let $\dot x_\mathrm{left}$ and $\dot x_\mathrm{right}$ be $\bbO_S^{L[S,X]} \times \bbO_S^{L[S,X]}$-names so that $\dot x_\mathrm{left}$ and $\dot x_\mathrm{right}$ is the evaluation of the $\bbO_S^{L[S,X]}$-name $\dot x_\mathrm{gen}$ according to the left and right $\bbO_S^{L[S,X]}$-generic filter, respectively, coming from an $\bbO_S^{L[S,X]} \times \bbO_S^{L[S,X]}$-generic filter.

Claim 1: $\HOD_{S}^{L[S,X]} \models (u,u) \forces_{\bbO_S^{L[S,X]}\times\bbO_S^{L[S,X]}} \neg (\dot x_\mathrm{left} \ E \ \dot x_\mathrm{right})$.

To see Claim 1: Suppose not. Then there is some $(v,w) \leq_{\bbO^{L[S,X]}_S \times \bbO_S^{L[S,X]}} (u,u)$ so that $\HOD_S^{L[S,X]} \models (v,w) \forces_{\bbO_S^{L[S,X]}\times\bbO_S^{L[S,X]}} \dot x_\mathrm{left} \ E \ \dot x_\mathrm{right}$.

Subclaim 1.1: Suppose $G_0$ and $G_1$ are $\bbO_S^{L[S,X]}$-generic filter over $\HOD_S^{L[S,X]}$ which belong to $V$ and contain the condition $v$, then $\dot x_\mathrm{gen}[G_0] \ E \ \dot x_\mathrm{gen}[G_1]$. (Note that $G_0$ and $G_1$ are not necessarily mutually generic.)

To see Subclaim 1.1: As before, since $\HOD_S^{L[S,X]} \models \AC$, $\mathscr{P}(\bbO_S^{L[S,X]})^{\HOD^{L[S,X]}_S}$ is countable in the real world which is a model of $\AD$. Therefore, one can find in the real world $V$, an $H \subseteq \bbO_S^{L[S,X]}$ which is $\bbO_S^{L[S,X]}$-generic over both $\HOD_S^{L[S,X]}[G_0]$ and $\HOD_S^{L[S,X]}[G_1]$ and $w \in H$. Applying the forcing theorem, one has that $\HOD_{S}^{L[S,X]}[G_0][H] \models \dot x_\mathrm{gen}[G_0] \ E \ \dot x_\mathrm{gen}[H]$ and $\HOD_S^{L[S,X]}[G_1][H] \models \dot x_\mathrm{gen}[G_1] \ E \ \dot x_\mathrm{gen}[H]$. As $E$ is defined by its $\infty$-Borel code $(S,\varphi)$, one has
$$\HOD_S^{L[S,X]}[G_0][H] \models L[S,\dot x_\mathrm{gen}[G_0],\dot x_\mathrm{gen}[H]] \models \varphi(S,\dot x_\mathrm{gen}[G_0], \dot x_\mathrm{gen}[H])$$
$$\HOD_S^{L[S,X]}[G_1][H] \models L[S,\dot x_\mathrm{gen}[G_1],\dot x_\mathrm{gen}[H]] \models \varphi(S,\dot x_\mathrm{gen}[G_1], \dot x_\mathrm{gen}[H])$$
In particular
$$V\models L[S,\dot x_\mathrm{gen}[G_0],\dot x_\mathrm{gen}[H]] \models \varphi(S,\dot x_\mathrm{gen}[G_0], \dot x_\mathrm{gen}[H])$$
$$V\models L[S,\dot x_\mathrm{gen}[G_1],\dot x_\mathrm{gen}[H]] \models \varphi(S,\dot x_\mathrm{gen}[G_1], \dot x_\mathrm{gen}[H])$$
Thus in the real world, $\dot x_\mathrm{gen}[G_0] \ E \ \dot x_\mathrm{gen}[H]$ and $\dot x_\mathrm{gen}[G_1] \ E \ \dot x_\mathrm{gen}[H]$. By the transitivity of the equivalence relation $E$, $\dot x_\mathrm{gen}[G_0] \ E \ \dot x_\mathrm{gen}[G_1]$. This proves Subclaim 1.1.

Observe that there must be some $a,b \in v$ so that $\neg(a \ E \ b)$ since otherwise $v$ would be an $\OD_S^{L[S,X]}$ $E$-component containing $a$. This is impossible since $v \leq_{\bbO^{L[S,X]}_S} u$ implies that $v \subseteq u$ and $u$ consists of those $a \in \reals$ which do not belong to any $\OD_S^{L[S,X]}$ $E$-component. Thus $p = v \times v \setminus E = \{(a,b) \in v \times v : \neg(a \ E \ b)\}$ is a nonempty $\OD_S^{L[S,X]}$ subset of $\reals^2$. Thus $p \in {}_2\bbO_S^{L[S,X]}$.

Let $\dot x_\mathrm{gen}^2$ denote the generic element of $\reals^2$ added by ${}_2\bbO_S^{L[S,X]}$. Let $\tau_0$ and $\tau_1$ be the ${}_2\bbO_S^{L[S,X]}$-names for the first and second coordinate of $\dot x_\mathrm{gen}^2$. Using the $\infty$-Borel code $(S,\varphi)$, the condition $q = \{(a,b) : \neg(a \ E \ b)\} = \{(a,b) : L[S,a,b] \models \neg \varphi(S,a,b)\}$ takes the form specified in Fact \ref{restate vopenka theorem}. Therefore, Fact \ref{restate vopenka theorem} implies that $q \forces_{{}_2\bbO_S^{L[S,X]}} L[S,\tau_0,\tau_1] \models \neg \varphi(S,\tau_0,\tau_1)$. 

As before, in the real world, find some $G \subseteq {}_2\bbO_S^{L[S,X]}$ containing $p$ which is ${}_2\bbO_S^{L[S,X]}$-generic over $\HOD_S^{L[S,X]}$. Since $p \leq_{{}_2\bbO_S^{L[S,X]}} q$, one has that $q \in G$. Hence 
$$\HOD_{S}^{L[S,X]}[G] \models L[S,\tau_0[G],\tau_1[G]] \models \neg\varphi(S,\tau_0[G],\tau_1[G]).$$
In particular, 
$$V \models L[S,\tau_0[G],\tau_1[G]] \models \neg\varphi(S,\tau_0[G],\tau_1[G]).$$ 
Since $(S,\varphi)$ is the $\infty$-Borel code for $E$, one has that $\neg(\tau_0[G] \ E \ \tau_1[G])$. However, $\tau_0[G]$ and $\tau_1[G]$ are just the two coordinates of $\dot x^2_\mathrm{gen}[G]$, which is the generic element of $\reals^2$ added by $G$. Fact \ref{individual vopenka generic} implies that $\tau_0[G]$ and $\tau_1[G]$ are individually $\bbO_S^{L[S,X]}$-generic filters over $\HOD_S^{L[S,X]}$. Then Subclaim 1.1 implies that $\tau_0[G] \ E \ \tau_1[G]$. Contradiction. This shows Claim 1.

As before, $\mathscr{P}(\bbO_S^{L[S,X]} \times \bbO_S^{L[S,X]})^{\HOD_S^{L[S,X]}}$ is countable in the real world. Let $\langle D_n : n \in \omega\rangle$ enumerate all the dense open subsets of $\bbO_S^{L[S,X]} \times \bbO_S^{L[S,X]}$ which belong to $\HOD_S^{L[S,X]}$. By intersecting, one may assume that $D_{n + 1} \subseteq D_n$ for all $n \in \omega$. Now in the usual manner, one will build a perfect tree of conditions whose distinct branches correspond to mutually generic filters. The details follow:

Let $p_\emptyset$ be any condition below $u$. Suppose for some $n \in \omega$, $p_\sigma$ has been defined for all $\sigma \in {}^n 2$. For each $\sigma \in {}^n 2$, $(p_\sigma,p_\sigma) \leq_{\bbO_S^{L[S,X]}\times\bbO_S^{L[S,X]}} (u,u)$ and Claim 1 implies that there must be some $q_0 \leq_{\bbO_S^{L[S,X]}} p_\sigma$ and $q_1 \leq_{\bbO_S^{L[S,X]}} p_\sigma$ which are incompatible conditions. Let $(\rho_{\sigma \hat{\ }0},\rho_{\sigma\hat{\ }1})$ denote the $\HOD_S^{L[S,X]}$-least such pair $(q_0,q_1)$.

For $\sigma \in {}^{n + 1} 2$, let $\rho_\sigma^{-1} = \rho_\sigma$. Let $\langle \sigma_i : i < 2^{n + 1}\rangle$ enumerate ${}^{n + 1}2$. 

Now suppose that $\rho_\sigma^j$ has been defined for all $\sigma \in {}^{n + 1}2$ and $-1 \leq j < 2^{n + 1}$. By considering all possible pairs and extending by density, find a collection $\{\rho_{\sigma_i}^{j + 1} : 0 \leq i < 2^{n + 1}\}$ with the property that  

\noindent (1) $\rho_\sigma^{j + 1} \leq_{\bbO_s^{L[S,X]}} \rho_\sigma^j$ 

\noindent (2) $(\rho_{\sigma_{j + 1}}^{j + 1},\rho_{\sigma_{\ell}}^{j + 1}) \in D_n$ and $(\rho_{\sigma_{\ell}}^{j + 1},\rho_{\sigma_{j + 1}}^{j + 1}) \in D_n$ whenever $\ell \neq j + 1$. 

For $\sigma \in {}^{n + 1}2$, let $p_\sigma = \rho_\sigma^{2^{n + 1} - 1}$. 

This defines a sequence $\langle p_\sigma : \sigma \in \finBinarySequence\rangle$. For each $f \in \cantorspace$, let $G_f$ be the $\leq_{\bbO_S^{L[S,X]}}$-upward closure of $\{p_{f \upharpoonright n} : n \in \omega\}$. By construction, if $f \neq g$, then $G_f \times G_g$ is an $\bbO_S^{L[S,X]} \times \bbO_S^{L[S,X]}$-generic filter over $\HOD_S^{L[S,X]}$ containing the condition $(u,u)$. Therefore, by Claim 1, for any $f,g \in \cantorspace$ with $f \neq g$, $\neg(\dot x_\mathrm{gen}[G_f] \ E \ \dot x_{\mathrm{gen}}[G_g])$. Thus $A = \{\dot x_\mathrm{gen}[G_f] : f \in \cantorspace\}$ is a perfect set of $E$-inequivalent reals. Thus $\Phi : \reals \rightarrow \reals \slash E$ defined by $\Phi(f) = [\dot x_\mathrm{gen}[G_f]]_E$ is an injection.
\end{proof}

\Begin{corollary}{perfect set result for surjective image of R}
Assume $\ZF + \AD + \mathsf{V = L(\reals)}$. If $X$ is a surjective image of $\reals$ (equivalently $X \in L_\Theta(\reals)$), then $X$ is either wellorderable or $\reals$ injects into $X$. 
\end{corollary}

In fact, as observed by \cite{A-Trichotomy-Theorem-in-Natural}, in $L(\reals)$, this dichotomy actually holds for all sets $X$ in $L(\reals)$ not just $X \in L_\Theta(\reals)$.

\bibliographystyle{amsplain}
\bibliography{references}

\end{document}